\documentclass{article}
\usepackage{graphicx} 
\usepackage{amsmath} 
\usepackage{amsthm}  
\usepackage{amssymb} 
\usepackage{mathrsfs} 
\usepackage{esint}
\newtheorem{theorem}{Theorem}[section] 
\newtheorem{proposition}[theorem]{Proposition} 
\newtheorem{lemma}[theorem]{Lemma} 
\newtheorem{corollary}[theorem]{Corollary} 
\newtheorem{definition}[theorem]{Definition} 
\newtheorem{remark}{Remark}[section] 
\usepackage{xcolor}  
\usepackage{hyperref}
\hypersetup{colorlinks=true,linkcolor=blue,anchorcolor=blue,citecolor=blue}
\title{Logarithmic Potential and its Applications in Conformal Geometry}
\author{Liming Liu, Shiguang Ma}
\date{\today}
\numberwithin{equation}{section}
\begin{document}
\maketitle 
\begin{abstract}
    We establish an intrinsic characterization of the end invariants $\nu_i$ in the Chern--Gauss--Bonnet formula for complete, simply connected, locally conformally flat 4-manifolds with non-negative scalar curvature and $Q$-curvature non-negative outside a compact set. For each end $M_i$, we show that $\nu_i$ is the cube root of its normalized asymptotic geodesic-ball volume growth:
    \[
        \nu_i=\left(\lim_{r\to\infty}
        \frac{V_g(B_r(p)\cap M_i)}{(\pi^2/2)r^4}\right)^{1/3}.
    \]
    This characterization uses only the Riemannian distance and volume of $g$; it is independent of the base point $p$, the choice of compact truncation of the end, and any conformal coordinates. The proof is based on asymptotic radial symmetry of the logarithmic potential outside a sparse $\mathscr{E}$-set, inspired by Hayman's work \cite{hayman1960slowly}, together with comparisons of radial lengths and geodesic distances.
\end{abstract}
\section{Introduction and Background}
\subsection{Background}
We begin with a brief review of the relevant background.\par
The well-known Gauss--Bonnet formula on a closed manifold $(M^2,g)$ is
\begin{align}
    \int_{M}K[g]dV[g]=2\pi\chi(M),\label{WellKnown:GaussBonnet}
\end{align}
where $K[g]$ is Gaussian curvature and $\chi(M)$ is the Euler number of $M$.
Later, Huber\cite{Huber1958subharmonic} extended this result and showed that if the negative part of the Gaussian curvature $K[g]$ is integrable, such a complete  surface can be compactified conformally by adding finitely many points. A sufficiently small punctured neighborhood of each added point can be conformally regarded as $\mathbb{R}^2\backslash B$.\par
In the 1960s, Finn \cite{Finn1965class} generalized (\ref{WellKnown:GaussBonnet}) based on the result of Huber \cite{Huber1958subharmonic}. Finn showed that if $(M^2,g=e^{2u}|dx|^2)$ is complete and $K[g]$ is absolutely integrable, then 
\begin{equation}
    \chi(M)-\frac{1}{2\pi}\int_{M}K[g]dV[g]=\sum_{i=1}^k\mu_i,
\end{equation}
where $$\mu_i=\lim_{r\to\infty}\frac{\mathscr{L}^2_i(r)}{4\pi \mathscr{A}_i(r)}$$
is the geometric quantity that describes the shape  of the i-th end $E_i:=(\mathbb{R}^2\backslash B,g_i=e^{2u_i(x)}|dx|^2)$ near infinity, and $\mathscr{A}_{i}(r)$ and $\mathscr{L}_i(r)$ are the area of the annulus $B^e_r\backslash B$ and the length of $ \partial B^e_{r}$ with respect to the metric $g_i$, respectively.\par
In recent years, $Q$-curvature has come to play a role in the study of high-dimensional conformal geometry,  similar to that of Gaussian curvature in two-dimensional manifolds.
In the conformal geometry of dimension 4, Paneitz \cite{PaneitzOperator} introduced the Paneitz operator 
\begin{equation}
    P_4[g]=\Delta^2_g-div_g\left[\left(\frac{2}{3}R[g]\cdot g-2Ric[g]\right)d\right],\label{Def:PaneitzOperator}
\end{equation}
on $(M^4,g)$, where $Ric[g]$ and $R[g]$ are the Ricci curvature tensor and scalar curvature of $(M^4,g)$, respectively. Later, Branson \cite{Branson1985QCurvatureInduce} introduced the associated $Q$-curvature
\begin{equation}
    Q[g]=\frac{1}{12}\left\{R[g]^2-3|Ric[g]|_g^2-\Delta_gR[g]\right\},\label{Def:QCurvature}
\end{equation}
where $|\cdot|_g$ denotes the norm with respect to the metric $g$, and $\Delta_g$ denotes the Laplace--Beltrami operator of $g$. We write $\Delta_0$ for the standard Euclidean Laplacian and may omit $g$ when it is the standard Euclidean metric. \par
It is well known that the conformal transformation of Gaussian curvature satisfies the following equation
$$\Delta _{g_0}w+K_{g_w}e^{2w}=K_{g_0}.$$
Moreover, for the conformal metric $g=e^{2u}\bar g$, the Paneitz operator $P[g]$ and the associated $Q$-curvature satisfy
\begin{equation}
    P_4[\bar g]u+2Q[\bar g]=2Q[g]e^{4u}.
\end{equation}
Furthermore, when the manifold has a boundary, Chang and Qing \cite{ChangQing1997BoundaryOperator} have defined the boundary operator $P_3$ together with the associated boundary curvature invariant $T$:
\begin{align}
    P_3[\bar g]u+T \left[\bar g\right]=T[g]e^{3u},\label{BoundaryOpt}
\end{align}
where 
\begin{align}
    P_3 = -\frac{1}{2} \partial_n \Delta
      - \tilde{\Delta} \partial_n
      - \frac{2}{3} H \tilde{\Delta}
      + L_{\alpha\beta} \tilde{\nabla}_\alpha \tilde{\nabla}_\beta
      + (\frac{1}{3}R - R_{\alpha N \alpha N}) \partial_n
      + \frac{1}{3} \tilde{\nabla} H \cdot \tilde{\nabla}.\label{DefOfP_3}
\end{align}
For an oriented hypersurface $\Sigma$ with a chosen unit normal $N$,
we use the convention
\[
L_N(X,Y):=g(\nabla_XN,Y),\qquad
H_N:=\operatorname{tr}_{g|_\Sigma}L_N.
\]
Thus the Euclidean unit sphere with its outward normal has
$L_N=g_{S^3}$ and $H_N=3$.
We define
\begin{align}
T_N[g]
={}&\frac{1}{12}\partial_N R[g]
+\frac{1}{6}R[g]H_N
-R_{aNbN}(L_N)_{ab}
+\frac{1}{9}H_N^3 \nonumber\\
&-\frac{1}{3}\operatorname{tr}(L_N^3)
-\frac{1}{3}\widetilde{\Delta}H_N.
\label{DefOfBoundaryCurvatureT}
\end{align}
Here $\widetilde{\Delta}$ is the Laplace--Beltrami operator
of the induced metric on $\Sigma$.
For the boundary of a compact domain, $N$ is taken to be
the outward unit normal unless otherwise specified. For a conformally flat end
$E=(\mathbb{R}^4\setminus B,g=e^{2u}|dx|^2)$,
write $r=|x|$. On each coordinate sphere, we set
\[
N_+[g]:=e^{-u}\partial_r,
\qquad
T_+[g]:=T_{N_+[g]}[g].
\]
Thus $N_+[g]$ points toward infinity along the end.
Reversing the unit normal changes the sign of the boundary
curvature:
\[
T_{-N}[g]=-T_N[g].
\]
At the interface between the end and a compact core $K$,
$N_+[g]$ is the outward unit normal of $K$, whereas the
outward unit normal of the end is $-N_+[g]$.
For further details on (\ref{Def:PaneitzOperator}), (\ref{Def:QCurvature}), and (\ref{BoundaryOpt}), readers may refer to \cite{ChangQing1997BoundaryOperator,Chang2004ComformalGeometryBook}.\par
Then Chang-Qing-Yang\cite{ChangQingYangDuke2000} made a successful generalization of Finn's work \cite{Finn1965class}. On the complete manifold  $(\mathbb{R}^4,g=e^{2u}|dx|^2)$, they showed that if $R[g]\geq0$ near infinity and $Q[g]\in L^1(\mathbb{R}^4, g)$, then 
$$
    \lim_{l\to\infty}\frac{A_g(\partial B^e_l(0))^{4/3}}{4(2\pi^2)^{1/3}V_g(B^e_l(0))}=1-\frac{1}{4\pi^2}\int_{\mathbb{R}^4}Q[g]dV[g]\geq0.
$$ 
Here, $A_g(\cdot)$ and $V_g(\cdot)$ denote the three-dimensional area and four-dimensional volume with respect to the metric $g$, respectively, and $B^e_l(0)$ denotes the four-dimensional Euclidean ball centered at the origin with radius $l$.

In the same year, they extended this result to a more general setting in \cite[Theorem 1.2]{ChangQingYangInvention2000} and generalized the two-dimensional result of Finn \cite{Finn1965class}. By employing methods similar to those in \cite{ChangQingYangDuke2000}, they studied the geometry of ends in the case of manifolds with finitely many ends and showed that the geometric quantity $\nu_i$ for a single end in the Chern–Gauss–Bonnet formula can be expressed as
\begin{align*}
\nu_i=\lim_{l\to\infty}\frac{A_g(\partial B^e_l)^{4/3}}{4(2\pi^2)^{1/3}V_g(B^e_l)}.
\end{align*}
However, such an expression for the geometric quantity depends on the structure of the Euclidean ball $B^e_r$.

In 2024, Zhang \cite{zhang2024sharp} removed the restriction on $R[g]$ and instead imposed the conditions that $Q[g](x)\geq e^{-o_+(1)|x|^2}$ near infinity and $Q_g\in L^1(\mathbb R^n,dV_g)$. He then obtained a result similar to the above result of Chang--Qing--Yang \cite{ChangQingYangDuke2000} on $(\mathbb{R}^n,g=e^{2u}|dx|^2)$, where $n\geq4$ is an even integer.

In 2025, Ma and Qing \cite{MaQing2025LinearPotential} weakened the assumption that $Q[g]$ is absolutely integrable to the assumption that $Q^-[g]$ is absolutely integrable, where $Q^-[g]$ denotes the negative part of $Q[g]$. Moreover, for a closed subset $S$ of a compact manifold $(M^4,\bar g)$, assume that the metric $g=e^{2u}\bar g$ on $D\backslash S$ is geodesically complete near $S$, that $R[g]\geq 0$, and that $$Q^-[g]\in L^1(D\backslash S,g),$$
where $D$ is an open neighborhood of $S$. Then Ma-Qing showed that $S$ consists of only finitely many points. \par

Li \cite{limingxiang2024VolumeEntropy} established some results on the asymptotic behavior of the logarithmic potential in certain settings. Additionally, together with Wei and Xu, Li studied the relationships between the nonnegativity of the scalar curvature, the $Q$-curvature, and the Ricci curvature on manifolds conformal to the entire Euclidean space $(\mathbb{R}^4, g = e^{2u}|dx|^2)$ in \cite{liWeiXu2025geometry}. Their work shows that if the scalar curvature is globally non-negative and the $Q$-curvature is globally non-negative on $(\mathbb{R}^4, g = e^{2u}|dx|^2)$, then the Ricci curvature is globally non-negative  (semi-positive definite).

Currently, the expression for the geometric quantity $\nu_i$ on conformally flat ends still relies on the structure of Euclidean balls and is not derived entirely from intrinsic geometric structures (such as geodesic balls). A natural question arises: can the geometric quantity $\nu_i$ be obtained using purely geometric structures, such as geodesic balls?

This question is precisely addressed by the results of this paper. In brief, inspired by the works of Chang–Qing–Yang \cite{ChangQingYangInvention2000} and Hayman \cite{hayman1960slowly, hayman1952minimum}, we discover that the logarithmic potential can be modified in such a way that it behaves as a radially symmetric function near infinity in a certain sense. Through a series of derivations, we are able to obtain the geometric quantity $\nu_i$ on the manifold ends by taking the cube root of the asymptotic ratio of the volume of geodesic balls to the fourth power of the geodesic radius. In this way, the geometric quantity $\nu_i$ depends solely on the intrinsic geometry of the manifold itself, thereby answering the question raised in the previous paragraph.
\subsection{Main results}
In this paper, by investigating the logarithmic potential and applying its results, we study the geometry of ends $(\mathbb{R}^4\backslash B, g = e^{2u}|dx|^2)$ of complete, simply connected, locally conformally flat manifolds with globally non-negative scalar curvature and $Q$-curvature non-negative outside a compact set, as well as the Chern–Gauss–Bonnet formula on such manifolds. The main results are as follows:
\begin{theorem}[Intrinsic characterization of the end invariants]\label{FiniteVolumeE11113}
    Suppose that $(M^4,g)$ is a complete, simply connected, locally conformally flat manifold, the scalar curvature $R[g]$ is non-negative, and the $Q$-curvature $Q[g]$ is non-negative near infinity. Then, we have the Chern-Gauss-Bonnet formula
        \begin{align}
            \chi(M)-\frac{1}{32\pi^2}\int_{M}8Q[g]dV[g]=\sum_{i=1}^m\nu_i,\label{chernGaussBonnet111}
        \end{align}
    where 
    \begin{align}
            \nu_i=\lim_{r\to\infty} \left(\frac{2}{\pi^2}\right)^{1/3}\frac{V_g(B_{r}(p)\cap M_i)^{1/3}}{r^{4/3}} \label{IthIsoperimetricRatioGeodsc111}
    \end{align}
    gives an intrinsic characterization of the invariant of the $i$-th end $M_i$. Here $p\in M$ is arbitrary, $B_r(p)$ is the geodesic ball of radius $r$ centred at $p$, and $V_g$ is the Riemannian volume. Thus $\nu_i$ is the cube root of the asymptotic volume ratio $V_g(B_r(p)\cap M_i)/((\pi^2/2)r^4)$. The limit is independent of $p$ and of changing $M_i$ within a compact set; in particular, it depends only on the intrinsic geometry of the end and does not require a choice of conformal coordinates.
\end{theorem}
When the fundamental group of the manifold is finite, in the non-simply connected case, the universal covering map ensures that the preimage of any compact set is a finite union of compact sets. Moreover, the $Q$-curvature on the universal covering space also remains non-negative near infinity. Consequently, an analogous conclusion to Theorem \ref{FiniteVolumeE11113} can be extended to the universal covering space of the manifold.
\begin{corollary}
            Suppose that $(M^4,g)$ is a complete,  locally conformally flat manifold, the scalar curvature $R[g]$ is non-negative, and the $Q$-curvature $Q[g]$ is non-negative near infinity. If the fundamental group $\pi_1(M)$ is finite, then the result in Theorem  \ref{FiniteVolumeE11113} is still valid for $\widetilde{M}$.
\end{corollary}

First, we provide an intrinsic characterization of the asymptotic geometry of ends that differs from that in \cite[Theorem 1.2]{ChangQingYangInvention2000}. In \cite[Theorem 1.2]{ChangQingYangInvention2000}, the geometric quantity $\nu_i$ appearing in the Chern–Gauss–Bonnet formula was obtained by Chang–Qing–Yang via the following asymptotic ratio of the metric surface area to the metric volume of the Euclidean ball $B^e_l$:
\begin{align*}
\nu_i=\lim_{l\to\infty}\frac{A_g(\partial B^e_l)^{4/3}}{4(2\pi^2)^{1/3}V_g(B^e_l)},
\end{align*}
where $A_g(\cdot)$ and $V_g(\cdot)$ denote the three-dimensional area and four-dimensional volume with respect to the metric $g$, respectively.

However, such an approach inherently depends on the structure of Euclidean balls. In this paper, we improve this point by replacing it with the asymptotic limit of the ratio between the volume of the intersection of a geodesic ball $B_r(p)$ centered at an arbitrary point $p$ on the manifold $M$ with the end $M_i$ and the fourth power of the geodesic radius $r$, thereby obtaining the geometric quantity $\nu_i$ on the manifold end. More precisely,
\begin{align*}
\nu_i=\lim_{r\to\infty} \left(\frac{2}{\pi^2}\right)^{1/3}\frac{V_g(B_{r}(p)\cap M_i)^{1/3}}{r^{4/3}}.
\end{align*}
This eliminates the dependence on the Euclidean ball structure present in \cite[Theorem 1.2]{ChangQingYangInvention2000}. Consequently, the geometric quantity $\nu_i$ now depends solely on the intrinsic geometry of the manifold itself.

Secondly, in the work of Li–Wei–Xu \cite{liWeiXu2025geometry}, the object of study is the even-dimensional manifold $(\mathbb{R}^n, g = e^{2u}|dx|^2)$. In contrast, this paper focuses on properties of complete, simply connected, locally conformally flat manifolds and the asymptotic behavior of the conformal factor on their conformally flat ends $(\mathbb{R}^4 \backslash B, g = e^{2u}|dx|^2)$, where $B$ denotes the Euclidean unit ball. Unlike the global Euclidean space $\mathbb{R}^4$ considered in \cite{liWeiXu2025geometry}, the background manifold of the end in our setting is the Euclidean space $\mathbb{R}^4$ with the unit ball removed. On even-dimensional geometries $(\mathbb{R}^n, g = e^{2u}|dx|^2)$, the Liouville theorem serves as a cornerstone and a key tool. However, on $(\mathbb{R}^4 \backslash B, g = e^{2u}|dx|^2)$, the Liouville theorem generally does not hold. This constitutes an important distinction between the two settings.

Furthermore, \cite[Theorem 1.2]{liWeiXu2025geometry} shows that on a complete conformally flat manifold $(\mathbb{R}^4, g = e^{2u}|dx|^2)$, if the scalar curvature is non-negative and the $Q$-curvature is non-negative, then the Ricci curvature is also non-negative (semi-positive definite). The assumptions on scalar curvature and $Q$-curvature in that result are global conditions. In contrast, our assumption in this paper is that the $Q$-curvature is non-negative only near infinity, which is a non-global condition. To highlight this distinction, in Section \ref{ChapterExample} we provide an example demonstrating that in the locally conformally flat setting, if the $Q$-curvature and scalar curvature are non-negative near a point, the Ricci tensor at that point may have a negative component and hence fail to be positive semidefinite.

The geometric results above are primarily derived from the study of the uniform asymptotic radial symmetry of the logarithmic potential $\mathcal{L}(f)(x)$ (see the definition in \eqref{defofLF}). More precisely, we obtain the following result for the logarithmic potential. Inspired by the work of Hayman \cite{hayman1960slowly,hayman1952minimum}, we obtain a uniform asymptotic radial symmetry result outside a sparse $\mathscr{E}$-set, which constitutes a refined radialization estimate:

\begin{theorem}\label{line:161HaymanCor222222}
    Let $f\in L^1(\mathbb{R}^4\backslash B)$ be non-negative, and let $\mathcal{L}(f)$ be defined by (\ref{defofLF}). Then the logarithmic potential $\mathcal{L}(f)$ satisfies
    \begin{align}
    \mathcal{L}(f)(x)=G_f(|x|)+o(1),\label{line:165111}
    \end{align}
uniformly for $x\in \mathbb{R}^4\backslash (B\cup E)$ as $|x|\to\infty$, where $E$ is an $\mathscr{E}$-set for which the sum of the areas of the projections of its constituent balls onto $\partial B_1^e$ is finite, and 
    $$G_f(r)=\inf_{|x|=r}\mathcal{L}(f)(x)$$
is the infimum of $\mathcal{L}(f)$ on the sphere $\partial B^e_r$.
\end{theorem}

Inspired by Chang–Qing–Yang \cite{ChangQingYangDuke2000, ChangQingYangInvention2000}, we will show that outside the $\mathscr{E}$-set $E$, the logarithmic potential $\mathcal{L}(f)$ not only converges uniformly to its minimum $G_f(|x|)$ on Euclidean spheres, but also converges uniformly to its spherical average $\mathcal{A}(|x|)$. Consequently, we have
    \begin{equation*}
        G_f(|x|) = \mathcal{A}(|x|) + o(1),~|x|\to\infty,
    \end{equation*}
where $o(1) \to 0$ as $|x| \to \infty$.

\begin{proposition}\label{theorem:111111111}
    Suppose that $\mathcal{L}(f)(x)$ is defined as in (\ref{defofLF}), that $f$ is non-negative, and that $f\in L^1(\mathbb{R}^4\backslash B)$. Then we have
    \begin{equation}
        \mathcal{L}(f)(x) = \mathcal{A}(|x|) + o(1),~|x|\to\infty,
    \end{equation}
    uniformly for $x\in \mathbb{R}^4\backslash (B\cup E)$, where
    \begin{align}
        \mathcal{A}(r):=\fint_{\partial B^e_{r}} {\mathcal{L}(f)}{d\sigma}=\frac{1}{|\partial B^e_r|}\int_{\partial B^e_r}\mathcal{L}(f)(y)d\sigma_y, \label{defofA111}
    \end{align}
    $|\partial B^e_r|:=2\pi^2r^{3}$ is the Euclidean area of $\partial B^e_r$, and $E$ is the $\mathscr{E}$-set mentioned in Theorem  \ref{line:161HaymanCor222222}.
\end{proposition}

\subsection{Organization of the paper}
The structure of this paper is as follows. We investigate complete, simply connected, locally conformally flat manifolds with globally non-negative scalar curvature and $Q$-curvature non-negative outside a compact set. In Section \ref{SectionPreliminary}, inspired by Hayman \cite{hayman1960slowly, hayman1952minimum}, we establish several analytic lemmas as preliminaries, which will be used in the proof of the radial symmetry of the logarithmic potential on conformally flat ends. In Section \ref{SectionAymRadialSymmetry}, we prove the asymptotic radial symmetry of the logarithmic potential and related properties outside a sparse $\mathscr{E}$-set on $\mathbb{R}^4\backslash B$, as stated in Theorem \ref{line:161HaymanCor222222} and Proposition \ref{theorem:111111111}. In Section \ref{ChapterMaintheorem}, we apply these analytic estimates to study the asymptotic geometry of conformally flat ends. Next, building upon the inspiring works \cite{ChangQingYangInvention2000, zhang2024sharp, MaQing2025LinearPotential} and combining the existing compactification results of Ma–Qing \cite{MaQing2025LinearPotential}, Schoen–Yau \cite{SchoenYau1988conformally}, and Chodosh–Li \cite{LiChaoChodosh2024DevelopingMapInjective} together with the results in Sections \ref{SectionAymRadialSymmetry} and \ref{ChapterMaintheorem}, we prove the main result of this paper, Theorem \ref{FiniteVolumeE11113}. In Section \ref{ChapterExample}, we provide an example to highlight the distinction between our work and that of \cite{liWeiXu2025geometry}, showing that in a local setting, the Ricci curvature has a negative component and hence  it may not be positive semidefinite even when the $Q$-curvature and scalar curvature are non-negative. Finally, in the appendix (Section \ref{ChapterAppendix}), we include several proofs to help readers fully grasp the structure of the argument.

\section{Preliminaries}\label{SectionPreliminary}
For convenience, we denote by $|\cdot|$ the Euclidean norm and set $B:=\{x\in\mathbb{R}^4:|x|<1\}$. With an explicit centre, $B_r^e(p):=\{x\in\mathbb{R}^4:|x-p|<r\}$ denotes an ambient open Euclidean ball. Throughout the paper, the notation without a centre denotes the end-truncated ball $B_r^e:=\{x\in\mathbb{R}^4\backslash B:|x|<r\}=B_r^e(0)\backslash B$; it is relatively open in the end and is empty for $r\leq1$. In particular, $B_r^e$ and $B_r^e(0)$ are distinct notations. For every $r>0$, we use $\partial B_r^e:=\{x\in\mathbb{R}^4:|x|=r\}$ for the coordinate sphere, rather than the full topological boundary of the truncated ball; $\partial B_r^e(p)$ denotes the boundary of the ambient ball $B_r^e(p)$. Closures of subsets of the end are taken in $\mathbb{R}^4\backslash B$. We use $V_g(\cdot)$ to denote the volume with respect to a metric $g$, and $d_g(\cdot,\cdot)$ for the corresponding geodesic distance under $g$. In particular, $g_e$ denotes the standard Euclidean metric, $d\sigma$ is the standard area element on $S^3$, and $d\sigma_x:=|x|^3d\sigma$ is the area element on the sphere of radius $|x|$; we also write $d\sigma_r:=r^3d\sigma$. We set $\omega_j:=|S^j|$, so that $\omega_2=4\pi$ and $\omega_3=2\pi^2$. The length of a curve $\gamma:[a,b]\to\mathbb{R}^4\backslash B$ with respect to $g$ is denoted by $L_g(\gamma(t))$. A unit direction in $S^3$ is represented by the vector $\vec n\in S^3$. Where necessary, vector notation is used to distinguish between a numerical value and a point. Finally, $L^1(M,g)$ denotes the space of Lebesgue integrable functions on $M$ with respect to the metric $g$, and $g$ may be omitted when $g=g_e$. Whenever an integral extends beyond the end, a function $f$ defined on $\mathbb{R}^4\backslash B$ is extended by zero to $B$. For such a function, we denote by
\begin{align}
\mathcal{L}(f) := \frac{1}{8\pi^2}\int_{\mathbb{R}^4\backslash B} \log\frac{|y|}{|x-y|} \, f(y) \, dy \label{defofLF}
\end{align}
the logarithmic potential, where $B$ is the open unit ball defined above.

For a nonnegative function $f\in L^1(\mathbb{R}^4\setminus B)$, we initially interpret $\mathcal{L}(f)$ as an extended-real-valued function. Indeed, since $|y|\geq1$,\[\log\frac{|y|}{|x-y|}\geq-\log\left(1+\frac{|x|}{|y|}\right)\geq-\log(1+|x|).\] Consequently,\[\mathcal{L}(f)(x)\geq-\frac{\|f\|_{L^1}}{8\pi^2}\log(1+|x|),\] so the defining integral is well defined in $(-\infty,+\infty]$, although it may equal $+\infty$ at some points. For every fixed $r>1$, the logarithmic singularity is integrable on the sphere, and\[\sup_{|y|\geq1}\fint_{\partial B_r^e}\left|\log\frac{|y|}{|x-y|}\right|\,d\sigma_x <\infty.\] To see this, one may first restrict $y$ to a bounded region and use local coordinates on the sphere: the singular contribution is bounded by a constant multiple of $\int_0^\delta \rho^2|\log\rho|\,d\rho<\infty$, uniformly in $y$. For sufficiently large $|y|$, the kernel is uniformly bounded because $|x|=r$ is fixed. Thus Tonelli's theorem yields \[\fint_{\partial B_r^e}|\mathcal{L}(f)(x)|\,d\sigma_x<\infty.\] In particular, the spherical average \[\mathcal{A}(r):=\fint_{\partial B_r^e}\mathcal{L}(f)(x)\,d\sigma_x\] and the spherical infimum \[G_f(r):=\inf_{|x|=r}\mathcal{L}(f)(x)\] are finite real numbers for every $r>1$.
 \par

\begin{lemma}[\cite{ChangQingYangInvention2000}.Prop 1.12;\cite{MaQing2025LinearPotential}.Thm 1.5] \label{NormalMetric}
Suppose the conformally flat end $M_i:=(\mathbb{R}^4\backslash B\, ,g=e^{2u}|dx|^2)$ satisfies $R[g]\geq0$ and $Q[g]\geq0$. If $M_i$ is an end of a complete manifold $(M,g)$, then 
  $Q[g]\in L^1(\mathbb{R}^4\backslash B\, ,g)$, and the conformal factor $u(x)$ satisfies
    \begin{align}
                 u(x)=\frac{1}{8\pi^2}\int_{\mathbb{R}^4\backslash B}\log\frac{|y|}{|x-y|}2Q_g(y)e^{4u(y)}dy+\gamma\log|x|+b(x/|x|^2),\label{Newline:115}
    \end{align}
    where $\gamma$ is a constant and $b(x)$ is a biharmonic function on the closed ball $\bar B$. More precisely, $b\in C^\infty(B)\cap C^3(\bar B)$ and $\Delta^2b=0$ in $B$; boundary values and derivatives through order three are understood as the continuous one-sided traces from $B$. A metric whose conformal factor has the form (\ref{Newline:115}) is called a normal metric.
\end{lemma}
\begin{proof}
    First, we show that $f:=2Q[g]e^{4u}\in L^1(\mathbb{R}^4\backslash B)$. Under the inversion $x=y/|y|^2$, the end becomes $(\bar B\backslash\{0\},g=e^{2w(y)}|dy|^2)$, where $w(y)=u(y/|y|^2)-2\log|y|$. The metric is complete near $y=0$ and smooth up to $\partial B$. Since the background metric is Euclidean, the curvature equation reads
    \[
        \Delta^2w=2Q[g]e^{4w}\qquad\text{on }B\backslash\{0\}.
    \]
    The normalization in \cite{MaQing2025LinearPotential} is $Q_4=2Q$ in dimension four. Thus the non-negativity of $Q[g]$ gives the required integrability of $Q_4^-[g]$ near the singularity. By the Radon measure extension argument in the proof of \cite[Theorem 1.5]{MaQing2025LinearPotential},
    \begin{equation}
        \Delta^2w=\widetilde f\label{BiLaplaceRadonMeasure}
    \end{equation}
    is a Radon measure in an open neighborhood $U\subset B$ of $0$, whose restriction to $U\backslash\{0\}$ is $2Q[g]e^{4w}\,dy$. Choose $0<\rho<1$ with $\overline{B^e_\rho(0)}\subset U$. Local finiteness on this compact subset gives
    \[
        \int_{0<|y|\le\rho}2Q[g]e^{4w}\,dy
        \le |\widetilde f|\bigl(\overline{B^e_\rho(0)}\bigr)<\infty,
    \]
    where $|\widetilde f|$ denotes the total variation measure. On the remaining compact annulus $\rho\le|y|\le1$, the metric and curvature are smooth, so the curvature density is integrable there as well. Returning to the original coordinates, we obtain
    \[
        Q[g]\in L^1(\mathbb{R}^4\backslash B,g),
        \qquad f\in L^1(\mathbb{R}^4\backslash B).
    \]

Second, applying \cite[Proposition 1.12]{ChangQingYangInvention2000}, we obtain (\ref{Newline:115}) for $|x|>1$, with $b$ biharmonic in $B$. We also record its boundary regularity. Near the inner boundary of the end, $u$ and $f$ are smooth up to the boundary. The logarithmic kernel and its derivatives through order three are locally integrable in dimension four; after extending the density by zero across the boundary, the corresponding potential has continuous derivatives through order three there. The contribution from large $|y|$ has the same regularity by $f\in L^1$ and differentiation under the integral sign. Consequently, the identity defining the remainder $b$ gives $b\in C^3(\bar B)$. The representation and its derivatives through order three extend to $|x|=1$ by one-sided limits. In particular, the boundary values of $D$ and its derivatives used below agree with those obtained by averaging the smooth conformal factor $u$.
\end{proof}
Regarding the logarithmic potential, certain techniques can be employed to make it behave asymptotically like a radially symmetric function near infinity. Before proceeding, we require several useful lemmas, inspired by the work of Hayman \cite{hayman1960slowly, hayman1952minimum}. The following notation will be used. For a positive measure $\mu$ on $U\subseteq\mathbb{R}^4$, extend $\mu$ by zero outside $U$ and denote 
$$
    n(x,t):=\int_{B^e_t(x)}d\mu(y), 
$$
$$
    N(x,t):=\int_{0}^{t}\frac{n(x,s)}{s}ds=\int_{B^e_t(x)}\log \frac{t}{|y-x|}d\mu(y).
$$
If $x$ is the origin, we omit $x$ and denote them simply as $n(t)$ and $N(t)$. Now, we will carry out some preliminary work. These preparations may appear intricate at first glance, but they will help to clarify and simplify the overall proof.
\begin{lemma}[\cite{hayman1952minimum}.Lemma 4]\label{haymanlem1} Suppose that $n(z_0,h)=n\geq0$, and let $0<d<h/2$. Then there exists a union $S$ of balls whose radii sum to at most $d$ such that, for every $z\in B^e_{{h}/{2}}(z_0)\backslash S$, we have
    \begin{align}
        N(z,h/2)\leq n \log\frac{16h}{d}.\label{line:72}
    \end{align}
\end{lemma}
\begin{proof}
    If $n=0$, then $N(z,h/2)=0$ for every $z\in B^e_{h/2}(z_0)$, so we take $S=\varnothing$. Assume henceforth that $n>0$. For each $k\in\mathbb{Z}^+$, choose a maximal finite family $\{z_{k,1},\ldots,z_{k,v_k}\}$ (possibly empty) satisfying
\begin{align}
    &|z_{k,v} - z_0| < h/2, \quad (v = 1, \cdots, v_k), \label{line:78} \\
    &|z_{k,v} - z_{k,v'}| \geq \frac{d}{4^k}, \quad(v\neq v'), \label{line:79} \\
    &n\left(z_{k,v},\frac{d}{2 \cdot 4^k}\right) > \frac{n}{2^k}. \label{line:80}
\end{align}
Such a maximal family is obtained by adding admissible points until none remain. Indeed, by (\ref{line:79}), the open balls $B^e_{d/(2\cdot4^k)}(z_{k,i})$ are pairwise disjoint, and by (\ref{line:78}) they are contained in $B_h^e(z_0)$. Hence every finite admissible family satisfies
    $$
        n=n(z_0,h)\geq\sum_{i=1}^{v_k}n\left(z_{k,i},\frac{d}{2\cdot4^k}\right)\geq v_k\frac{n}{2^k},
    $$
    so $v_k\leq2^k$ and the selection terminates after finitely many steps. Denote
    $$
        S:=\bigcup_{k=1}^\infty\bigcup_{v=1}^{v_k}B^e_{\frac{d}{4^k}}(z_{k,v}).  
    $$
    It is clear that the sum of the radii of the balls in $S$ is at most
    $$
        \sum_{k=1}^{\infty}\frac{dv_k}{4^k}\leq \sum_{k=1}^{\infty}\frac{d}{2^k}=d.
    $$
    Finally, let $z\in B^e_{h/2}(z_0)\backslash S$. For every $k$, this point is at distance at least $d/4^k$ from each $z_{k,v}$. By maximality, it cannot satisfy (\ref{line:80}), since otherwise it could be added to the family. Therefore
    \begin{equation}
        n(z,\frac{d}{2\cdot4^k})\leq\frac{n}{2^k}.\label{line:96}
    \end{equation}
    Therefore, by (\ref{line:96}), we obtain directly
    \begin{align*}
        N(z,h/2)=&\int_{0}^{h/2}n(z,t)\frac{dt}{t}\\
        =&\left\{\sum_{k=1}^{\infty}\int_{\frac{d}{2\cdot4^{k+1}}}^{\frac{d}{{2\cdot4^{k}}}}+\int_{d/8}^{h/2}\right\}\frac{n(z,t)}{t}dt\\
        \leq&\sum_{k=1}^{\infty}\int_{\frac{d}{2\cdot4^{k+1}}}^{\frac{d}{{2\cdot4^{k}}}}n(z,\frac{d}{2\cdot4^{k}})\frac{dt}{t}+n(z,\frac{h}{2})\int_{d/8}^{h/2}\frac{dt}{t}\\
        \leq&n\log\frac{4h}{d}+\sum_{k=1}^{\infty}\frac{n}{2^k}\log4=n\log\frac{16h}{d}.
    \end{align*}
    This completes the proof.
\end{proof}\par
The following Lemma \ref{haymanLem2} is a corollary of Lemma \ref{haymanlem1} and is repeatedly used in the proof of Lemma \ref{Lemma4}.
\begin{lemma}\label{haymanLem2}
Let $\sigma$ be a non-negative finite measure on $\mathbb{R}^4\backslash B$ that vanishes outside a compact set $K$. Then given $\rho>0$, there exists a collection of balls $S:=\bigcup _{i\in I}B^e_{r_i}(\cdot)$ such that for every $z \in (\mathbb{R}^4\backslash B) \backslash S$, we have
\begin{equation}
\int_{\mathbb{R}^4\backslash B} \log|z-y| d\sigma(y) \geq n \log \rho , \label{line:109}
\end{equation}
where $n := \sigma(\mathbb{R}^4\backslash B)$, and the radii of the balls in $S$ satisfy $\sum_{i\in I} r_i \le  32\rho$.
\end{lemma}
\begin{proof}
    Fix $\rho>0$. Choose $R$ sufficiently large that $K\subset B_R^e$ and $R>15\rho$. Thus $\sigma$ vanishes outside $B_R^e$.\par
    On the one hand,  for $|z| \geq R+\rho$, inequality (\ref{line:109}) follows directly from the following computation:
\begin{align*}
            \int_{\mathbb{R}^4\backslash B}\log|z-y|d\sigma(y)\geq&\int_{\mathbb{R}^4\backslash B}\log(|z|-|y|)d\sigma(y)\geq\int_{\mathbb{R}^4\backslash B}\log(|z|-R)d\sigma(y)\nonumber\\
            =&\int_{B^e_R}\log(|z|-R)d\sigma(y)\nonumber=n\log(|z|-R)\geq n\log\rho.\label{line:117}
\end{align*}\par
On the other hand, we consider those $z$ such that $|z|\leq R+\rho$. Applying Lemma \ref{haymanlem1} with $h:=4(R+\rho)$ and $d:=32\rho<\frac{h}{2}$, we find that, for $z\in B^e_{h/4}\backslash S$, we have
\begin{equation}
    N(z,\frac{h}{2})=\int_{B^e_{h/2}(z)}\log\frac{h}{2}+\log\frac{1}{|z-y|}d\sigma(y)\leq  n\log\frac{16h}{d}, \label{line:121}
\end{equation}
where $S:=\bigcup_{i\in I}B^e_{r_i}(\cdot)$ satisfies $\sum_{i\in I}r_i\leq d$.
Using
$$\int_{B^e_R(0)}\log\frac{h}{2}d\sigma(y)\leq n\log\frac{h}{2}, $$
together with (\ref{line:121}), we obtain
$$\int_{\mathbb{R}^4\backslash B}\log|z-y|d\sigma(y)\geq n\log\rho.$$
In summary, we conclude that (\ref{line:109}) holds for all $z \in (\mathbb{R}^4\backslash B) \backslash S$, thereby completing the proof.
\end{proof}
To simplify the notation, for $0<r_1<r_2$ we set $A_{r_1}^{r_2}:=B_{r_2}^e\backslash B_{r_1}^e=\{x\in\mathbb{R}^4\backslash B:r_1\leq|x|<r_2\}$. Thus these annuli use the same end-truncated convention, including when $r_1<1$.
\begin{lemma}\label{Lemma4}
Let $\mu$ be a positive finite measure on $\mathbb{R}^4\backslash B$ with total mass $\mu(\mathbb{R}^4\backslash B) = n < \infty$. Then, given any $k \geq 7$, there exists a set of balls $E_k$ depending on $k$ such that for every $z\in(\mathbb{R}^4\backslash B)\backslash E_k$ with $t:=|z|$, we have
\begin{align}
\int_{\frac{t}{2} < |y| < 2t} \log\frac{|z-y|}{t + |y|}  d\mu(y) \geq -nk, \label{line:189A}
\end{align}
where the projection of $E_k$ on $\partial B_1^e$ has area at most $M_0 e^{-k}$. Here $M_0$ is a constant independent of $k$.
\end{lemma}
\begin{proof}
    If $n$ vanishes, the Lemma holds trivially, since the integral in (\ref{line:189A}) vanishes for sufficiently large $t$. In what follows, we may assume without loss of generality that $n>0$.
    
    The outline of the proof is as follows. We partition $\mathbb{R}^4\backslash B = \bigcup_{v=0}^{+\infty} A_{2^v}^{2^{v+1}}$. For a fixed $k \geq 7$ and for each annulus $A_{2^v}^{2^{v+1}}$, we apply Lemma \ref{haymanLem2} to extract a special collection of balls $S_{v,k}$ that intersect $A_{2^v}^{2^{v+1}}$. Finally, we demonstrate that the projection of the union $E_k = \bigcup_{v=0}^{\infty} S_{v,k}$ onto $\partial B_1^e$ has area at most $M_0 e^{-k}$.
    
    For $v\in\mathbb{Z}$, denote $\mu_v:=\mu(A^{2^{v+2}}_{2^{v-1}})$, $v=0,1,2,\cdots$. We extend the measure $\mu$ by zero to $B$. Then $\sum_{v=0}^{+\infty}\mu_v\leq3n$. In applying Lemma \ref{haymanLem2}, set  $K:=\overline{A^{2^{v+2}}_{2^{v-1}}}$,  $\sigma:=\mu\cdot\chi(A^{2^{v+2}}_{2^{v-1}})$, and choose $\rho_v$ such that $32\rho_v<2^{v-2}$. Then there exists a collection of balls $S_{v}=\bigcup_{i\in I_v}B^e_{r_{i}}(\cdot)$ such that $\sum_{i\in I_v}r_{i}\leq 32\rho_v $, and for every $z \in (\mathbb{R}^4\backslash B) \backslash S_v$, it holds that
    \begin{equation*}
        \int_{A_{2^{v-1}}^{2^{v+2}}}\log|y-z|d\mu(y)\geq\mu_v\log \rho_v. 
    \end{equation*}\par
  In particular, under the assumption that $32\rho_v < 2^{v-2}$, we can select a suitable $\rho_v$ such that $\frac{\mu_v}{n} \log\frac{6 \cdot 2^v}{\rho_v}$ equals the prescribed $k \geq 7$; that is,
    \begin{align}
    k = \frac{\mu_v}{n} \log\frac{6 \cdot 2^v}{\rho_v}. \label{kZhenghaoshizhege}
    \end{align}
    We then select those balls in $S_v$ that intersect $A^{2^{v+1}}_{2^{v}}$. Denote this subcollection by $S_{v,k} = \bigcup_{i \in I_{v,k}} B^e_{r_i}(\cdot)$. The set $S_{v,k}$ still satisfies
\begin{align}
\sum_{i \in I_{v,k}} r_i \leq 32\rho_v < 2^{v-2}.
\end{align}\par
Since $t+|y| \leq 2^{v+1}+2^{v+2}=6\cdot 2^v$, for every $z \in A_{2^{v}}^{2^{v+1}} \backslash S_{v,k}$, we have
    \begin{align}
        \int_{A_{2^{v-1}}^{2^{v+2}}}\log \frac{|z-y|}{t+|y|}d\mu(y)
        =&\int_{A_{2^{v-1}}^{2^{v+2}}}\log\left[|z-y|\cdot\frac{1}{t+|y|}\right]d\mu(y)\label{20260904777}\\
        \geq&\mu_v\left(\log\rho_v+\log\frac{1}{6\cdot2^v}\right)\nonumber\\
        =&-\mu_v\left(\log\frac{6\cdot2^v}{\rho_v}\right)=-nk. 
    \end{align}\par
    Note that if $\mu_v=0$ for some $v$, then (\ref{20260904777}) vanishes for every $z \in A_{2^{v}}^{2^{v+1}}$ without removing any set. Henceforth, without loss of generality, we assume that $\mu_v>0$ for every $v\in \mathbb{Z}^+$.
    
    Now, we estimate the projected area $\theta_{v,k}$ of $S_{v,k}$. Here, since the radii satisfy $\sum_{i\in I_{v,k}}r_{i}\leq 32\rho_v <2^{v-2}$, each ball $B^e_h(\cdot)\in S_{v,k}$ intersecting $A^{2^{v+1}}_{2^{v}}$ lies entirely in $\mathbb{R}^4\backslash B^e_{2^{v-1}}(0)$. Therefore, for every radius $h\in\{r_{i}\}_{i\in I_{v,k}}$, the corresponding arc length $L(B^e_h(\cdot))$ on $\partial B_1^e$ is at most 
    \begin{align}
        L(B^e_h(\cdot))\leq2\arcsin\frac{h}{2^{v-1}} < \frac{2\pi h}{2^v}.\label{yigexianxingshangjie}    
    \end{align}
     This bound is obtained by fixing the radius $h$ and translating the center of $B_h^e(\cdot)$ to $z_0\in\partial B_{2^{v-1}}^e$ to increase the projected arc length, and then estimating the angular span.
    \par
    Denote by $P(B^e_h(\cdot))$ the area of the projection of a ball $B^e_h(\cdot)$ onto $\partial B_1^e$. We will show below that the total projected area of $S_{v,k}$ satisfies
    \begin{align}
    \sum_{i\in I_{v,k}}P\left(B^e_{r_i}(\cdot)\right)\leq\sum_{i\in I_{v,k}}P\left(B^e_{r_i}(z_0)\right)\leq P\left(B^e_{\sum_{i\in I_{v,k}}r_i}(z_0)\right), \label{jiandanjihe}
\end{align}
where $z_0 \in \partial B^e_{2^{v-1}}$.\par
    Denote by $L(B^e_h(\cdot))$ the supremum of the lengths of geodesic arcs in the dome-shaped projection of $B^e_h(\cdot)$ onto $\partial B_1^e$. From (\ref{yigexianxingshangjie}), we know that
    \begin{align}
        L\left(B^e_{\sum_{i\in I_{v,k}}r_i}(z_0)\right)=2\arcsin\frac{\sum_{i\in I_{v,k}}r_i}{2^{v-1}}<&\frac{2\pi\cdot 32 \rho_v}{2^{v}}\left(<\frac{\pi}{2}\right)\nonumber\\:=&C_1 \frac{\rho_v}{6\cdot2 ^v}.\label{LUpperbd}
    \end{align}\par
    
   In fact, the quantities $P$ and $L$ defined above for a single ball satisfy the relation $P(L)=\pi\bigl(L-\sin L\bigr)$ for $0\le L\le2\pi$. This follows from the polar-coordinate integral
    $$
        P(L)=4\pi\int_0^{L/2}\sin^2\varphi\,d\varphi=\pi L-\pi\sin(L).
    $$
    Thus, for $z_0\in\partial B^e_{2^{v-1}}$, we move the center of each ball $B_h^e(\cdot)\in S_{v,k}$ to $z_0$, resulting in $B_h^e(z_0)$, which has larger projected area and arc length. We then compute the projected arc length $L(h)$ and the projected area $P(L(h))$ of $B_h^e(z_0)$ as follows:
    \begin{align}
        L(h)=&2\arcsin \frac{h}{2^{v-1}},\nonumber\\
         P(L(h))=&\pi L(h)-\pi \sin(L(h))
          \nonumber\\
          =&2\pi\arcsin\frac{h}{2^{v-1}}-\pi \sin(2\arcsin\frac{h}{2^{v-1}})\label{20260827}
    \end{align}
and we know that the function $f(x)=2\pi\arcsin x - \pi \sin(2\arcsin x)$ is increasing and strictly convex for $x\in[0,1]$ and $f(0)=0$. Such functions are superadditive. From this, we deduce that the last inequality in  (\ref{jiandanjihe}) holds. As for the first inequality in (\ref{jiandanjihe}), consider the ball $B^e_r(y)$. A calculation similar to (\ref{20260827}) shows that the projected area of $B^e_h(y)$ on the unit sphere is
\begin{align*}
    2\pi\arcsin\frac{h}{|y|}-\pi \sin(2\arcsin\frac{h}{|y|}),
\end{align*}
which shows that moving the center $y$ closer to the origin while keeping the radius $h$ fixed increases the projected area on the unit sphere.\par

  Obviously, there exists $C_0>0$ such that for all $L>0$,  $P(L) = \pi(L - \sin L) \le C_0 L$. Consequently, using (\ref{jiandanjihe}) and (\ref{LUpperbd}), the total projected area $\theta_{v,k}$ of $S_{v,k}$ on $\partial B_1^e$ satisfies
    \begin{align}
    \theta_{v,k}=&\sum_{i\in I_{v,k}}P\left(B^e_{r_i}(\cdot)\right)\leq\sum_{i\in I_{v,k}}P\left(B^e_{r_i}(z_0)\right)\leq P\left(B^e_{\sum_{i\in I_{v,k}}r_i}(z_0)\right)\nonumber\\=&\pi\left(L\left(B^e_{\sum_{i\in I_{v,k}}r_i}(z_0)\right)\right)-\pi \sin\left(L\left(B^e_{\sum_{i\in I_{v,k}}r_i}(z_0)\right)\right)\nonumber\\
   \leq&C_0 L\left(B^e_{\sum_{i\in I_{v,k}}r_i}(z_0)\right) \le C_0 C_1 \left( \frac{\rho_v}{6\cdot2^v} \right). \label{upperboundsdsd}
    \end{align}
     Then by (\ref{upperboundsdsd}), (\ref{kZhenghaoshizhege}) and $e^{-xy}\leq y^{-1}e^{-x}$, ($x\geq1$, $y\geq1$), we obtain
    \begin{align}
        \theta_{v,k}\leq C_0 C_1 \left( \frac{\rho_v}{6\cdot2^v} \right)=C_0C_1\exp\left(-\frac{nk}{\mu_v}\right)\leq C_0C_1\frac{\mu_v}{n}e^{-k}.\label{UpperBoundofThetaVK}
    \end{align}
    Finally, we take the union of all $S_{v,k}$ to obtain $E_k = \bigcup_{v=0}^{\infty} S_{v,k}$. It is clear that for every $k\geq7$ and any $z\in(\mathbb{R}^4\backslash B)\backslash E_k$, there exists a $v_z\in\mathbb{Z}_{\geq0}$ such that $A_{t/2}^{2t} \subseteq A_{2^{v_z-1}}^{2^{v_z+2}}$ and
\begin{equation*}
0 \geq \int_{A_{t/2}^{2t}} \log\frac{|y-z|}{t+|y|}  d\mu(y)
\geq \int_{A_{2^{v_z-1}}^{2^{v_z+2}}} \log\frac{|y-z|}{t+|y|} d\mu(y)
\geq -nk,
\end{equation*}
where $E_k$ has total projected area bounded by
\begin{align}
    \sum_{v=0}^{\infty}\theta_{v,k}\leq\sum_{v=0}^{\infty}C_0C_1\frac{\mu_v}{n}e^{-k}\label{UpperboundProjAreaOfSingleK}\\
    \leq3C_0C_1e^{-k}:=M_0e^{-k}
\end{align}
    on $\partial B_1^e$. Here $M_0$ does not rely on $k$.
\end{proof}
\section{Asymptotic Radial Symmetry of the Conformal Factor of a Normal Metric}\label{SectionAymRadialSymmetry}
\begin{definition}
An $\mathscr{E}$-set for the end $\mathbb{R}^4\backslash B$ is an ambient set $E\subseteq\mathbb{R}^4$ of the form
\[
E := \bigcup_{i=1}^{\infty} B^e_{r_i}(z_i)
\]
whose intersection with the end is the exceptional region, and such that the sum of the areas of the projections of its constituent balls onto the unit sphere $\partial B_1^e$ is finite:
\[
\sum_{i=1}^\infty Area_{S^3}(\pi(B^e_{r_i}(z_i))) <\infty.\]

\end{definition}

In this section, using the lemmas in Section \ref{SectionPreliminary} dealing with $\mathcal{L}(f)$, we will prove three key points:\par

1\label{11112}. The logarithmic potential $\mathcal{L}(f)$ possesses favorable radial symmetry: outside the $\mathscr{E}$-set $E$, the logarithmic potential $\mathcal{L}(f)$ converges uniformly to a radially symmetric function---the infimum of $\mathcal{L}(f)$ on the sphere $\partial B_r^e$---as the radius $r$ tends to infinity.\par

2. The $\mathscr{E}$-set $E$ is sparse near infinity: $$|\partial B^e_r\cap E|/|\partial B^e_r|\to 0, $$
as $r \to \infty$.

3. Even when the $\mathscr{E}$-set $E$ is included, the average of the logarithmic potential over spheres also converges uniformly to the same radially symmetric function mentioned in item 1, as the radius $r=|x|$ tends to infinity. 
\begin{theorem}
    \label{line:161HaymanCor}
    Let $f\in L^1(\mathbb{R}^4\backslash B)$ be non-negative, and let $\mathcal{L}(f)$ be defined by (\ref{defofLF}). Then the logarithmic potential $\mathcal{L}(f)$ satisfies
    \begin{align}
    \mathcal{L}(f)(x)=G_f(|x|)+o(1),\label{line:165}
    \end{align}
uniformly for $x\in \mathbb{R}^4\backslash (B\cup E)$ as $|x|\to\infty$, where $E$ is an $\mathscr{E}$-set for which the sum of the areas of the projections of its constituent balls onto $\partial B_1^e$ is finite, and 
    $$G_f(r)=\inf_{|x|=r}\mathcal{L}(f)(x).$$
\end{theorem}
\begin{proof}
    For any $x=t\vec n$ with $\vec n:=x/|x|$ and $t:=|x|$, we consider
    \begin{align}
        v(x):=-\mathcal{L}(f)(x)=\frac{1}{8\pi^2}\int_{\mathbb{R}^4\backslash B}\log\frac{|x-y|}{|y|}f(y)dy.        
    \end{align}
    Denote 
    $$S_v(r):=\sup_{|x|=r}v(x).$$
    Fix $t:=|x|$. Then $v(x)$ satisfies
    \begin{align*}
        v(x)=&\frac{1}{8\pi^2}\int_{\mathbb{R}^4\backslash B}\log\frac{|x-y|}{|y|}f(y)dy\\
        \leq&\frac{1}{8\pi^2}\int_{\mathbb{R}^4\backslash B}\log\frac{t+|y|}{|y|}f(y)dy\\
        =&\frac{1}{8\pi^2}\int_{\mathbb{R}^4\backslash B}\log(1+\frac{t}{|y|})f(y)dy.
    \end{align*}
    Taking the supremum of $v$ over the sphere $\partial B_t^e$, we obtain
    \begin{align}
            S_v(t)\leq\frac{1}{8\pi^2}\int_{\mathbb{R}^4\backslash B}\log(1+\frac{t}{|y|})f(y)dy.\label{S_vUpperBound}
    \end{align}
    We now use (\ref{S_vUpperBound}) to estimate the difference between $v(x)$ and $S_v(|x|)$.
    \begin{align}
        0\geq& v(x)-S_v(t)\nonumber\\
        \geq&\frac{1}{8\pi^2}\int_{\mathbb{R}^4\backslash B}\log\left( \frac{|t\vec n-y|}{|y|}\frac{|y|}{t+|y|} \right) f(y)dy\nonumber\\
        =&\frac{1}{8\pi^2}\int_{\mathbb{R}^4\backslash B}\log\left( \frac{|t\vec n-y|}{t+|y|}\right)f(y)dy\nonumber\\
        =&:J_1(t)+J_2(t)+J_3(t),\label{line:186}
    \end{align}
    where $$J_i(t):=\frac{1}{8\pi^2}\int_{I_i^t}\log\left( \frac{|t\vec n-y|}{t+|y|}\right)f(y)dy,~(i=1,2,3),$$
    $I_1^{t}:=\left\{y:|y|\leq t/2\right\}$, $I_2^{t}:=\{y:t/2<|y|<2t\}$, $I_3^{t}:=\left\{y:|y|\geq2t\right\}$. \par
    Since $f$ is non-negative, $\log\frac{t+|y|}{|t\vec n-y|}f(y)\leq\log\frac{t+|y|}{t-|y|}f(y)$ for $y\in I_1^{t}$. Combining this with
    \begin{align}
        \log{\frac{1+z}{1-z}}\leq3z,~0<z<\frac{1}{2},\label{LogBianchengxianxing}
    \end{align}
    we have
    \begin{align}
        0\leq-J_1(t)=&\frac{1}{8\pi^2}\int_{I_1^{t}}\log  \frac{t+\left| {y} \right|}{\left|t\vec n- {y} \right|} f\left( y \right) dy    \nonumber\\
        \leq&\frac{1}{8\pi^2}\int_{I_1^{t}}\log  \frac{1+\left| \frac{y}{t} \right|}{1-\left| \frac{y}{t} \right|} f\left( y \right) dy\nonumber\\
        \leq&\frac{3}{8\pi^2t}\int_{I_1^{t}}|y|f(y)dy=\frac{3}{t}\int_{0}^{t/2}rd\rho(r),\label{line:193}
    \end{align}
    where $\rho(r):=\rho(B_r^e\backslash B):=\frac{1}{8\pi^2}\int_{B_r^e\backslash B}f(y)dy$.
    Since $f\in L^1(\mathbb{R}^4\backslash B)$ and $f\geq0$, we may define $\rho_0:=\lim_{l\to\infty}\rho(l)<\infty$. For every $\varepsilon>0$, there exists $t_0>1$ such that $\rho(t)>\rho_0-\varepsilon$ for all $t/2\geq t_0$. Then, the right-hand side of (\ref{line:193}) satisfies 
    \begin{align}
        \frac{3}{t}\int_{1}^{t/2}r~d\rho(r)=&\frac{3}{t}\left\{\int_{1}^{t_0}rd\rho(r)+\int_{t_0}^{t/2}rd\rho(r)\right\}\nonumber\\
        \leq&\frac{3}{t}\left\{ t_0\cdot\rho_0+\frac{t}{2}\int_{t_0}^{t/2}d\rho(r) \right\}\nonumber\\
        \leq&\frac{3\rho_0 t_0}{t}+\frac{3}{2}\varepsilon. \label{line:200}
    \end{align}
    Combining (\ref{line:193}) and (\ref{line:200}), we see that
    \begin{equation}
    J_1(t) \to 0 ,\quad \text{uniformly as } t \to \infty. \label{line:204:J1}
    \end{equation}\par
Similarly, $J_3(t)$ can be treated with an analogous approach. For every $\varepsilon > 0$, choose $t_0$ such that for all $t = |x| > 2t_0$, the measure $\int_{2t}^{\infty} d\rho(s) < \varepsilon$.
Since $f$ is non-negative, for $y \in I_3^{t}$, we have $\log\frac{|y|+t}{|t\vec n-y|}f(y) \leq \log\frac{|y|+t}{|y|-t} f(y)$. Together with formula (\ref{LogBianchengxianxing}), this yields
    \begin{align}
        0\leq-J_3(t)=&\frac{1}{8\pi^2}\int_{I_3^{t}}\log\frac{t+|y|}{|t\vec n-y|} f\left( y \right) dy
        \nonumber\\
        \leq&\frac{1}{8\pi^2}\int_{I_3^{t}}\log \frac{1+\frac{t}{|y|}}{1-\frac{t}{|y|}} f\left( y \right) dy\nonumber\\
        \leq&\frac{3t}{8\pi^2}\int_{I_3^{t}}\frac{1}{|y|}f(y)dy=3t\int_{2t}^{\infty}\frac{1}{s}d\rho(s)\nonumber\\
        \leq&3t\int_{2t}^{\infty}\frac{1}{2t}d\rho(s)<\frac{3}{2}\varepsilon,\nonumber
    \end{align}
    which means 
    \begin{equation}
        J_3(t)\to0,\label{line:215:j3}
    \end{equation}
    uniformly as $t\to\infty$.\par
    Next, we consider $J_2(t)$. Set $\varepsilon:=1/N^2$. There exists $t_N>0$ such that $\rho(t) > \rho_0 - \frac{1}{N^2}$. Consequently, $\rho(I_2^{t}) \leq N^{-2}$ for any $t/2 > t_N$. Applying Lemma \ref{Lemma4} with $k := N \geq 7$, and $\mu:=\rho\cdot\chi(\mathbb{R}^4\backslash B^e_{t_N}(0))$, we obtain a set $E_N$ associated with $N$. Then, for any $x \in I_2^t \backslash E_N$, we have 
    \begin{equation}
        0\geq J_2(t)=\frac{1}{8\pi^2}\int_{I_2^{t}}\log\frac{|x-y|}{|x|+|y|}f(y)dy\geq-\frac{1}{N^2}\cdot N=-\frac{1}{N},\label{line:219}
    \end{equation}
    where 
    $$E_N=\bigcup_{v=0}^{\infty} S_{v,N}$$
    has a sum of projected areas on $\partial B_1^e$ of at most $M_0e^{-N}$.
    Denote $E := \bigcup_{N=7}^{\infty} E_N$, which is an $\mathscr{E}$-set. Evidently, for any $x \notin E$ with $t/2 > t_N$, $J_2(t)$ satisfies
    \begin{align*}
        0\geq J_2(t)\geq -\frac{1}{N},
    \end{align*}
    which shows that for $x\notin E$,
    \begin{align}
        J_2(t)\to0,~as~t\to\infty, \label{line:230}
    \end{align}
    where $E$ has a finite sum of projected areas of at most
    $$\sum_{N=7}^{\infty}M_0e^{-N}\leq M<\infty$$
    on $\partial B_1^e$. \par
    In summary, combining the results (\ref{line:204:J1}), (\ref{line:230}), and (\ref{line:215:j3}) for $J_1(t)$, $J_2(t)$, $J_3(t)$ with formula (\ref{line:186}), we conclude that
    $$v(x)=S_v(|x|)+o(1),~as~|x|\to\infty,$$
     uniformly for $x\in \mathbb{R}^4\backslash (B\cup E)$, where the $\mathscr{E}$-set $E$ has finite projected area on $\partial B_1^e$.\par 
    Finally, combining $v(x) = -\mathcal{L}(f)(x)$ and the identity $\sup\left(-\mathcal{L}(f)\right)= -\inf \mathcal{L}(f)$, we see immediately that the theorem holds.
\end{proof}
\begin{remark}\label{FinitnfinitInters}
     Since the sum of the projected areas of the balls forming the $\mathscr{E}$-set on $S^3$ is finite, we know that given an $\varepsilon>0$, there exists an $R_\varepsilon>0$ such that the sum of projected areas of $E\bigcap (\mathbb{R}^4\backslash B_r^e(0))$ on $S^3$ is no more than $\varepsilon$. Consequently, the set of directions $\theta_\infty\in S^3$ for which the ray  $\gamma_{\theta_\infty}:=\{r\theta_\infty:r\in (1,+\infty)\}$ intersects the $\mathscr{E}$-set $E$ outside every sufficiently large compact set has measure zero. We denote the collection of all such $\theta_\infty$ by $S^3_\infty$. We also define $S_0^3:=S^3\backslash S^3_\infty$, in which every element $\theta_0\in S_0^3$ has the property that the ray $\gamma_{\theta_0}:=\{r\theta_0:r\in (1,+\infty)\}$ does not intersect the $\mathscr{E}$-set for $r\geq R_{\theta_0}$, for some $R_{\theta_0}>0$. 
     
     Moreover, the ratio $|\partial B^e_r\backslash E|/|\partial B^e_r|\to 1 $, as $r\to\infty$, which means the $\mathscr{E}$-set $E$ is sparse near infinity.
\end{remark}
\begin{proof}
Still using the notation in Theorem \ref{line:161HaymanCor}, we know that
\begin{align*}
E_N=\bigcup_{v=0}^{\infty} S_{v,N}.
\end{align*}
Recalling the construction of $E_N$, we use the cut-off measure $\mu_N:=\frac{f(y)}{8\pi^2} \chi_{|y|>t_N}dy.$
Recall the proof in Lemma \ref{Lemma4}: each $S_{v,N}$ intersects the $A_{2^v}^{2^{v+1}}$-annulus defined before Lemma \ref{Lemma4}. If $n_N=0$, we know that $E_N=\varnothing$. Without loss of generality, we assume the sums below are taken for those $N$ such that $n_N>0$. By (\ref{UpperBoundofThetaVK}), the sum of projected areas $\theta_{v,N}$ of $S_{v, N}$ onto $\partial B^e_1$ satisfies
\begin{align}
\theta_{v,N}\leq C_0C_1\frac{\mu_{v,N}}{n_N}e^{-N},\label{Righetwwerwar}
\end{align}
where $\mu_{v,N}:=\mu_N\left({{A^{2^{v+2}}_{2^{v-1}}}}\right)$, and $n_N:=\mu_N(\mathbb{R}^4\backslash B)$. First, for every $\varepsilon>0$, we take $N_*$ large enough such that 
\begin{align}
    \sum_{N=N_*+1}^\infty 3C_0C_1 e^{-N}\leq \varepsilon/2.
\end{align}
For each $N_i\in\{7,\ldots,N_*\}$ with $n_{N_i}>0$,
choose $m_i\in\mathbb Z^+$ sufficiently large so that
\begin{align}
\sum_{v=m_i }^\infty\frac{\mu_{v,N_i}}{n_{N_i}}<\frac{\varepsilon}{2N_* C_0 C_1\sum_{N=7}^\infty e^{-N}}.\label{leftaeraerdf}
\end{align} For indices with $n_{N_i}=0$, we set $m_i=1$; no estimate involving division by $n_{N_i}$ is needed, since $E_{N_i}=\varnothing$.
Since $E:=\bigcup_{N=7}^{\infty} E_{N}$, we set $R_i:=2^{m_i+8}$ and remove $B^e_{R_i}(0)$ for each $E_N$. Thus for each $r>R_i$, the set $E_{N_i}\bigcap (\mathbb{R}^4\backslash B_r^e(0))$ lies entirely in $E_{N_i}\bigcap (\mathbb{R}^4\backslash B_{R_i/{2^4}}^e(0))=E_{N_i}\backslash B^e_{{2^{m_i+4}}}(0)$.
Let $R_\varepsilon:=\max\{R_i\}_{i=7}^{N_*}$ and $m_\varepsilon:=\max\{m_i\}_{i=7}^{N_*}$. Thus, for every $i=7,\cdots, N_*$, all sets $E_{N_i}\bigcap (\mathbb{R}^4\backslash B_r^e(0))$ lie in $(\bigcup_{i=7}^{N_*}E_{N_i})\bigcap (\mathbb{R}^4\backslash B_{R_\varepsilon}^e(0))$. By the construction of $\{S_{v,N}\}_{N=7}^\infty$, the set $\{\{S_{v,N}\}_{N=7}^\infty,v<m_\varepsilon\}$ lies entirely in $B_{2^{m_\varepsilon+4}}^e(0)$.

So, using (\ref{Righetwwerwar}) and (\ref{leftaeraerdf}), for every  $r>{R_\varepsilon}$, it follows that the sum of projected areas $\theta_r$ of $E\bigcap (\mathbb{R}^4\backslash B_r^e(0))$ on $\partial B^e_1$ satisfies
\begin{align*}
\theta_r &\leq \sum_{N=7}^\infty \sum_{v=m_\varepsilon+4}^\infty \theta_{v,N}
\leq \sum_{N=7}^\infty \sum_{v=m_\varepsilon}^\infty \theta_{v,N} \leq (\sum_{N=7}^{N_*}+\sum_{N=N_*+1}^\infty) \sum_{v=m_\varepsilon}^\infty C_0C_1\frac{\mu_{v,N}}{n_N}e^{-N}\\
&\leq \sum_{N=7}^{N_*} \sum_{v=m_\varepsilon}^\infty C_0C_1\frac{\mu_{v,N}}{n_N}e^{-N}+\sum_{N=N_*+1}^\infty  C_0C_1\cdot 3 \cdot e^{-N} \\
&\leq \sum_{N=7}^{N_*}  C_0C_1\frac{\varepsilon}{2N_* C_0 C_1\sum_{i=7}^\infty e^{-i}}e^{-N}+\sum_{N=N_*+1}^\infty 3 C_0C_1e^{-N}\\
&< \frac{\varepsilon}{2}+\frac{\varepsilon}{2}=\varepsilon,
\end{align*}
where we should notice that for every $N\geq 7$, 
\begin{align*}
    \sum_{v=1}^{\infty}\frac{\mu_{v,N}}{n_N}\leq 3.
\end{align*}
It means $\theta_r \to 0$ as $r\to \infty$.

Next, we discuss ${|\partial B^e_r \cap E|}/{|\partial B^e_r|}$ as $r\to\infty$. First, from the construction of $E$, we know that for $r>R_\varepsilon$, the balls contained in $E$ that intersect with  $\partial B^e_r$ lie entirely in $E\cap(\mathbb{R}^4\backslash B^e_{2^{m_\varepsilon+1}}(0))$. Thus, the ratio  satisfies
\begin{align*}
    \frac{|\partial B^e_r \cap E|}{|\partial B^e_r|}&\leq \sum_{N=7}^\infty \sum_{v=m_\varepsilon+1}^\infty \theta_{v,N}
\leq \sum_{N=7}^\infty \sum_{v=m_\varepsilon}^\infty \theta_{v,N} \\
&\leq (\sum_{N=7}^{N_*}+\sum_{N=N_*+1}^\infty) \sum_{v=m_\varepsilon}^\infty C_0C_1\frac{\mu_{v,N}}{n_N}e^{-N}\\
&\leq \sum_{N=7}^{N_*} \sum_{v=m_\varepsilon}^\infty C_0C_1\frac{\mu_{v,N}}{n_N}e^{-N}+\sum_{N=N_*+1}^\infty  C_0C_1\cdot 3 \cdot e^{-N} \\
&\leq \sum_{N=7}^{N_*}  C_0C_1\frac{\varepsilon}{2N_* C_0 C_1\sum_{i=7}^\infty e^{-i}}e^{-N}+\sum_{N=N_*+1}^\infty 3 C_0C_1e^{-N}\\
&< \frac{\varepsilon}{2}+\frac{\varepsilon}{2}=\varepsilon.
\end{align*}
Thus, we directly obtain that ${|\partial B^e_r \cap E|}/{|\partial B^e_r|}\to 0$, 
as $r \to \infty$.
\end{proof}
 In fact, due to the uniformity in (\ref{line:165}), the infimum term can be replaced by the spherical average outside the $\mathscr{E}$-set. Furthermore, inspired by \cite{ChangQingYangInvention2000}, we can show that the infimum differs from the spherical average by at most $o(1)$, as the radius $r$ tends to $\infty$.
\begin{proposition}\label{theorem:1}
    Suppose that $\mathcal{L}(f)(x)$ is defined as in (\ref{defofLF}), that $f$ is non-negative, and that $f\in L^1(\mathbb{R}^4\backslash B)$. Then we have
    \begin{equation}
        \mathcal{L}(f)(x) = \mathcal{A}(|x|) + o(1),~|x|\to\infty,
    \end{equation}
    uniformly for $x\in \mathbb{R}^4\backslash (B\cup E)$, where
    \begin{align}
        \mathcal{A}(r):=\fint_{\partial B^e_{r}} {\mathcal{L}(f)}{d\sigma}=\frac{1}{|\partial B^e_r|}\int_{\partial B^e_r}\mathcal{L}(f)(y)d\sigma_y, \label{defofA}
    \end{align}
    $|\partial B^e_r|:=2\pi^2r^{3}$ is the Euclidean area of $\partial B^e_r$, and $E$ is the $\mathscr{E}$-set mentioned in Theorem  \ref{line:161HaymanCor}.
\end{proposition}
 
 \begin{proof}
    Denote $r:=|x|$. From Theorem  \ref{line:161HaymanCor}, we have for $x\in \mathbb{R}^4\backslash (B\cup E)$, 
    $$\mathcal{L}(f)(x)=\inf_{|y|=|x|}\mathcal{L}(f)(y)+o(1),$$
    uniformly as $|x|\to\infty$. Due to this property of uniformity, we know that 
\begin{align}
    \mathcal{L}(f)(x)&=\frac{1}{|\partial B^e_{|x|}\backslash E|}\int_{\partial B^e_{|x|}\backslash E}\mathcal{L}(f)(y){dS_y}+o(1)\nonumber\\
    :&=\fint_{\partial B^e_{r}\backslash E}\mathcal{L}(f)(y)d\sigma_y+o(1).
\end{align}
    We only need to prove that
    \begin{align}
                \fint_{\partial B^e_{r}\backslash E}\mathcal{L}(f)(y)d\sigma_y=\fint_{\partial B^e_{r}}\mathcal{L}(f)(y)d\sigma_y+o(1).\label{yougeyizhix}
    \end{align}\par
First, we split the integration domain of $\mathcal{L}(f)(x)$ into two parts.
    \begin{align}
        \mathcal{L}(f)(x)=&
        \frac{1}{8\pi^2}\int_{\mathbb{R}^4\backslash B}\log\frac{|y|}{|x-y|}f(y)dy\nonumber\\
        =&\frac{1}{8\pi^2}\int_{B^e_{r/2}\backslash B}\log\frac{|y|}{|x-y|}f(y)dy\nonumber\\
        &+\frac{1}{8\pi^2}\int_{\mathbb{R}^4\backslash  B^e_{r/2}(0)}\log\frac{|y|}{|x-y|}f(y)dy\nonumber\\
        :=&\mathcal{L}_1(x)+\mathcal{L}_2(x). \label{Neewline:355}
    \end{align}
Then we split $\mathcal{L}_1(x)$ into two parts
\begin{align}
    \mathcal{L}_1(x)=&\frac{1}{8\pi^2}\int_{B^e_{r/2}\backslash B}\log\frac{|y|}{|x-y|}f(y)dy\nonumber\\
    =&\frac{1}{8\pi^2}\int_{B^e_{r/2}\backslash B}\log\left[\frac{|y|}{|x|}\frac{|x|}{|x-y|}\right]f(y)dy\nonumber\\
    =&\frac{1}{8\pi^2}\int_{B^e_{r/2}\backslash B}\log\frac{|y|}{|x|}f(y)dy\nonumber\\&+\frac{1}{8\pi^2}\int_{B^e_{r/2}\backslash B}\log\frac{|x|}{|x-y|}f(y)dy\nonumber\\
    :=&\mathcal{L}_{11}(|x|)+\mathcal{L}_{12}(x), \label{Neewline:362}
\end{align}
where $\mathcal{L}_{11}$ is a radially symmetric function. We will show that
\begin{align}
    \mathcal{L}_{12}(x)=o(1),\label{Neewline:366}    
\end{align}
 uniformly as $|x|\to \infty$. To show this, we divide $B^e_{|x|/2}$ into $A^{0.5|x|}_{\log|x|}\bigcup B^e_{\log|x|}$, where 
\begin{align*}
A^{0.5|x|}_{\log|x|}&=\left\{y\in\mathbb{R}^4\backslash B:\log|x|\leq|y|<0.5|x|\right\},\\
B^e_{\log|x|}&=\left\{y\in\mathbb{R}^4\backslash B:|y|<\log|x|\right\},
\end{align*}
with $|x|\geq e$.
For $y\in A^{0.5|x|}_{\log|x|}$, we have
$$0.5|x|\leq|x|-|y|\leq|x-y|\leq |x|+|y|\leq|x|+0.5|x|.$$
Dividing the above inequalities by $|x-y|$, we obtain $2/3\leq|x|/|x-y|\leq 2$, which means
\begin{align}
    \left|\log\frac{|x|}{|x-y|} \right|\leq\log 2,~ \forall~y\in  A^{0.5|x|}_{\log|x|}.\label{Neewline:371}
\end{align}
For $y\in B^e_{\log|x|}$, it is easy to see that
\begin{align}
    \frac{|x|}{|x|+\log|x|}\leq\frac{|x|}{|x|+|y|}\leq\frac{|x|}{|x-y|}\leq\frac{|x|}{|x|-|y|}\leq\frac{|x|}{|x|-\log|x|}.
\end{align}
    It follows that
    \begin{align}
        \left|\log\frac{|x|}{|x-y|} \right|\leq\max\left\{-\log\frac{|x|}{|x|+\log|x|},\log\frac{|x|}{|x|-\log|x|}\right\}=o(1), \label{Neewline:379}
    \end{align}
    uniformly as $|x|\to \infty$ for $y\in B^e_{\log|x|}$.
   For sufficiently large $|x|$, using (\ref{Neewline:371}), (\ref{Neewline:379}), and the integrability of $f$, which implies $||f||_{L^1(A^{0.5|x|}_{\log|x|})}=o(1)$, as $|x|\to\infty$, we obtain
    \begin{align}
        |\mathcal{L}_{12}(x)|=&\left|\frac{1}{8\pi^2}\int_{B^e_{r/2}\backslash B}\log\frac{|x|}{|x-y|}f(y)dy\right| \nonumber\\
        \leq&\frac{1}{8\pi^2}\int_{B^e_{r/2}\backslash B}\left|\log\frac{|x|}{|x-y|}\right|f(y)dy\nonumber\\
        \leq&\frac{1}{8\pi^2}\int_{A^{0.5|x|}_{\log|x|}\bigcup B^e_{\log|x|}}\left|\log\frac{|x|}{|x-y|}\right|f(y)dy\nonumber\\
        \leq&\frac{\log 2}{8\pi^2}\int_{A^{0.5|x|}_{\log|x|}}f(y)dy\nonumber+\frac{o(1)}{8\pi^2}\int_{B^e_{\log|x|}}f(y)dy\nonumber\\
        =&o(1), \label{Neewline:3922}
    \end{align}
which means (\ref{Neewline:366}) is valid. \par
Second, we prove that 
\begin{align}
    \fint_{ \partial B^e_{r}\backslash E}  \mathcal{L}_2(y)d\sigma_y=o(1), \label{Neewline:397}
\end{align}
as $r\to\infty$.
Using Fubini's theorem, we have
\begin{align}
        ~~~~&\fint_{ \partial B^e_{r}\backslash E}  \mathcal{L}_2(y)d\sigma_y\nonumber\\
        &=\frac{1}{8\pi^2}\int_{\mathbb{R}^4\backslash B^e_{r/2}(0)}f(y)\fint_{\partial B^e_{r}\backslash E}\log\frac{|y|}{|x-y|}d\sigma_xdy\nonumber\\
        &:=\frac{1}{8\pi^2}\int_{\mathbb{R}^4\backslash B^e_{r/2}(0)}f(y)\bar Pdy,\label{Neewline:403}
\end{align}
where
$$\bar P(y)=\fint_{\partial B^e_{r}\backslash E}\log\frac{|y|}{|x-y|}d\sigma_x.$$ We will show that $\bar P(y)$ is bounded for $|y|\geq r/2$. Again, we divide the integral region $\partial B^e_{r}\backslash E$ into two parts
 \begin{align}
     |\bar P(y)|\leq&\fint_{\partial B^e_{r}\backslash E}\left|\log\frac{|y|}{|x-y|}\right|d\sigma_x \nonumber\\
     =&\frac{1}{|\partial B^e_{r}\backslash E|}\int_{(\partial B^e_{r}\backslash E)\backslash\{|x-y|\leq0.25|y|\}}\left|\log\frac{|y|}{|x-y|}\right|d\sigma_x
     \nonumber\\&+\frac{1}{|\partial B^e_{r}\backslash E|}\int_{(\partial B^e_{r}\backslash E)\bigcap\{|x-y|\leq0.25|y|\}}\left|\log\frac{|y|}{|x-y|}\right|d\sigma_x\nonumber\\
     :=&\bar P_1(y)+\bar P_2(y).\label{Neewline:411}
 \end{align}\par
For the term $\bar P_1(y)$, note that the integration variable $x$ lies in $\partial B^e_{r}$. Hence, for every $x \notin \{x:|x-y|\leq 0.25|y|\}$, i.e., for $|x-y|\geq 0.25|y|$, we have
\begin{align}
     \log\frac{|y|}{|x-y|}\leq \log4.\label{shangjieLOg4}
\end{align}
Since $|y|\geq0.5|x|$, we have
\begin{align}
        \log\frac{|y|}{|x-y|}\geq\log\frac{|y|}{|x|+|y|}\geq\log\frac{|y|}{2|y|+|y|}=-\log 3.\label{xiajielog3}
\end{align}
So by the bound (\ref{shangjieLOg4}) and (\ref{xiajielog3}), we have 
\begin{align}
    \bar P_1(y)&=\frac{1}{|\partial B^e_{r}\backslash E|}\int_{(\partial B^e_{r}\backslash E)\backslash\{|x-y|\leq0.25|y|\}}\left|\log\frac{|y|}{|x-y|}\right|d\sigma_x\nonumber\\
    &\leq \frac{|\partial B^e_{r}\backslash E\backslash \{|x-y|\leq0.25|y|\}|}{|\partial B^e_{r}\backslash E|}\cdot \log 4\leq\log 4, \label{Neeewline:425}
\end{align}
for $x\in (\partial B^e_{r} \backslash E) \backslash\{x:|x-y|\leq 0.25|y|\}$.\par
For $\bar P_2(y)$, take any $x \in (\partial B^e_r \backslash E)\bigcap \{x:|x-y|$ $\leq 0.25|y|\}$. Since $r=|x|$, $|x|$ and $|y|$ satisfy 
\begin{align}
    |y|\leq |x|+|y-x|\leq 0.25|y|+|x|,\label{zheyishi}
\end{align}
and
\begin{align}
    |x|\leq |y|+|x-y|\leq |y|+ 0.25|y| .\label{zheershi}
\end{align}
Combining (\ref{zheyishi}) and (\ref{zheershi}), it follows that the ratio $r/|y|$ has the bounds
\begin{align}
    \frac{3}{4}\leq r/|y| \leq \frac{5}{4}. \label{labubu}
\end{align}
Using an integral transformation $z=x/|y|$, we have 
\begin{align}
        \bar P_2(y)&=\frac{1}{|\partial B^e_{r}\backslash E|}\int_{(\partial B^e_{r}\backslash E)\bigcap\{x:|x-y|\leq0.25|y|\}}\left|\log\frac{|y|}{|x-y|}\right|d\sigma_x\nonumber\\
        &=\frac{1}{|(\partial B^e_{r}\backslash E)/|y||}\int_{((\partial B^e_{r}\backslash E) /|y| 
        )\bigcap\{z:|z-y/|y||\leq0.25\}}\left|\log\frac{1}{|z-y/|y||}\right|d\sigma_z\nonumber\\
        &=\frac{|\partial B^e_{r}/|y||}{|(\partial B^e_{r}\backslash E)/|y||}\cdot \frac{1}{|\partial B^e_{r}/|y||}\nonumber\\
        &~~~~\cdot \int_{((\partial B^e_{r}\backslash E) /|y| 
        )\bigcap\{z:|z-y/|y||\leq0.25\}}\left|\log\frac{1}{|z-y/|y||}\right|d\sigma_z\nonumber\\
        &\leq \tilde C \int_{(\partial B^e_{r} /|y| 
        )\bigcap\{z:|z-y/|y||\leq0.25\}}\left|\log\frac{1}{|z-y/|y||}\right|d\sigma_z
        \leq C. \label{Neewline:430}
 \end{align}
Here the boundedness of $\tilde C$ follows from Remark \ref{FinitnfinitInters} together with (\ref{labubu}). The constant $C$ is independent of $y$: by (\ref{labubu}), the radius $r/|y|$ of the sphere $\partial B^e_{r/|y|}$ containing the integration domain in (\ref{Neewline:430}) lies between $3/4$ and $5/4$, and the term $\log\frac{1}{|z-y/|y||}$ is integrable on this domain.
Consequently, from the boundedness estimates (\ref{Neeewline:425}) and (\ref{Neewline:430}) for $\bar P_1(y)$ and $\bar P_2(y)$, we conclude that $\bar P(y)$ is bounded. Hence, letting $r \to \infty$ in (\ref{Neewline:403}) and using the integrability of $f$, we see that (\ref{Neewline:397}) holds. \par
Similarly, following the same method above for $\fint_{ \partial B^e_{r}} \mathcal{L}_2(y)\,d\sigma_y$, we also obtain
\begin{align}
\fint_{ \partial B^e_{r}} \mathcal{L}_2(y)\,d\sigma_y = o(1), \qquad \text{as } r \to \infty, \label{Newline:436}
\end{align}
since in this case we may simply take $E = \varnothing$.

With the above preparations, we can now prove the proposition. The references below indicate the equalities applied at each step:
Using the decomposition $\mathcal{L}(f)=\mathcal{L}_1+\mathcal{L}_2$ from (\ref{Neewline:355}),
the $o(1)$ property (\ref{Neewline:397}) of the spherical mean of $\mathcal{L}_2$ on $\partial B_r^e\backslash E$,
the further splitting $\mathcal{L}_1=\mathcal{L}_{11}+\mathcal{L}_{12}$ from (\ref{Neewline:362}),
the radial symmetry of $\mathcal{L}_{11}(|x|)$,
together with $\mathcal{L}_{12}=o(1)$ from (\ref{Neewline:3922}),
the fact that the spherical mean of $\mathcal{L}_{2}$ is $o(1)$ as shown in (\ref{Newline:436}),
and the definition (\ref{Neewline:355}) as stated earlier,
we obtain
\begin{align}
    &\fint_{ \partial B^e_{r}\backslash E}  \mathcal{L}(f)(y)d\sigma_y=\fint_{ \partial B^e_{r}\backslash E}  \mathcal{L}_1(y)+\mathcal{L}_2(y)d\sigma_y\nonumber\\
    =&\fint_{ \partial B^e_{r}\backslash E}  \mathcal{L}_1(y)d\sigma_y+o(1)=\fint_{ \partial B^e_{r}\backslash E} \mathcal{L}_{11}(|y|)+\mathcal{L}_{12}(y)d\sigma_y+o(1)\nonumber\\
=&\fint_{ \partial B^e_{r}\backslash E} \mathcal{L}_{11}(|y|)d\sigma_y+o(1)=\fint_{ \partial B^e_{r}} \mathcal{L}_{11}(|y|)d\sigma_y+o(1)\nonumber\\
=&\fint_{ \partial B^e_{r}} \mathcal{L}_{11}(|y|)+\mathcal{L}_{12}(y)d\sigma_y+o(1)=\fint_{ \partial B^e_{r}}\mathcal{L}_1(y) d\sigma_y+o(1)\nonumber\\
=&\fint_{ \partial B^e_{r}}  \mathcal{L}_1(y)+\mathcal{L}_2(y)d\sigma_y+o(1)=\fint_{ \partial B^e_{r}}  \mathcal{L}(f)(y)d\sigma_y+o(1),\label{pingjunchayidian}
\end{align}
where $o(1)\to0$, as $r=|x|\to\infty$. \par
Finally, from (\ref{pingjunchayidian}) and (\ref{yougeyizhix}), we now complete the proof.
\end{proof}

\section{The Chern--Gauss--Bonnet Theorem and a Huber-Type Theorem}\label{ChapterMaintheorem}
In this section, we mainly consider the conformally flat  end $(\mathbb{R}^4\backslash B,g=e^{2u(x)}|dx|^2)$ of the complete, locally conformally flat manifold $(M,g)$ that satisfies 
 \begin{equation}
     R[g]\geq0,\label{line:244}
 \end{equation} 
 and
 \begin{equation}
     Q[g]\geq 0\label{line:248}
 \end{equation}
 near infinity. We first truncate the end at a sufficiently large coordinate sphere so that both inequalities hold everywhere on the remaining end, and then rescale the coordinates to make its inner boundary the unit sphere. Unless otherwise stated, the end and the coordinates in this section refer to this normalized end. Thus $R[g]\geq0$ and $Q[g]\geq0$ hold on all of $\mathbb{R}^4\backslash B$ when Lemma \ref{NormalMetric} and the non-negative-potential estimates are applied. The final step in the proof of Theorem \ref{MainTheorem} shows that the volume limit and the curvature formula are unaffected by the choice of the original inner boundary.\par
Under the assumptions (\ref{line:244}) and (\ref{line:248}), we introduce notation for several basic quantities and describe the relations among them. These quantities and relations will be used frequently in this section. By Lemma \ref{NormalMetric}, the conformal factor $u(x)$ admits an expansion such that
    \begin{equation}
        u(x)=\frac{1}{8\pi^2}\int_{\mathbb{R}^4\backslash B}\log\frac{|y|}{|x-y|}2Q_g(y)e^{4u(y)}dy+\gamma\log|x|+b(x/|x|^2),\label{line:263}
    \end{equation}
    where $\gamma$ is a constant and $b(x)$ is biharmonic on $\bar B$ in the sense specified in Lemma \ref{NormalMetric}. For the logarithmic potential term appearing in (\ref{line:263}), we denote $$f(z):=2Q_g(z)e^{4u(z)},$$
    $$
        \mathcal{L}(f)(x):=\frac{1}{8\pi^2}\int_{\mathbb{R}^4\backslash B}\log\frac{|y|}{|x-y|}2Q_g(y)e^{4u(y)}dy.
    $$
    Applying Theorem \ref{line:161HaymanCor} and Proposition \ref{theorem:1}, for any $x\notin E$, we have 
    \begin{align}
        u(x)=G(|x|)+o(1)=D(|x|)+o(1),~|x|\to\infty, \label{line:269}
    \end{align}
    where 
    \begin{equation}
        G(|x|):=G_f(|x|)+\gamma\log |x|+\inf_{|y|=|x|}b(y/|y|^2),\label{line:272DefOfG}
    \end{equation}
    and $G_f(r)$ is the infimum of the logarithmic potential on $\partial B_r^e$ as defined in Theorem \ref{line:161HaymanCor}, and
    \begin{equation}
        D(|x|):=\mathcal{A}(|x|)+\gamma\log |x|+\frac{1}{|\partial B^e_r|}\int_{\partial B^e_r}  b(\frac{y}{|y|^2})d\sigma_r\label{line:272DefOfD}
    \end{equation}
    is the spherical average of $u$, and
    $\mathcal{A}(s)$ is defined in Proposition \ref{theorem:1}.
    In particular, the latter equality of (\ref{line:269})
    \begin{align}
        G(r)=D(r)+o(1) \label{GxDx}
    \end{align}
    is valid for sufficiently large $r:=|x|$, since $b(z)$ is biharmonic on $\bar B$ in the sense of Lemma \ref{NormalMetric} and hence $b(x/|x|^2)\to b(\vec0)$, as $|x|\to\infty$. It is not hard to see that 
    \begin{align}
        G(|x|)&\leq D(|x|),\label{LleqD}\\
        G(|x|)&\leq u(x), \label{Llequ}
    \end{align}
    hold for every $x\in\mathbb{R}^4\backslash B$. Indeed, for $r=|x|$, we have $G_f(r)\leq\mathcal{L}(f)(x)$ and $\inf_{|y|=r}b(y/|y|^2)\leq b(x/|x|^2)$. Adding these inequalities and $\gamma\log r$ gives $G(r)\leq u(x)$; averaging over the sphere gives $G(r)\leq D(r)$. Thus $G$ is a radial lower bound for $u$, and need not equal the spherical infimum of $u$.
    \par  

\begin{lemma}\label{wuxianchangxingzhi}
    Suppose a conformally flat end $M_i:=(\mathbb{R}^4\backslash B,g=e^{2u}|dx|^2)$ satisfies $R[g]\geq0$ and $Q[g]\geq0$, and assume that $M_i$ is an end of a complete manifold $(M,g)$. Then, for every $R\geq1$, 
    \begin{align}
        \int_{R}^\beta e^{D(s)}ds\to+\infty,\,
        \int_{R}^\beta e^{G(s)}ds\to+\infty,\label{wanbeiDG}
    \end{align}
    as $\beta\to\infty$.
\end{lemma}
\begin{proof}
    Fix $R>0$. Since $(M,g)$ is complete, for every $\theta\in S^3$ we have
    \begin{align}
        \int_{R}^{\beta}e^{u(s\theta)}ds\to\infty, \label{temp1}
    \end{align}
    as $\beta\to\infty$. From (\ref{line:269}), for $x\notin E$, we know that 
    \begin{align}
        u(x)=G(|x|)+o(1)=D(|x|)+o(1),
    \end{align}
    where $o(1)\to0$, as $|x|\to\infty$.
    Hence, for $\varepsilon = 1/2$, there exists $R_1 > 0$ such that for every $x\notin E$,  $|x| \ge R_1$,
    \begin{align}
        u(x)\leq G(|x|)+\frac{1}{2}.\label{temp2}
    \end{align}
    By Remark \ref{FinitnfinitInters}, choose $\theta_0\in S_0^3$ and $R_2\geq1$ such that $s\theta_0\notin E$ for every $s\geq R_2$.
    Set $R_0=\max\{R_1,R_2\}$. Consequently,
    \begin{align}
        \int_{R_0}^{\beta}e^{G(s)+\frac{1}{2}}ds\geq\int_{R_0}^{\beta}e^{u(s\theta_0)}ds\to\infty,
    \end{align}
    as $\beta\to\infty$.\par
    Similarly, using the fact that $u(x)=D(|x|)+o(1)$ for $x\notin E$, we can show that for every $R\ge 1$,
    \begin{equation*}
        \int_{R}^\beta e^{D(s)}ds\to\infty,
    \end{equation*}
    as $\beta\to\infty$.
\end{proof}    

  The following result shows that the $g$-volume of the $\mathscr{E}$-set has little influence on our discussion, which makes a great contribution to our main theorem. 
  \begin{lemma}\label{FiniteVolumeE}
     Suppose a conformally flat end $M_i:=(\mathbb{R}^4\backslash B,g=e^{2u}|dx|^2)$ satisfies $R[g]\geq0$ and $Q[g]\geq0$, and assume that $M_i$ is an end of a complete manifold $(M,g)$. Then, in the infinite-volume case, the $g$-volume of the $\mathscr{E}$-set $E$ defined in Proposition  \ref{theorem:1} is negligible compared with the volumes of the Euclidean ball measured with respect to the metrics  $\widetilde{g}=e^{2D(|x|)}|dx|^2$ and $\hat{g}=e^{2G(|x|)}|dx|^2$. That is, if $V_g( B_R^e)\to \infty$ as $R\to \infty$, then we have $V_g(B_R^e\cap E)=o(V_{\widetilde{g}}(B_R^e))$ and $V_g(B_R^e\cap E)=o(V_{\hat{g}}(B_R^e))$ as $R\to\infty$.
  \end{lemma}
\begin{proof}
All volumes below are taken over the end-truncated balls $B_R^e$
fixed in the global notation above.
For $s>1$, define
\[
\begin{aligned}
a(s)&:=\int_{\partial B_s^e}e^{4u(x)}\,d\sigma_x,\\
b(s)&:=\int_{\partial B_s^e\setminus E}
             e^{4u(x)}\,d\sigma_x,\\
w(s)&:=|\partial B_s^e|e^{4D(s)},\\
z(s)&:=|\partial B_s^e|e^{4G(s)}.
\end{aligned}
\]
These functions are nonnegative and locally integrable.

The spherical exponential estimate used above, obtained from
\cite[Lemma 1.8, (1.20)]{ChangQingYangInvention2000}, gives
\begin{equation}
\frac{a(s)}{w(s)}=\frac{\int_{\partial B_s^e}e^{4u(x)}\,d\sigma_x}{|\partial B_s^e|e^{4D(s)}}\longrightarrow1
\qquad\text{as }s\to\infty.
\label{20260904444}
\end{equation}
We next show that $$\frac{b(s)}{w(s)}=\frac{\int_{\partial B_s^e\setminus E}e^{4u(x)}\,d\sigma_x}{|\partial B_s^e|e^{4D(s)}}\to1,$$ as $s\to \infty$.
By \eqref{line:269}, for every $\varepsilon>0$ there exists
$S_\varepsilon>1$ such that
\[
|u(x)-D(|x|)|\le\varepsilon
\qquad
\text{whenever }x\notin E,\quad |x|\ge S_\varepsilon.
\]
Consequently, for $s\ge S_\varepsilon$,
\[
e^{-4\varepsilon}
\frac{|\partial B_s^e\setminus E|}{|\partial B_s^e|}
\le
\frac{b(s)}{w(s)}=\frac{\int_{\partial B_s^e\setminus E}e^{4u(x)}\,d\sigma_x}{|\partial B_s^e|e^{4D(s)}}
\le
e^{4\varepsilon}
\frac{|\partial B_s^e\setminus E|}{|\partial B_s^e|}.
\]
Remark \ref{FinitnfinitInters} implies that the area ratio
tends to $1$. Letting first $s\to\infty$ and then
$\varepsilon\to0$, we obtain
\begin{align}
    \frac{b(s)}{w(s)}\longrightarrow1.\label{20260912diyi}
\end{align}
Moreover, \eqref{GxDx} yields
\begin{align}
    \frac{z(s)}{w(s)}
=e^{4(G(s)-D(s))}\longrightarrow1. \label{20260912dier}
\end{align}

Set
\[
W(R):=\int_1^R w(s)\,ds
     =V_{\widetilde g}(B_R^e).
\]
By \eqref{20260904444}, there exists $S>1$ such that
$a(s)\le2w(s)$ for all $s\ge S$. Hence
\begin{align}
    V_g(B_R^e)
=\int_1^R a(s)\,ds&= \left(\int_1^S+\int_S^R\right)a(s)\,ds \leq \int_1^S a(s)\,ds+2\int^R_S w(s)ds.\nonumber\\
&\le \int_1^S a(s)\,ds+2W(R),
\qquad(R\ge S).
\end{align}

The hypothesis $V_g(B_R^e)\to\infty$ therefore implies
\begin{equation}
W(R)=V_{\widetilde g}(B_R^e)\longrightarrow\infty.
\label{2026090411}
\end{equation}

We now pass from the spherical estimates to volume estimates
by direct integration. Let $h(s)$ be any one of $a(s)$, $b(s)$, and $z(s)$.
From (\ref{20260904444}), (\ref{20260912diyi}) and  (\ref{20260912dier}), we have $h(s)/w(s)\to1$. For every $\varepsilon>0$ there exists
$S_\varepsilon>1$ such that
\[
|h(s)-w(s)|\le\varepsilon w(s)
\qquad(s\ge S_\varepsilon).
\]
Thus, for $R\ge S_\varepsilon$,
\[
\begin{aligned}
\left|\int_1^R h(s)\,ds-W(R)\right|
&\le
\int_1^{S_\varepsilon}|h(s)-w(s)|\,ds
+\int_{S_\varepsilon}^R|h(s)-w(s)|\,ds\\
&\le C_{\varepsilon,h}+\varepsilon W(R),
\end{aligned}
\]
where $C_{\varepsilon,h}<\infty$ is independent of $R$.
Dividing by $W(R)$ and using \eqref{2026090411}, we obtain
\[
\limsup_{R\to\infty}
\left|\frac{\int_1^R h(s)\,ds}{W(R)}-1\right|
\le\varepsilon.
\]
Since $\varepsilon>0$ is arbitrary,
\[
\frac{\int_1^R h(s)\,ds}{W(R)}\longrightarrow1.
\]

Polar integration, justified by Tonelli's theorem, now gives
\begin{equation}
\frac{V_g(B_R^e)}{W(R)}\to1,\qquad
\frac{V_g(B_R^e\setminus E)}{W(R)}\to1,\qquad
\frac{V_{\hat g}(B_R^e)}{W(R)}\to1.
\label{IntegratedVolumeComparison}
\end{equation}
In particular,
\begin{equation}
V_{\hat g}(B_R^e)\longrightarrow\infty,
\label{20260904}
\end{equation}
and
\begin{equation}
\frac{V_g(B_R^e\setminus E)}{V_g(B_R^e)}
\longrightarrow1.
\label{20260904555}
\end{equation}
Subtracting the first two normalized volume limits yields
\begin{equation}
\frac{V_g(B_R^e\cap E)}{V_{\widetilde g}(B_R^e)}
=
\frac{V_g(B_R^e)}{W(R)}
-\frac{V_g(B_R^e\setminus E)}{W(R)}
\longrightarrow0.
\label{line:539}
\end{equation}
Finally,
\[
\frac{V_g(B_R^e\cap E)}{V_{\hat g}(B_R^e)}
=
\frac{V_g(B_R^e\cap E)/W(R)}
     {V_{\hat g}(B_R^e)/W(R)}
\longrightarrow0.
\]
This proves both assertions.
\end{proof}
We can now prove the following important Lemma.
\begin{lemma}\label{gptlemma}
   Let $M_i:=(\mathbb{R}^4\backslash B, g=e^{2u}|dx|^2)$ be a conformally flat end satisfying $R[g]\ge 0$ and $Q[g]\ge 0$, and assume that $M_i$ is an end of the complete manifold $(M,g)$. Then, as for $u(x)$ and $D(s)$, for any $T\geq 2R$, we have
    \begin{align}
    \sup_{\theta\in S^3}\left|
    \frac{\displaystyle\int_R^{T}e^{u(s\theta)}\,ds}{\displaystyle\int_R^{T}e^{D(s)}\,ds}-1\right|\longrightarrow 0,\qquad R\to\infty.\label{gpt:3}
    \end{align}

    Moreover, 
    \begin{align}
        \displaystyle\int_1^{T}e^{u(s\theta)}\,ds=(1+o_T(1))\displaystyle\int_1^{T}e^{D(s)}\,ds,\label{haoshiziizzi}
    \end{align}
    where $o_T(1)\to 0$, as $T\to \infty$ uniformly for $\theta\in S^3$. 
\end{lemma}
\begin{proof}
    From Lemma \ref{NormalMetric}, we know that $f:=2Q[g]e^{4u}\in L^1(\Omega)$, where $\Omega:=\mathbb{R}^4\backslash B$. Let's simplify the notation here. We denote
    \begin{align}
        d\mu(y)=\frac{f(y)}{8\pi^2}\,dy, \qquad m_\mu:=\mu(\Omega)<\infty, \label{gpt:5}
    \end{align}
    and 
    \begin{align}
    w(s):=e^{D(s)}>0,\qquad\bar b(s):=\fint_{S^3}b(\omega/s)\,d\sigma(\omega).\label{gpt:6}
    \end{align}
To simplify notation, we omit the domain $\Omega$ when integrating over it; that is,
\begin{align*}
    \displaystyle   \int g(y)dy:=\int_{\Omega}g(y)dy.
\end{align*}

    Fix $y\in\Omega$ and set $t:=|y|$. We define
    \begin{align}
    h(s,t):=\fint_{S^3}\log|s\omega-y|\,d\sigma(\omega). \label{gpt:7}
    \end{align}
By spherical symmetry, the expression in (\ref{gpt:7}) depends only on $s$ and $t$ and is independent of the direction of $y$. Next, we will show that
    \begin{align}
        h(s,t)=
    \begin{cases}
    \displaystyle\log t+\frac{s^2}{4t^2},&0\le s\le t,\\[5pt]
    \displaystyle\log s+\frac{t^2}{4s^2},&s\ge t.
    \end{cases}\label{gpt:8}
    \end{align}
    
    First, we consider $s<t$. It is easy to compute that
    \begin{align}
        \Delta_x F(x)=\frac2{|x-y|^2},\label{gpt:9}
    \end{align}
    where $F(x)=\log|x-y|$.
    For example, the Laplacian of the radially symmetric function $\log r$ in dimension four is
    \begin{align}
        \frac{d^2}{dr^2}\log r+\frac3r\frac d{dr}\log r=-\frac1{r^2}+\frac3{r^2}=\frac2{r^2}.
    \end{align}
We could obtain (\ref{gpt:9}) by moving the origin to $y$.

Note that the function $v(x)=|x-y|^{-2}$ is harmonic in the domain $\{x\in\mathbb{R}^4:|x|<t\}$. Thus, using the mean value property of harmonic functions, we have 
\begin{align}
    \fint_{S^3}|s\omega-y|^{-2}\,d\sigma(\omega)=|y|^{-2}=t^{-2},\qquad s<t.\label{gpt:10}
\end{align}
Now, we take the spherical average of (\ref{gpt:9}). In polar coordinates, for the smooth function $F$, it follows that
\begin{align*}
    \Delta F=F_{ss}+\frac3sF_s+\frac1{s^2}\Delta_{S^3}F.
\end{align*}
Since $S^3$ has no boundary, the integral of the Laplacian on the sphere $S^3$ vanishes:
\begin{align*}
    \int_{S^3}\Delta_{S^3}F(s\omega)\,d\sigma(\omega)=0.
\end{align*}
Therefore, applying (\ref{gpt:10}), we know that
\begin{align}
    h_{ss}(s,t)+\frac3s h_s(s,t)=\frac2{t^2},\qquad0<s<t.\label{gpt:11}
\end{align}
Multiplying both sides of (\ref{gpt:11}) by $s^3$, we have
\begin{align*}
\frac d{ds}\bigl(s^3h_s(s,t)\bigr)=\frac{2s^3}{t^2}.
\end{align*}
Integrate the above equation over $[0,s]$. Since $F$ is smooth near the origin, $h_s$ is bounded near $s=0$, so $s^3 h_s\to 0$, as $s\to 0$. So, we have
\begin{align*}
    s^3h_s(s,t)=\int_0^s\frac{2r^3}{t^2}\,dr=\frac{s^4}{2t^2},
\end{align*}
that is, 
\begin{align*}
    h_s(s,t)=\frac{s}{2t^2}.
\end{align*}
Integrating $h_s$ from $0$ to $s$ and using $h(0,t)=\log t$ and the Newton--Leibniz formula, we obtain
\begin{align}
    h(s,t)=\log t+\int_0^s\frac r{2t^2}\,dr=\log t+\frac{s^2}{4t^2},
\label{gpt:12}
\end{align}
for $0\le s<t$.
Then, we consider $s>t$. We claim that the function $h$ has the property $h(s,t)=h(t,s)$. Since $h(s,t)$ depends only on the norm $|y|$ of $y$, we change the direction $y/|y|$ of $y$ to a fixed $e_1\in S^3$ with the norm $|y|$ of $y$ unchanged. Then, 
\begin{align}
    |s\omega-te_1|^2=s^2+t^2-2st\omega_1=|t\omega-se_1|^2.
\end{align}
So, we have $|s\omega-te_1|=|t\omega-se_1|$. Hence,
\begin{align}
    h(s,t)=\fint_{S^3}\log|s\omega-te_1|\,d\sigma(\omega)=\fint_{S^3}\log|t\omega-se_1|\,d\sigma(\omega)=h(t,s).\label{gpt:13}
\end{align}
In the case $s>t$, we can calculate the right hand side of (\ref{gpt:13}) from (\ref{gpt:12}), which is
\begin{align}
    h(s,t)=h(t,s)=\log s+\frac{t^2}{4s^2},\label{gpt:14}
\end{align}
where $s>t$.

Now, we consider the case of $s=t$. Here the logarithmic kernel only has one singular point. And the kernel is certainly integrable. In detail, let $\varphi$ denote the angle between $\omega\in S^3$ and $e_1$. We can write $\omega=(\cos\varphi,\sin \varphi\cdot \xi)$, $0\leq\varphi\leq \pi$, $\xi\in S^2$. So, we have  $d\sigma_{S^3}=\sin^2\varphi d\varphi d\sigma_{S^2}.$ Since $|S^2|=4\pi$ and $|S^3|=2\pi^2$, we have
\begin{align*}
    h(s,t)=\fint_{S^3}\log|s\omega-te_1|d\sigma(\omega)=\frac{2}{\pi}\int_{0}^{\pi}\log d_s(\varphi)\sin^2\varphi d\varphi,
\end{align*}
where
\begin{align*}
    d_s(\varphi)=\sqrt{(s-t)^2+2st(1-\cos(\varphi))}.
\end{align*}
When $s=t$, $d_t(\varphi)=\sqrt{2st(1-\cos(\varphi))}=2t\sin (\varphi/2)$. The only zero of $d_t$ is at $\varphi=0$. We simplify $\log d_t(\varphi)$ as follows:
\begin{align*}
    \log d_t(\varphi)=&\log 2t\sin (\varphi/2)=\log t +\log 2 +\log \sin(\varphi/2)+\log (\varphi)-\log (\varphi)\\
    =&\log (\varphi)+\log t +\log \frac{\sin (\varphi/2)}{\varphi/2}.
\end{align*}
    The last term above tends to $0$ as $\varphi\to 0$. Hence, 
    \begin{align*}
        \log d_t(\varphi)=\log \varphi+ \log t+o(1).
    \end{align*}
    Thus, 
    \begin{align*}
        h(t,t)=\frac{2}{\pi}\int_{0}^{\pi}\log d_t(\varphi)\sin^2\varphi d\varphi
    \end{align*}
    is integrable.
Then, we show that $h(s,t)$ is continuous near $s=t$. It suffices to consider $t/2\leq s\leq3t/2$ when taking the limit $s\to t$.
    By the triangle inequality,
    \begin{align*}
        d_s(\phi)=|s\omega-te_1|\leq s+t\leq \frac{5}{2}t.
    \end{align*}
    On the other hand, for $0\leq\varphi\leq\pi$, there holds $\sin(\varphi/2)\geq \varphi/\pi$. Combining this with $2st\geq t^2$, we obtain 
    \begin{align*}
        d_s(\varphi)^2=(s-t)^2+2st(1-\cos(\varphi))\geq   2t^2 \sin^2 (\varphi/2)\geq \frac{2t^2}{\pi^2}\varphi^2.
    \end{align*}
Thus, for all $0\leq\phi\leq \pi$, we have
\begin{align*}
    \frac{\sqrt{2}t}{\pi}\varphi\leq d_s(\varphi)\leq \frac{5}{2}t.
\end{align*}
Taking logarithms in the above bound, we have 
\begin{align*}
    \left|\log d_s(\varphi)\right|\leq C_t+|\log \varphi|,
\end{align*}
where $C_t$ depends on $t$ but not on $s$. We split $[0,\pi]=[0,1]\cup[1,\pi]$. For $\varphi\in [0,1]$, $\sin^2\phi\leq \phi^2$. Therefore,
\begin{align*}
    \int_{0}^{1}\sin^2 \varphi |\log \varphi|d\varphi\leq\int_{0}^{1}\varphi^2 |\log \varphi|d\varphi<\infty.
\end{align*}
As for $\varphi\in [1,\pi]$, the dominating function $C_t+|\log \varphi|$ is bounded. Hence $(C_t+|\log \varphi|)\sin^2 \varphi$ is integrable on $[0,\pi]$.

Note that for every $\varphi>0$, 
$\log d_s(\varphi)\to\log d_t(\varphi)$. The single point $\varphi=0$ does not affect the value of the integral. By the dominated convergence theorem,
\begin{align*}
    \lim_{s\to t}h(s,t)=&\lim_{s\to t}\frac{2}{\pi}\int_{0}^{\pi}\log d_s(\varphi)\sin ^2 \varphi d\varphi=\frac{2}{\pi}\int_{0}^{\pi}\log d_t(\varphi)\sin ^2 \varphi d\varphi\\
    =&\frac{2}{\pi}\int_{0}^{\pi}\log (2t \sin \frac{\varphi}{2})\sin ^2 \varphi d\varphi=h(t,t)
\end{align*}

Next, using $h(s,t)$, for $s<t$, we have 
\begin{align*}
    h(t,t)=\lim_{s\to t^-}\log t+\frac{s^2}{4t^2}=\log t+\frac{1}{4}.
\end{align*}

Now, we start to analyze $D(s)$. We first differentiate (\ref{gpt:8}) with respect to $s$ and multiply the result by $s$. We obtain that 
\begin{align}
    s h_s(s,t)=
\begin{cases}
\displaystyle\frac{s^2}{2t^2},&s\le t,\\[5pt]
\displaystyle1-\frac{t^2}{2s^2},&s\ge t.
\end{cases}\label{gpt:15}
\end{align}
Both expressions above in (\ref{gpt:15}) are equal to $1/2$ when $s=t$. And it is not hard to see that 
\begin{align}
    0\le s h_s(s,t)\le1.\label{gpt:16}
\end{align}
Next, taking the spherical average of (\ref{line:263}), we have 
\begin{align}
    D(s)=\int_\Omega\bigl(\log|y|-h(s,|y|)\bigr)d\mu(y)+\gamma\log s+\bar b(s).\label{gpt:17}
\end{align}
Here we explain the interchange of the order of integration to avoid hidden integrability issues. Fix \(s\). Since \(|y|\ge 1\),
\begin{align*}
    \log\frac{|y|}{|s\omega-y|}\ge-\log\left(1+\frac{s}{|y|}\right)\ge-\log(1+s).
\end{align*}
Therefore, after adding \(\log(1+s)\) to the kernel, we obtain a nonnegative function, so Tonelli's theorem can be used to interchange the order of integration. Then, by (\ref{gpt:8}), \(\log|y|-h(s,|y|)\) is a bounded function and is no greater than zero. Hence the integral after interchanging the order is finite; after subtracting the constant back, we obtain (\ref{gpt:17}). In particular, we did not assume that \(\int \log|y|\,d\mu\) is finite on its own.

By (\ref{gpt:16}), on any compact interval in \(s\), \(h_s\) is dominated by a constant independent of \(y\). Since \(\mu\) is also a finite measure, we can differentiate under the integral sign in (\ref{gpt:17}):
\begin{align}
    sD'(s)=-\int_\Omega s h_s(s,|y|)d\mu(y)+\gamma+s\bar b'(s).\label{gpt:18}
\end{align}
Choose a sufficiently large fixed \(s_0\) such that, for all \(s\ge s_0\), \(\omega/s\) lies in a fixed small closed ball on which \(b\) is smooth. Take \(L\ge 0\) such that \(|\nabla b|\le L\) on this small closed ball. By the chain rule,
\begin{align*}
    \bar b'(s)=\fint_{S^3}\nabla b(\omega/s)\cdot\left(-\frac\omega{s^2}\right)d\sigma,\qquad |s\bar b'(s)|\le\frac Ls\le L.
\end{align*}
Therefore, defining
\begin{align*}
    M:=m_\mu+|\gamma|+L,
\end{align*}
we have 
\begin{align}
    |sD'(s)|&=|-\int_\Omega s\cdot h_s(s,|y|)d\mu_y+\gamma+s\bar b'(s)|\nonumber\\
    &\leq \int_\Omega |s\cdot h_s(s,|y|)|d\mu_y+|\gamma| +|s\bar b'(s)|\nonumber\\
    &\leq m_\mu +|\gamma|+L= M,\qquad(s\ge s_0).\label{gpt:19}
\end{align}
For \(s\ge r\ge s_0\), integrating from \(r\) to \(s\) gives
\begin{align*}
|D(s)-D(r)|&\le\int_r^s|D'(t)|dt\nonumber\\
&\le M\int_r^s\frac{dt}{t}=M\log\frac sr.
\end{align*}
Exponentiating the above inequalities, we have
\begin{align}
    \left(\frac sr\right)^{-M}\le\frac{w(s)}{w(r)}\le\left(\frac sr\right)^M.\label{gpt:20}
\end{align}
In particular, for arbitrary $s,t\in[a,2a]$, $s/t\leq 2$. It follows that
\begin{align}
    \max_{s\in[a,2a]}w(s)\le2^M\min_{t\in[a,2a]}w(t).\label{gpt:21}
\end{align}
The purpose of this step is that later, when extracting \(\max w\) from the integral, we will not lose a coefficient that grows with \(a\).

In this step, we analyze the exact expression for \(u-D\) on \([a,2a]\). Now fix $a\geq s_0$ and take arbitrary  $s\in[a,2a]$,  $\theta\in S^3$.
Subtract (\ref{gpt:17}) from (\ref{line:263}). The \(\gamma\log s\) terms cancel completely, and the \(\log|y|\) terms also cancel within the same integral. Therefore,
\begin{align}
    u(s\theta)-D(s)=\int_\Omega q(s,\theta,y)d\mu(y)+b(\theta/s)-\bar b(s),\label{gpt:22}
\end{align}
where
\begin{align}
q(s,\theta,y):=h(s,|y|)-\log|s\theta-y|.\label{gpt:23}
\end{align}

Next, we estimate the two parts \(|y|\le a/2\) and \(|y|>a/2\) separately. The difference between them is that the first part remains at a certain distance from \(s\theta\) and has no singularity; the second part may contain \(y=s\theta\), so the logarithmic singularity must be retained.

On the one hand, we estimate the region \(|y|\le a/2\).
Let \(t:=|y|\). Since \(t\le a/2\le s/2\), by (\ref{gpt:8}), it follows that
\begin{align*}
h(s,t)=\log s+\frac{t^2}{4s^2}.
\end{align*}
Substituting into (\ref{gpt:23}), we obtain
\begin{align}
    q(s,\theta,y)&=h(s,|y|)-\log|s\theta-y|\nonumber\\
&=\log s+\frac{t^2}{4s^2}-\log |s\theta-y|\nonumber\\
&=\frac{t^2}{4s^2}-\log\frac{|s\theta-y|}{s}.
\label{gpt:24}
\end{align}
By the triangle inequality,
\begin{align*} 
1-\frac ts\le\frac{|s\theta-y|}{s}\le1+\frac ts.
\end{align*}
Let \(z=t/s\in[0,1/2]\). Direct integration gives
\begin{align*}   
\log(1+z)=\int_0^z\frac{d\xi}{1+\xi}\le z,\qquad-\log(1-z)=\int_0^z\frac{d\xi}{1-\xi}\le2z.
\end{align*}
Therefore, it follows that
\begin{align*}
    \log(1-\frac{t}{s})\leq \log \frac{|s\theta-y|}{s}\leq \log (1+\frac{t}{s}).
\end{align*}
Then, 
\begin{align*}
    \left| \log \frac{|s\theta-y|}{s} \right|&\leq \max\{-\log(1-\frac{t}{s}),\log(1+\frac{t}{s})\}\\
    &\le 2\frac{t}{s}.
\end{align*}
Using \(t/s\le 1/2\), \(s\in [a,2a]\)  we have
\begin{align}
|q(s,\theta,y)|
&\le\frac{t^2}{4s^2}+\frac{2t}{s}=(\frac{t}{4s}+2)\frac{t}{s}\nonumber\\
&\le\left(\frac18+2\right)\frac ts
\le3\frac ta.\label{gpt:25}
\end{align}
Therefore, the integral over the region \(|y|\le a/2\) satisfies
\begin{align}
    \left|\int_{|y|\le a/2}q(s,\theta,y)d\mu(y)\right|\le\frac3a\int_{|y|\le a/2}|y|d\mu(y).\label{gpt:26}
\end{align}
Notice that the right hand side of \ref{gpt:26} is independent of $s\in[a,2a]$ and $\theta$.

Now, we estimate the integration domain $\{|y|>a/2\}$.

We define a constant 
\begin{align*}
    C_0:=\log4+\frac{1}{4}.
\end{align*}
We will show that
\begin{align}
-C_0\le q(s,\theta,y)\le C_0+\log^+\frac a{|s\theta-y|},\quad |y|>a/2,\label{gpt:27}
\end{align}
where \(\log^{+} z = \max\{\log z, 0\}\).

We first show the lower bound in (\ref{gpt:27}). 
By (\ref{gpt:8}), \(h(s,t)\ge \log \max\{s,t\}\); meanwhile, 
\begin{align*}
|s\theta-y|\le s+t\le2\max\{s,t\}.
\end{align*}
Thus, it follows that
\begin{align}
q(s,\theta,y)&=h(s,t)-\log|s\theta-y|\nonumber\\
&\ge\log\max\{s,t\}-\log(s+t)\nonumber\\
&\ge\log\max\{s,t\}-\log 2\max\{s,t\}\nonumber\\
&=-\log2\ge-C_0.\label{gpt:28}
\end{align}
Here we obtain the lower bound of (\ref{gpt:27}).

As for the upper bound of (\ref{gpt:27}), we divide the region into two parts: $a/2 \leq t\leq 4a$ and $t>4a$. 

We firstly consider $a/2 \leq t\leq 4a$. Here \(s\le 2a\), \(t\le 4a\), and the fractional term in (\ref{gpt:8}) does not exceed \(1/4\) (regardless of which case in (\ref{gpt:8}) applies). Therefore,
\begin{align*}
    h(s,t)\le\log \max\{s,t\}+\frac{1}{4}
    \leq \log 4a +\frac{1}{4}.
\end{align*}
Therefore, we have 
\begin{align}
q(s,\theta,y)&=h(s,|y|)-\log |s\theta-y|\nonumber\\
&\le\log(4a)+\frac{1}{4}-\log|s\theta-y|\nonumber\\
&=C_0+\log\frac a{|s\theta-y|}\nonumber\\
&\le C_0+\log^+\frac a{|s\theta-y|}.\label{gpt:29}
\end{align}

Next, we consider $t>4a$. Here $s\leq 2a<t/2$, so we have 
\begin{align*}
|s\theta-y|\ge t-s>t/2.
\end{align*}
Using $h(s,t)\leq\log t+1/4$, we obtain that 
\begin{align}
    q(s,\theta,y)\le\log t+\frac14-\log(t/2)=\log2+\frac14\le C_0.\label{gpt:30}
\end{align}
Combining (\ref{gpt:28}), (\ref{gpt:29}), and (\ref{gpt:30}), we obtain the upper bound of  (\ref{gpt:27}).

Putting the above estimates together, we obtain a small error term and a nonnegative singular term.
Denote
\begin{align}
    I_a:=\frac1a\int_{|y|\le a/2}|y|\,d\mu(y),\qquad\tau_a:=\mu(\{|y|>a/2\}).\label{gpt:31}
\end{align}
Since \(\mu\) is finite, \(\tau_a\to 0\), as $a\to\infty$. Below, we show in an \(\varepsilon\)-style argument that \(I_a\to 0\).
Given \(\varepsilon>0\), first choose a fixed \(A\) sufficiently large so that \(\mu(\{|y|>A\})<\varepsilon\). When \(a>2A\),
\begin{align}
I_a&=\frac{1}{a}\int_{|y|\le A}|y|d\mu+\frac{1}{a}\int_{A<|y|\le a/2}|y|d\mu\nonumber\\
&\le\frac{A}{a} m_\mu+\frac{a/2}{a}\int_{A<|y|\le a/2}d\mu\nonumber\\
&\le\frac{A}{a} m_\mu+\frac{1}{2}\mu(\{|y|>A\})\nonumber\\
&\le\frac{A}{a} m_\mu+\frac\varepsilon2.\label{gpt:32}
\end{align}
Fixing \(A\) and letting \(a\to\infty\), we obtain \(\limsup_{a\to\infty} I_a\le \varepsilon/2\). Then letting \(\varepsilon\to 0\), we get \(I_a\to 0\). For the remainder term of \(b\), since \(|\nabla b|\le L\),
\begin{align*}
    |b(\theta/s)-b(0)|\le L/s,\qquad |\bar b(s)-b(0)|\le L/s.
\end{align*}
Adding them together, using $s\geq a$, we obtain that
\begin{align}
|b(\theta/s)-\bar b(s)|\le2L/s\le2L/a.\label{gpt:33}
\end{align}

Now we define a direction-independent error quantity.
\begin{align}
\eta_a:=3I_a+C_0\tau_a+\frac{2L}{a}\longrightarrow0,\label{gpt:34}
\end{align}
as $a\to \infty$, and the non-negative function
\begin{align}   
H_a(s,\theta):=\int_{|y|>a/2}\log^+\frac a{|s\theta-y|}\,d\mu(y)\ge0.\label{gpt:35}
\end{align}

Substituting (\ref{gpt:26}), (\ref{gpt:27}), and (\ref{gpt:33}) into (\ref{gpt:22}) and adding term by term, we obtain the lower bound
\begin{align*}
    u(s\theta)-D(s)=&\int_\Omega q(s,\theta,y)d\mu(y)+b(\theta/s)-\bar b(s)\\
    &\geq -\frac{3}{a}\int_{|y|\leq a/2}|y|d\mu_y - C_0 \int_{|y|> a/2} d\mu(y) - \frac{2L}{a}\\
    &\ge-3I_a-C_0\tau_a-2L/a=-\eta_a,
\end{align*}
and the upper bound
\begin{align*}
u(s\theta)-D(s)\le3I_a+C_0\tau_a+2L/a+H_a(s,\theta)=\eta_a+H_a(s,\theta).
\end{align*}
That is 
\begin{align}
    -\eta_a\le u(s\theta)-D(s)\le\eta_a+H_a(s,\theta)\quad\text{for all }s\in[a,2a],\ \theta\in S^3.\label{gpt:36}
\end{align}
At this point, all possibly large positive deviations have been absorbed into \(H_a\). We do not claim that \(H_a\) tends to zero uniformly pointwise; below we only prove that its contribution to the integral over the ray length is small. In the following, we will apply Jensen's inequality to show it.

Take \(a\) sufficiently large so that \(\tau_a\le 1/2\). If \(\tau_a=0\), then (\ref{gpt:35}) directly gives \(H_a=0\), and this part has no error. Hence we first discuss \(0<\tau_a\le 1/2\).

Define 
\begin{align}
    d\nu_a(y):=\frac{\mathbf1_{\{|y|>a/2\}}}{\tau_a}d\mu(y).\label{gpt:37}
\end{align}
It is a probability measure because
\[
\int d\nu_a=\frac{1}{\tau_a}\mu(\{|y|>a/2\})=1.
\]

Let
\[
\ell(s,\theta,y):=\log^+\frac{a}{|s\theta-y|}.
\]
Then \(H_a=\tau_a\int \ell\,d\nu_a\). Since the exponential function is convex, Jensen's inequality gives
\begin{align}
    e^{H_a(s,\theta)}
&=\exp\left(\int \tau_a\ell(s,\theta,y)\,d\nu_a(y)\right)\nonumber\\
&\le \int \exp\bigl(\tau_a\ell(s,\theta,y)\bigr)\,d\nu_a(y).\label{gpt:38}
\end{align}
Note that
\[
e^{\tau_a\log^+(a/d)}=\max\{1,a/d\}^{\tau_a},
\]
where \(d=|s\theta-y|\).
Using \(\int d\nu_a=1\) again, subtracting \(1\) from both sides of (\ref{gpt:38}) gives
\begin{align}
    e^{H_a(s,\theta)}-1
\le\int\left[
\max\left\{1,\frac{a}{|s\theta-y|}\right\}^{\tau_a}-1
\right]d\nu_a(y).\label{gpt:39}
\end{align}
Subtracting \(1\) here is important: what we want to control is the ``additional contribution'' relative to the original weight, not the integral of the entire exponential function.

If one is worried that the logarithmic kernel is unbounded, one can first replace \(\ell\) by \(\min\{\ell,N\}\), apply Jensen's inequality to the bounded function, and then let \(N\to\infty\), using monotone convergence to obtain the same conclusion.

We now estimate the contribution of $H_a$ to the length along the ray. Fix \(y\) and \(\theta\), and set
\[
c:=y\cdot\theta,\qquad y_\perp:=y-c\theta.
\]
Since \(|\theta|=1\) and \(y_\perp\cdot\theta=0\), we have
\begin{align}
    s\theta-y=(s-c)\theta-y_\perp,\qquad|s\theta-y|^2=(s-c)^2+|y_\perp|^2\ge(s-c)^2.\label{gpt:40}
\end{align}
Thus \(|s\theta-y|\ge |s-c|\). A smaller denominator makes the kernel larger, hence
\begin{align}
    \max\left\{1,\frac{a}{|s\theta-y|}\right\}^{\tau_a}-1\le\max\left\{1,\frac{a}{|s-c|}\right\}^{\tau_a}-1.\label{gpt:41}
\end{align}
Both sides of (\ref{gpt:41}) are non-negative. Therefore, after integrating over \(s\in[a,2a]\), we can extend the right-hand side to the entire real line:
\begin{align}
&\int_a^{2a}\left[
\max\left\{1,\frac{a}{|s\theta-y|}\right\}^{\tau_a}-1
\right]ds\nonumber\\
&\le\int_{\mathbb{R}}\left[
\max\left\{1,\frac{a}{|s-c|}\right\}^{\tau_a}-1
\right]ds. \label{gpt:42}
\end{align}
When \(|s-c|\ge a\), the integrand on the right-hand side equals \(1-1=0\). Hence the only interval that actually contributes is \((c-a,c+a)\). Let \(r=s-c\); then, using even symmetry, we obtain
\begin{align}
    \int_{\mathbb R}\left[\max\left\{1,\frac a{|s-c|}\right\}^{\tau_a}-1
\right]ds
&=\int_{-a}^a\left[\left(\frac a{|r|}\right)^{\tau_a}-1\right]dr\nonumber\\
&=2\int_0^a\left(a^{\tau_a}r^{-\tau_a}-1\right)dr\nonumber\\
&=2\left[a^{\tau_a}\frac{a^{1-\tau_a}}{1-\tau_a}-a\right]\nonumber\\
&=\frac{2a\tau_a}{1-\tau_a}.\label{gpt:43}
\end{align}
The integral above is finite, because \(\tau_a\le 1/2<1\). This is precisely where the vanishing of the tail mass comes into play.
Combining both (\ref{gpt:42}) and (\ref{gpt:43}), we have
\begin{align}
    \int_a^{2a}\left[\max\left\{1,\frac a{|s\theta-y|}\right\}^{\tau_a}-1\right]ds\le\frac{2a\tau_a}{1-\tau_a}.\label{gpt:44}
\end{align}
The right-hand side contains neither \(y\) nor \(\theta\). Therefore, this estimate covers all directions $\theta\in S^3$, even including the worst case where \(y\) lies exactly on the line.

With the preparation above, we could estimate the length $$\int_a^{2a}w(s)ds:=\int_a^{2a}e^{D(s)}ds.$$
To simplify the expression, denote 
\begin{align}
    W_a:=\int_a^{2a}w(s)ds,\qquad w_a^+:=\max_{s\in[a,2a]}w(s),\qquad w_a^-:=\min_{s\in[a,2a]}w(s).\label{gpt:45}
\end{align}
Since \(w\) is continuous and positive, these two extrema exist, and \(w_a^{-}>0\).
From (\ref{gpt:39}), extracting \(w_a^{+}\) again, we obtain
\begin{align}
&~~~~\int_a^{2a}w(s)(e^{H_a(s,\theta)}-1)ds\nonumber\\
&\leq \int_a^{2a}w(s) \left(e^{H_a(s,\theta)}-1\right)ds\nonumber\\
&\le\int_a^{2a}w(s)\int\left[
\max\left\{1,\frac{a}{|s\theta-y|}\right\}^{\tau_a}-1
\right]d\nu_a(y)ds\nonumber\\
&\le w_a^+\int_a^{2a}\int\left[\max\left\{1,\frac a{|s\theta-y|}\right\}^{\tau_a}-1\right]d\nu_a(y)ds.\label{gpt:46}
\end{align}
The integrand is nonnegative, so Tonelli's theorem can be used to interchange the order of integration in \(s\) and \(y\), without having to prove in advance that the double integral is finite:
\begin{align}
    \text{RHS of (\ref{gpt:46})}
&=w_a^+\int\left\{\int_a^{2a}
\left[\max\left\{1,\frac{a}{|s\theta-y|}\right\}^{\tau_a}-1\right]ds\right\}d\nu_a(y)\nonumber\\
&\le w_a^+\int\frac{2a\tau_a}{1-\tau_a}\,d\nu_a(y)\nonumber\\
&=w_a^+\frac{2a\tau_a}{1-\tau_a},\label{gpt:47}
\end{align}
where the second line uses (\ref{gpt:44}), and the last line uses the fact that \(\nu_a\) is a probability measure.

On the other hand, the length of \([a,2a]\) equals \(a\), so
\(
W_a\ge a w_a^-.
\)
Combining (\ref{gpt:21}), we have 
\begin{align}
    a w_a^+\le2^M a w_a^-\le2^M W_a. \label{gpt:48}
\end{align}
Again, since \(1/(1-\tau_a)\le 2\), (\ref{gpt:46}), (\ref{gpt:47}) and (\ref{gpt:48}) give
\begin{align}  
\int_a^{2a}w(s)(e^{H_a(s,\theta)}-1)ds&\le w_a^+\int_a^{2a}\int\left[\max\left\{1,\frac a{|s\theta-y|}\right\}^{\tau_a}-1\right]d\nu_a(y)ds\nonumber\\&\leq w_a^+\frac{2a\tau_a}{1-\tau_a}\leq 2^M W_a \frac{2\tau_a}{1-\tau_a}\nonumber\\
&\le2^{M+2}\tau_a W_a.\label{gpt:49}
\end{align}
When \(\tau_a=0\), \(H_a=0\), the left-hand side is zero, and (\ref{gpt:49}) also holds. 
Although $H_a$ may be very large at some points, the resulting increase in length is at most a constant times $\tau_a$ times the original radial length.

Here, we estimate the error
\begin{align}
\int_a^{2a}|e^{u(s\theta)}-e^{D(s)}|ds.\nonumber
\end{align} 
Denote \(
z(s,\theta):=u(s\theta)-D(s)
\). We start to estimate $|e^z-1|$.
By (\ref{gpt:36}), \(-\eta_a\le z\le \eta_a+H_a\).

If \(z\le 0\), then
\begin{align}
    |e^z-1|=1-e^z\le 1-e^{-\eta_a}.
\label{gpt:50}
\end{align}
If \(z\ge 0\), then
\begin{align}
|e^z-1|=e^z-1
&\le e^{\eta_a+H_a}-1\nonumber\\
&=(e^{\eta_a}-1)+e^{\eta_a}(e^{H_a}-1).\label{gpt:51}
\end{align}
Since the right-hand sides in both cases are non-negative, they can be added together to give a common upper bound:
\begin{align}
    |e^z-1|\le(1-e^{-\eta_a})+(e^{\eta_a}-1)+e^{\eta_a}(e^{H_a}-1).\label{gpt:52}
\end{align}
Take \(a\) larger so that \(0\le \eta_a\le 1\). By an elementary inequality,
\begin{align*}
    1-e^{-\eta_a}\le\eta_a,\qquad e^{\eta_a}-1\le e\eta_a,\qquad e^{\eta_a}\le e.
\end{align*}
Therefore, (\ref{gpt:52}) becomes \begin{align}
    |e^z-1|&\le 1-e^{-\eta_a}+e^{\eta_a}-1+e^{\eta_a}(e^{H_a}-1)\nonumber\\&\le(1+e)\eta_a+e(e^{H_a}-1).
\label{gpt:53}
\end{align}
Note that \(e^{u(s\theta)}=e^{D(s)}e^{u(s\theta)-D(s)}=w(s)e^{z(s,\theta)}\). Multiplying (\ref{gpt:53}) by \(w(s)\), integrating over \([a,2a]\), and using (\ref{gpt:49}), we obtain
\begin{align}
\int_a^{2a}|e^{u(s\theta)}-w(s)|ds
&=\int_a^{2a}w(s)|e^{z(s,\theta)}-1|ds\nonumber\\
&\le(1+e)\eta_a W_a
+e\int_a^{2a}w(s)(e^{H_a}-1)ds\nonumber\\
&\le\bigl[(1+e)\eta_a+e2^{M+2}\tau_a\bigr]W_a.\label{gpt:54}
\end{align}
In fact, defining the fixed constant \(C_1=1+e+e2^{M+2}\), we can rewrite (\ref{gpt:54}) as
\begin{align}
\int_a^{2a}|e^{u(s\theta)}-e^{D(s)}|ds\le C_1(\eta_a+\tau_a)\int_a^{2a}e^{D(s)}ds.\label{gpt:55}
\end{align}
Since \(\eta_a\) and \(\tau_a\) are independent of \(\theta\), (\ref{gpt:55}) holds simultaneously for all \(\theta\). This completes the proof for a single interval $[a,2a]$.

Since \(\eta_a+\tau_a\to 0\), define
\begin{align}
    \delta_R:=\sup_{a\ge R}(\eta_a+\tau_a).\label{gpt:56}
\end{align}
Then \(\delta_R\to 0\). This follows directly from the definition of the limit: given \(\varepsilon>0\), there exists \(A_\varepsilon\) such that every \(a\ge A_\varepsilon\) satisfies \(\eta_a+\tau_a<\varepsilon\). Hence, when \(R\ge A_\varepsilon\), this bound holds for all \(a\ge R\), and taking the supremum still gives \(\delta_R\le\varepsilon\).
Thus, for any \(a\ge R\), (\ref{gpt:55}) gives the uniform version
\begin{align}
    \int_a^{2a}|e^{u(s\theta)}-w(s)|\,ds
\le C_1\delta_R\int_a^{2a}w(s)\,ds.\label{gpt:57}
\end{align}
Next, regardless of how many doubling intervals like $[a,2a]$ \([R,T]\) contains, all intervals use the same \(\delta_R\).
Given \(T\ge 2R\), take the unique integer \(n\ge 1\) such that
\begin{align}
    2^nR\le T<2^{n+1}R.
\label{gpt:58}
\end{align}
Let \(a=2^n R\). Then \([R,T]\) is partitioned into
\begin{align}
    [R,2R],\ [2R,4R],\ldots,
[2^{n-1}R,2^nR],\ [a,T].\label{gpt:59}
\end{align}
If \(T=a\), the last segment has length zero. Overlap at the endpoints does not affect the one-dimensional integral.
For the first \(n\) complete intervals, replace \(a\) in (\ref{gpt:57}) by \(2^j R\), respectively, and then sum:
\begin{align}
\int_R^a|e^{u(s\theta)}-w(s)|ds&=\sum_{j=0}^{n-1}\int_{2^jR}^{2^{j+1}R}|e^{u(s\theta)}-w(s)|ds\nonumber\\
&\le C_1\delta_R\sum_{j=0}^{n-1}\int_{2^jR}^{2^{j+1}R}w(s)ds\nonumber\\
&=C_1\delta_R\int_R^a w(s)ds.\label{gpt:60}
\end{align}

The last segment is \([a,T]\subset[a,2a]\). Since the absolute value of the integrand is nonnegative,
\begin{align}
\int_a^T |e^{u(s\theta)}-w(s)|\,ds
&\le \int_a^{2a} |e^{u(s\theta)}-w(s)|\,ds \nonumber\\
&\le C_1\delta_R \int_a^{2a} w(s)\,ds. \label{gpt:61}
\end{align}
The problem is that the right-hand side integrates over \([T,2a]\) as well. We cannot directly say that it is less than \(\int_R^T w\). Next, we use the weight comparison between adjacent intervals.

For any \(s\in[a,2a]\) and \(r\in[a/2,a]\), we have \(1\le s/r\le 4\). By (\ref{gpt:20}),
\[
w(s)\le 4^M w(r).
\]
Therefore,
\begin{align}
\max_{s\in[a,2a]}w(s)\le 4^M\min_{r\in [a/2,a]}w(r). \label{gpt:62}
\end{align}
The interval \([a,2a]\) has length \(a\), while \([a/2,a]\) has length \(a/2\). Therefore,
\begin{align}
\int_a^{2a} w(s)\,ds
&\le a\max_{r\in[a,2a]} w(r) \nonumber\\
&\le a4^M\min_{s\in[a/2,a]} w(s) \nonumber\\
&=2\cdot4^M\left(\frac a2\min_{s\in[a/2,a]} w(s)\right) \nonumber\\
&\le 2\cdot4^M\int_{a/2}^a w(s)\,ds. \label{gpt:63}
\end{align}

Since \(n\ge 1\), \(a/2=2^{n-1}R\ge R\); and since \(a\le T\), we have
\[
[a/2,a]\subset[R,T].
\]
Thus (\ref{gpt:63}) further gives
\begin{align}
\int_a^{2a} w(s)\,ds \le 2\cdot 4^M \int_R^T w(s)\,ds. \label{gpt:64}
\end{align}
Substituting (\ref{gpt:64}) back into (\ref{gpt:61}), we obtain the estimate for the last segment
\begin{align}
\int_a^T |e^{u(s\theta)}-w(s)|\,ds
&\le C_1 \delta_R \int_a^{2a}w(s)ds \nonumber\\ &\leq 2\cdot 4^M C_1\delta_R \int_R^T w(s)\,ds. \label{gpt:65}
\end{align}
This also explains why we require \(T\ge 2R\) here: it guarantees that before the last segment there is at least one complete adjacent interval, which can be used to control the error on the last segment. The present proof does not claim that the result holds for arbitrarily short intervals \([R,T]\).
Adding (\ref{gpt:60}) and (\ref{gpt:65}), and using \(\int_R^a w\le \int_R^T w\), we obtain
\begin{align}
\int_R^T |e^{u(s\theta)}-w(s)|\,ds
&=\int_R^a |e^{u(s\theta)}-w(s)|\,ds
+\int_a^T |e^{u(s\theta)}-w(s)|\,ds \nonumber\\
&\le C_1\delta_R\int_R^a w(s)\,ds
+2\cdot4^M C_1\delta_R\int_R^T w(s)\,ds \nonumber\\
&\le C_2\delta_R\int_R^T w(s)\,ds, \label{gpt:66}
\end{align}
where $C_2=C_1(1+2\cdot4^M)$ is a fixed constant independent of $R$, $T$, and $\theta$, and $\delta_R\to0$ as $R\to\infty$.

Since \(w>0\) and \(T>R\), the denominator \(\int_R^T w\) is strictly positive. Hence
\begin{align}
\left|\frac{\int_R^T e^{u(s\theta)}\,ds}{\int_R^T e^{D(s)}\,ds}-1\right|
&=\frac{\left|\int_R^T(e^{u(s\theta)}-e^{D(s)})\,ds\right|}
{\int_R^T e^{D(s)}\,ds} \nonumber\\
&\le\frac{\int_R^T|e^{u(s\theta)}-e^{D(s)}|\,ds}
{\int_R^T e^{D(s)}\,ds} \nonumber\\
&\le C_2\delta_R. \label{gpt:67}
\end{align}

Since the right-hand side of (\ref{gpt:67}) does not contain \(\theta\), we can now take the supremum:
\begin{align}
\sup_{\theta\in S^3}
\left|\frac{\int_R^T e^{u(s\theta)}\,ds}{\int_R^T e^{D(s)}\,ds}-1\right|
\le C_2\delta_R\longrightarrow0. \label{gpt:68}
\end{align}

Moreover, since the right-hand side also does not contain \(T\), this estimate holds uniformly for all \(T\ge 2R\).

Finally, we show (\ref{haoshiziizzi}). For fixed $R$ large enough, we define  
\begin{align}
    C_R:=\sup_{\theta\in S^3}\int_1^R|e^{u(s\theta)}-e^{D(s)}|ds<\infty.
\end{align} 
From (\ref{gpt:66}), it follows that
\begin{align}
    \sup_{\theta\in S^3}\left|\frac{\displaystyle\int_1^Te^{u(s\theta)}ds}{\displaystyle\int_1^T e^{D(s)}ds}-1\right|&\leq \frac{\displaystyle\sup_{\theta\in S^3}\left(\displaystyle\int_1^R+\displaystyle\int_R^T\right)|e^{u(s\theta)}-e^{D(s)}|ds}{\displaystyle\int_1^T e^{D(s)}ds}\nonumber\\
    &\leq \frac{C_R}{\displaystyle\int_1^T e^{D(s)}ds}+\sup_{\theta\in S^3}\frac{\displaystyle\int_R^T|e^{u(s\theta)}-e^{D(s)}|ds}{\displaystyle\int_R^T e^{D(s)}ds}\nonumber\\
    &\leq \frac{C_R}{\displaystyle\int_1^T e^{D(s)}ds} +C_2\delta_R,
\end{align}
where $\delta_R \to 0$, as $R\to\infty$. We choose $R=R_\varepsilon$ large enough such that $C_2 \delta_R<\varepsilon/2$. Then we take $T\geq 2R_\varepsilon$ large enough such that $C_R/\int_1^T e^{D(s)}ds <\varepsilon/2$, which means (\ref{haoshiziizzi}) is valid. 
\end{proof}

We now focus on distances under different metrics. For any $p, q \in (\mathbb{R}^4\backslash B, g)\subseteq (M,g)$, where $(M,g)$ is a complete manifold, we denote by $d_g(p,q)$ the distance between $p$ and $q$  under the metric $$g = e^{2u(x)}|dx|^2$$. We introduce the following notation for several quantities and sets that will be used frequently.  For every $z\in\mathbb{R}^4\backslash B$, define
\begin{align}
        l_z:=d_g&\left(z,\frac{z}{|z|}\right),~\widetilde{l_z}:=\int_{1}^{|z|}e^{D(s)}ds,~\hat{l_z}:=\int_{1}^{|z|}e^{G(s)}ds,\label{DefOfgeodesic}\\
        B_r:=&\left\{q\in\mathbb{R}^4\backslash B: d_g(q,\frac{q}{|q|})\leq r \right\},\label{def:B_r}\\
    \widetilde{B}_r:=&\left\{ q\in\mathbb{R}^4\backslash B: \widetilde{l_q}\leq r \right\}.\label{Bwanr}
\end{align}
Readers should distinguish the end-truncated Euclidean ball $B_r^e$, the cylinder $\widetilde{B}_r$ for the radial metric, the geodesic ball $B_r(p)$ centred at $p\in M$, and the set $B_r$ defined above. We call $B_r$ a radial geodesic cylinder.

First, we consider the relationship between $V_g(B_r^e)$ and $V_{\widetilde{g}}(B_r^e)$.

\begin{align}
    V_{\widetilde{g}}(B_r^e)&=\int_{ B_r^e\backslash B}\exp\left\{ 4\fint_{\partial B_{|x|}^e}u(y)d\sigma_y\right\}dx\nonumber\\
    &=\int_1^r |\partial B_s^e| \exp\left\{ 4\fint_{\partial B_s^e}u(y)d\sigma_y\right\}ds.
\end{align}
Note that $\frac{1}{|\partial B_s^e|}d\sigma_y$ is a probability measure and the function $s\mapsto e^s$ is convex. Applying Jensen's inequality, we have
\begin{align}
    V_{\widetilde{g}}(B_r^e)&=\int_1^r |\partial B_s^e| \exp\left\{ 4\fint_{\partial B_s^e}u(y)d\sigma_y\right\}ds\nonumber\\
    &\leq\int_1^r |\partial B_s^e|  \int_{\partial B_s^e}\exp \left\{4u(y)\right\}\frac{1}{|\partial B_s^e|}d\sigma_yds\nonumber\\
    &=\int_1^r \int_{\partial B_s^e} e^{4u(y)}d\sigma_y ds=V_g(B_r^e).\label{20260912ee}
\end{align}
\begin{lemma}\label{Kscunzai}
    There exists a continuous, strictly increasing function  $K:[1,\infty)\to [1,\infty)$ such that
    \begin{align}
        \frac{\int_{1}^R e^{D(s)}ds}{\int_{R}^{K(R)}e^{D(s)}ds}\to 0,\qquad R\to\infty.\nonumber
    \end{align}
\end{lemma}
\begin{proof}
    From Lemma \ref{wuxianchangxingzhi}, we know that 
 \begin{align*}
     \int_{1}^R e^{D(s)}ds\to \infty, 
 \end{align*} as $R\to \infty$.
 Let
\[
A(t):=\int_{1}^t e^{D(s)}\,ds.
\]
\(A\) is continuous and strictly increasing, and since \(A(1)=0\) and \(A(t)\to\infty\), it has a continuous inverse defined on \([0,\infty)\). For \(R\geq1\), take
\begin{align}
K(R):=\max\{2R,A^{-1}((R+1)A(R))\}. \label{gpt:69}
\end{align}
Both functions in the maximum are continuous and strictly increasing on $[1,\infty)$, so $K$ is continuous and strictly increasing, with $K(1)=2$. Moreover, $K(R)\geq 2R$, and hence $K(R)\to\infty$. Second, because \(A\) is increasing,
\[
A(K(R))\ge (R+1)A(R),
\]
hence
\[
\int_R^{K(R)} e^{D(s)}\,ds=A(K(R))-A(R)\ge R A(R).
\]
Therefore,
\begin{align}
\frac{\int_{1}^R e^{D(s)}\,ds}{\int_R^{K(R)} e^{D(s)}\,ds}
\le \frac{1}{R}\longrightarrow 0. \label{gpt:70}
\end{align}
Here \(K\) depends only on \(D\), not on \(\theta\in S^3\).
\end{proof}
\begin{lemma}\label{line:399lemmaGeodesic}
    Under the same assumptions and notation as in Lemma \ref{FiniteVolumeE}, fix $\varepsilon>0$. There exists $R_0\in\mathbb{R}^+$ such that, for any $z:=K(|x|)x/|x| \in\mathbb{R}^4\backslash B$ such that $|x|>R_0$, we have  $\widetilde{l_z}\leq e^{\varepsilon}\cdot l_{z}$, and $l_z\leq e^{\varepsilon}\cdot\widetilde{l_z}$, where $K(s)$ is defined in Lemma \ref{Kscunzai}.
\end{lemma}
\begin{proof}
Using Lemma \ref{gptlemma}, for $T>2R$ and sufficiently large $R$, we have 
\begin{align}
    \int_R^Te^{u(s \theta)}ds=(1+o_R(1))\int_R^Te^{D(s)}ds.\label{20260912111}
\end{align}
Note this $o_R(1)$ is independent of $\theta$. From Lemma \ref{Kscunzai}, there exists a function $K(s)$ such that 
    \begin{align}
        \frac{\displaystyle\int_{1}^R e^{D(s)}ds}{\displaystyle\int_{R}^{K(R)}e^{D(s)}ds}\to 0,\nonumber
    \end{align} as $R\to \infty$.
Then, using Lemma \ref{gptlemma}, for every $x\in \mathbb{R}^4\backslash B$, we have
\begin{align}
       \frac{\displaystyle\int_{1}^R e^{u(s\cdot x/|x|)}ds}{\displaystyle\int_{R}^{K(R)}e^{D(s)}ds}=\frac{(1+o_R(1))\displaystyle\int_{1}^R e^{D(s)}ds}{\displaystyle\int_{R}^{K(R)}e^{D(s)}ds}\to 0,\nonumber
\end{align} as $R\to\infty.$
Recall the property $G(r) = D(r) + o(1)$ in (\ref{line:269}).

Therefore, for any $\varepsilon > 0$, let $z:=K(|x|)x/|x|$. We will show that $l_z \le e^\varepsilon \widetilde{l}_z$ holds for sufficiently large $|x|$. Using the minimizing property of $l_z$ and setting $T:=K(|x|)$ and $R:=|x|$, we obtain
\begin{align}
    l_z\leq&\int_{1}^{K(|x|)}e^{u(s\frac{x}{|x|})}ds=\int_{1}^{|x|}+\int_{|x|}^{K(|x|)}e^{u(s\frac{x}{|x|})}ds\nonumber\\
    =&\int_1^{|x|}{e^{u(s\frac{x}{|x|})}ds}+(1+o_{|x|}(1))\int_{|x|}^{K(|x|)}{e^{D(s)}}ds\nonumber\\
    =&\left(\frac{\int_1^{|x|}{e^{u(s\frac{x}{|x|})}ds}}{\int_{|x|}^{K(|x|)}e^{D(s)}ds}+1+o_{|x|}(1)\right)\int_{|x|}^{K(|x|)}e^{D(s)}ds\nonumber\\
    \leq&(1+o_{|x|}(1))\int_{|x|}^{K(|x|)}e^{D(s)}ds \nonumber\\
     \leq& (1+\tilde o_{|x|}(1))\int_{1}^{K(|x|)}e^{D(s)}ds\leq e^{\varepsilon}\widetilde{l_z}.\label{zhongjijielun}
\end{align}
where we choose $|x|$ sufficiently large so that $1+\tilde o_{|x|}(1)\leq e^\varepsilon$, which gives the last inequality.

    Conversely, for any $\varepsilon>0$, we show that  $\widetilde{l}_z\leq e^\varepsilon l_z$ for sufficiently large $|x|$. For any $z:=K(|x|)\cdot x/|x|\in\mathbb{R}^4\backslash B$ with $|x|$ sufficiently large, there exists a curve $\gamma$ on $M$ connecting $z/|z|=x/|x|$ to $z=K(|x|)\cdot x/|x|$ such that its $g$-length equals $l_z := d_g(K(|x|)x/|x|,x/|x|)$. It is permissible for $\gamma(s)$ to leave $\mathbb{R}^4\backslash B$. If this happens, there exists the last time $t_0$ such that $\gamma(t_0) \in \partial B$ and thereafter $\gamma$ stays entirely in the end $(\mathbb{R}^4\backslash B, g)$. Denote the Euclidean length of this part by $L_{g_e}(\eta)$, and parameterize this part of the curve by Euclidean arc length as $\eta(s):[0,L_{g_e}(\eta)]\to \mathbb{R}^4\backslash B$. Set $\rho(s):=|\eta(s)|$. Then $|\rho'|\leq 1$, and
    \begin{align}
        \hat{l_z}=\int_{1}^{K(|x|)}e^{G(s)}ds\leq&\int_{0}^{L_{g_e}(\eta)}e^{G(|\eta(s)|)}|\rho'(s)|ds\leq\int_{0}^{L_{g_e}(\eta)}e^{u(\eta(s))}ds\nonumber\\\leq&L_g(\eta)\leq
        L_g(\gamma) = l_z. \label{line:406}
    \end{align}
    If $\gamma$ stays entirely in the end, set $\eta=\gamma$. The same argument shows that \ref{line:406} remains valid.\par
    Finally, we compare $\widetilde{l_z}$, $\hat l_z$, and $l_z$. By (\ref{GxDx}), for any $\varepsilon>0$, we can take the same $R_0>1$ as above sufficiently large so that $D(s)\leq G(s)+0.5\varepsilon$ for $s\in[R_0,K(|x|)]$. Choosing $\vec{n_1}\in S^3_0$ and using (\ref{wanbeiDG}), we obtain $D(s)\leq u(s\vec{n_1})+0.5\varepsilon$ for $s\in[R_0,K(|x|)]$, provided that $|x|$ is sufficiently large. Thus, we have
    \begin{align}
        \widetilde{l_z}=\int^{K(|x|)}_{1}e^{D(s)}ds&= 
        \int_{1}^{R_0}e^{D(s)}ds+\int^{K(|x|)}_{R_0}e^{D(s)}ds\nonumber\\
        &\leq C_{R_0}+\int^{K(|x|)}_{R_0}e^{G(s)+0.5\varepsilon}ds\nonumber\\
        &\leq C_{R_0}+\int^{K(|x|)}_{1}e^{G(s)+0.5\varepsilon}ds\nonumber\\
        &\leq e^{\varepsilon}\int^{K(|x|)}_{1}e^{G(s)}ds=e^{\varepsilon}\hat{l_z},\label{NewLine659}
    \end{align}
    where the last inequality follows by choosing $|x|>R_1$ sufficiently large that
    \[
        C_{R_0}\leq\bigl(e^{\varepsilon}-e^{\varepsilon/2}\bigr)
        \int_1^{K(|x|)}e^{G(s)}\,ds.
    \]
    Such a choice is possible by the divergence in (\ref{wanbeiDG}). 
    Combining it with (\ref{line:406}), naturally, we have
    \begin{align}
        \widetilde{l_z}\leq e^{\varepsilon}\hat{l_z}\leq e^{\varepsilon}{l_z},\label{line:410}
    \end{align}
    as $|x|>\max\{R_0,R_1\}$.\par
In summary, in view of the results (\ref{line:410}) and (\ref{zhongjijielun}), the lemma holds.
\end{proof}

\begin{theorem}
    Under the same assumptions as in Lemma \ref{FiniteVolumeE}, the $g$-volume of the Euclidean ball $B^e_{R}$ satisfies
\begin{align}
V_g(B_{R}^e) = (C_0+o(1))\bigl( \widetilde{l_1^R} \bigr)^4, \label{line:380ResultEuclid}
\end{align}
where $\widetilde{l_1^R}$ is defined in (\ref{tildeldingyi}), $o(1)\to 0$, as $R\to\infty$,   
    \begin{align*}
         C_0:=\omega_3\frac{(1+m)^3}{4}, 
    \end{align*}
    and
    $$1+m:=\frac{1}{4\pi^2}\int_{\partial B}T_+[g]e^{3u(y)}d\sigma_y-\frac{1}{4\pi^2}\int_{\mathbb{R}^4\backslash B}Q_{g}(y)e^{4u(y)}dy.$$
\end{theorem}
\begin{proof}
Applying Lemma \ref{MainLemma} directly yields
\begin{align}
\lim_{r\to\infty} r D'(r) := m, \label{line:286}
\end{align}
where $m$ is related to the integral of the $Q$-curvature and the boundary invariant $T$ defined in Lemma \ref{MainLemma}.
Set $\ell(R):=\int_1^R e^{D(r)}\,dr$. By Lemma \ref{wuxianchangxingzhi}, $\ell(R)\to\infty$. This implies $m\geq-1$: if $m<-1$, choose $\varepsilon>0$ with $m+\varepsilon<-1$. For all sufficiently large $r$, (\ref{line:286}) gives $D'(r)\leq(m+\varepsilon)/r$, and integration yields $e^{D(r)}\leq Cr^{m+\varepsilon}$. This would make $\int_1^\infty e^{D(r)}\,dr$ finite, a contradiction.

We first treat the finite-volume case. Since $D$ is the spherical average of $u$, Jensen's inequality, as in (\ref{20260912ee}), gives
\[
    \omega_3\int_1^\infty r^3e^{4D(r)}\,dr
    \leq V_g(\mathbb{R}^4\backslash B)<\infty.
\]
If $m>-1$, choose $0<\varepsilon<m+1$. For $r\geq R_0$ sufficiently large, $D'(r)\geq(m-\varepsilon)/r$, so $e^{D(r)}\geq Cr^{m-\varepsilon}$. Consequently,
\[
    \int_1^\infty r^3e^{4D(r)}\,dr
    \geq C'\int_{R_0}^\infty r^{3+4m-4\varepsilon}\,dr=\infty,
\]
because $3+4m-4\varepsilon>-1$. This contradiction shows that $m=-1$ and hence $C_0=\omega_3(1+m)^3/4=0$. The finite total volume and $\ell(R)\to\infty$ now imply
\begin{align}
    \lim_{R\to\infty}\frac{V_g(B_R^e)}{\ell(R)^4}=0=C_0.
    \label{20260912JXwei0}
\end{align}
This proves (\ref{line:380ResultEuclid}) in the finite-volume case.

For the rest of the proof, assume $V_g(B_R^e)\to\infty$ as $R\to\infty$.
    From (\ref{line:269}), we know that for any $\varepsilon>0$ and for $x\notin E$ with $R := |x|$ sufficiently large, there holds
    \begin{align}
        D(R)-\varepsilon\leq& u(x)=D(R)+o(1)
        \leq D(R)+\varepsilon.\label{line:291GandULookSame}
    \end{align}
    Formula (\ref{line:291GandULookSame}) shows that $u(x)$ is close to a radially symmetric function $D(|x|)$ outside the set $E$ near infinity. Now, we define a new radially symmetric metric $\widetilde{g}:=e^{2D(|x|)}|dx|^2$.     Denote\begin{align}
        \widetilde{{l}_{a}^{b}}:=\int_{a}^{b}e^{D(s)}ds,~(b\geq a\geq 1).\label{tildeldingyi}
    \end{align}
    \par
We next estimate the contribution of $E$ under the radial metric
$\widetilde g$ on the end-truncated balls $B_R^e$.
Set
\[
q(s):=\frac{|\partial B_s^e\cap E|}{|\partial B_s^e|},
\qquad
w(s):=|\partial B_s^e|e^{4D(s)}.
\]
By Remark \ref{FinitnfinitInters}, $q(s)\to0$ as $s\to\infty$.
Since $\widetilde g$ is radial, polar integration gives
\[
V_{\widetilde g}(B_R^e\cap E)=\int_1^R |\partial B_s^e\cap E|e^{4D(s)}ds
=\int_1^R q(s)w(s)\,ds.
\]
Fix $\varepsilon>0$ and choose $S_\varepsilon>1$ such that
$0\le q(s)\le\varepsilon$ for all $s\ge S_\varepsilon$.
Then, for $R\ge S_\varepsilon$,
\begin{align}
    0\le V_{\widetilde g}(B_R^e\cap E)=&\left(\int_1^{S_\varepsilon}+\int_{S_\varepsilon}^R\right) q(s)w(s)ds\leq C_\varepsilon +\varepsilon \int_{S_\varepsilon}^R w(s)ds \nonumber\\
    \le& C_\varepsilon+\varepsilon \int_{1}^R w(s)ds = C_\varepsilon
+\varepsilon V_{\widetilde g}(B_R^e),\label{202609121212}
\end{align}
where
\[
C_\varepsilon:=\int_1^{S_\varepsilon}q(s)w(s)\,ds<\infty.
\]
In the present infinite-volume case where $V_g(B_R^e)\to\infty$, we know that  
$V_{\widetilde g}(B_R^e)\to\infty$ by the proof of
Lemma \ref{FiniteVolumeE}.
Dividing (\ref{202609121212}) by $V_{\widetilde g}(B_R^e)$, letting $R\to\infty$, and then
letting $\varepsilon\to0$, we obtain
\[
\frac{V_{\widetilde g}(B_R^e\cap E)}
     {V_{\widetilde g}(B_R^e)}
\longrightarrow0.
\]
Consequently,
\begin{equation}
\frac{V_{\widetilde g}(B_R^e\setminus E)}
     {V_{\widetilde g}(B_R^e)}=1-\frac{V_{\widetilde g}(B_R^e\cap E)}
     {V_{\widetilde g}(B_R^e)}
\longrightarrow1.
\label{haodegainian}
\end{equation}
as $R\to \infty$.

The integrated estimate (\ref{IntegratedVolumeComparison}), proved by separating a fixed compact part from the region where the spherical errors are uniformly small, gives
\[
    \frac{V_g(B_R^e)}{V_{\widetilde g}(B_R^e)}\longrightarrow1.
\]
Indeed, the fixed compact part has finite volume, whereas $V_{\widetilde g}(B_R^e)\to\infty$ by (\ref{2026090411}). Thus no uniform smallness of the pointwise error on the entire expanding ball is required. It follows that
\begin{align}
    \frac{V_g(B_R^e)}{\left(\widetilde{l_1^R}\right)^4}
    &=(1+o(1))\frac{V_{\widetilde g}(B_R^e)}{\left(\widetilde{l_1^R}\right)^4}
    =e^{o(1)}\frac{V_{\widetilde g}(B_R^e)}{\left(\widetilde{l_1^R}\right)^4}.
    \label{line:306}
\end{align}

We use the version of L'Hôpital's rule in which the denominator tends to $+\infty$, its derivative is positive, and the derivative ratio has a limit; the numerator need not tend to infinity. First, (\ref{line:286}) and $\ell(R)\to\infty$ give
\begin{align}
\lim_{R\to\infty}\frac{Re^{D(R)}}{\int_1^Re^{D(t)}\,dt}
    &=\lim_{R\to\infty}\frac{e^{D(R)}+Re^{D(R)}D'(R)}{e^{D(R)}}\nonumber\\
    &=1+m. \label{line:320}
\end{align}
No prior divergence of $Re^{D(R)}$ is needed. Applying the same rule to the radial volume ratio, whose denominator $\ell(R)^4$ tends to infinity and has positive derivative, yields
\begin{align}
    \lim_{R\to\infty}\frac{V_{\widetilde g}(B_R^e)}{\ell(R)^4}
    &=\lim_{R\to\infty}
    \frac{\omega_3\int_1^R t^3e^{4D(t)}\,dt}{\ell(R)^4}\nonumber\\
    &=\frac{\omega_3}{4}\lim_{R\to\infty}
    \left(\frac{Re^{D(R)}}{\ell(R)}\right)^3.\label{line:314}
\end{align}
Substituting (\ref{line:320}) into (\ref{line:314}) yields
\begin{align}
    \lim_{R\to\infty}\frac{V_{\widetilde{g}}(B_{R}^e)}{\left(\widetilde{l_{1}^{R}}\right)^{4}}=\frac{\omega_3(1+m)^3}{4}:=C_0.\label{line:324}
\end{align}
Hence, by (\ref{line:324}) and (\ref{line:306}), we conclude that for sufficiently large $R:=|x|$,
\begin{align}
V_g(B_R^e) = (C_0+o(1))  \bigl( \widetilde{l_{1}^{R}} \bigr)^4, \label{line:380ResultEuclid111}
\end{align}
which completes the proof. 
\end{proof} \par

\begin{theorem}\label{MainTheorem}
    Suppose that the conformally flat simple end $(\mathbb{R}^4\backslash B,g=e^{2u}|dx|^2)$ is a submanifold of a complete manifold $(M,g)$ satisfying $R[g]\geq0$ and $Q[g]\geq0$ near infinity. Then, for the radial geodesic cylinder $B_{r}$ defined in (\ref{def:B_r}), we have
    $$
        \lim_{r\to \infty} \frac{V_g(B_{r})}{r^4}=\omega_3\frac{(1+m)^3}{4}:=C_0,    
    $$
    where 
    $$1+m:=\frac{1}{4\pi^2}\int_{\partial B}T_+[g]e^{3u(y)}d\sigma_y-\frac{1}{4\pi^2}\int_{\mathbb{R}^4\backslash B}Q_{g}(y)e^{4u(y)}dy.$$
    Here $T_+[g]$ is the boundary curvature computed with
respect to the unit normal $N_+[g]=e^{-u}\partial_r$,
which points toward infinity along the end. 
\end{theorem}

\begin{proof}
We first prove the theorem when $R[g]\geq0$ and $Q[g]\geq0$ on the entire displayed end. The reduction from non-negativity only near infinity is given at the end of the proof.

If $\lim_{r\to\infty}V_g(B_r^e)$ is finite, then the corresponding limit for the radial geodesic cylinders is also finite: $$\lim_{r\to\infty}V_g(B_r)<\infty.$$ We have already shown that if $\lim_{r\to\infty}V_g(B_r)$ is finite, then $C_0=0$, which is consistent with 
\begin{align}
            \lim_{r\to \infty} \frac{V_g(B_{r})}{r^4}=0=C_0, 
\end{align}
when  $\lim_{r\to\infty}V_g(B_r)$ is finite.

In the following proof, we assume that $\lim_{r\to\infty}V_g(B_r^e)=\infty$.
Recall the definition (\ref{DefOfgeodesic}) of $l_z$ and $\widetilde{l_z}$. We denote $z:=K(|x|)x/|x|$. By Lemma \ref{line:399lemmaGeodesic}, for every $\varepsilon>0$, there exists $R_0>0$ such that if $|x|\geq R_0$,
\begin{align}
e^{-\varepsilon}\widetilde{l_z}\leq l_z\leq e^{\varepsilon}\widetilde{l_z}.\label{Cedixianchabuduo}
\end{align}
Since $K$ is continuous and strictly increasing, $K(R)\geq2R$, and $K(R)\to\infty$, its restriction to $[R_0,\infty)$ maps onto $[K(R_0),\infty)$. Thus every $y$ with $|y|\geq K(R_0)$ has the form $y=K(|x|)x/|x|$ for some $|x|\geq R_0$, and the comparison (\ref{Cedixianchabuduo}) applies to all such $y$.

Fix a direction $\theta_0\in S^3$. The function $t\mapsto l_{K(t)\theta_0}=d_g(K(t)\theta_0,\theta_0)$ is continuous. Moreover, $\widetilde l_{K(t)\theta_0}=\int_1^{K(t)}e^{D(s)}\,ds\to\infty$ by Lemma \ref{wuxianchangxingzhi}, so (\ref{Cedixianchabuduo}) implies $l_{K(t)\theta_0}\to\infty$. The intermediate value theorem therefore shows that $l_{K(t)\theta_0}$ takes every sufficiently large value, even when $t$ is restricted to any sufficiently large tail interval. This justifies the changes of the length parameter used below.

First, we claim that 
\begin{align}
    B_{e^{-2\varepsilon}l_z}\subseteq \widetilde{B}_{\widetilde{l_z}}.\label{BallDeBaoHan}
\end{align}
Taking $R_*:=K(R_0)+1$, since $\widetilde{l_z}\to \infty$, we can choose $|x|$ sufficiently large so that
\begin{align}
B^e_{R_*} \subseteq {\widetilde{B}}_{\widetilde{l_z}}, \label{Limianyougeoushiqiu}
\end{align}
because $\overline{B^e_{R_*}}$ is compact in the end and $\widetilde l_y$ is bounded on this closure, whereas $\widetilde l_z\to\infty$.

For any $y \in B_{e^{-2\varepsilon}l_z} \backslash B^e_{R_*}$, we have $l_y \le e^{-2\varepsilon} l_z$. Then, $\widetilde{l_y} \le e^{\varepsilon} l_y \le e^{-\varepsilon} l_z \le \widetilde{l_z}$. The last inequality follows from (\ref{zhongjijielun}). Hence, we have
\begin{align}
    B_{e^{-2\varepsilon}l_z} \backslash B^e_{R_*}\subseteq{\widetilde{B}}_{\widetilde{l_z}}.\label{1233}
\end{align}
This establishes the claim by formulas (\ref{1233}) and (\ref{Limianyougeoushiqiu}).\par

For the cylinder $\widetilde{B}_{\widetilde{l_z}}$ associated with the radial metric $\widetilde g=e^{2D(|x|)}|dx|^2$, the function $r\mapsto\int_1^r e^{D(s)}ds$ is strictly increasing. Hence, with $\widetilde R=K(|x|)$,
\begin{align}
    \overline{B^e_{\widetilde R}}\backslash N=\widetilde{B}_{\widetilde{l_z}}\backslash N,\label{ZhengHaoxiangdeng}
\end{align}
where $N$ is a fixed large compact subset of the end. The added outer sphere $\partial B^e_{\widetilde R}$ has zero four-dimensional volume for both $g$ and $\widetilde g$, so $B^e_{\widetilde R}\backslash N$ and $\widetilde{B}_{\widetilde{l_z}}\backslash N$ have the same volume for either metric. Moreover,
\begin{align}
    \widetilde{l_z}=\int_{1}^{K(|x|)}e^{D(|s|)}ds=\widetilde{l_1^{\widetilde{R}}}.\label{y1}
\end{align}
From (\ref{line:306}) we have
\begin{align}
\frac{V_g(B_R^e)}{(\widetilde{l_1^{R}})^4}=(1+o(1))\frac{V_{\widetilde{g}}(B_R^e)}{(\widetilde{l_1^{R}})^4}=e^{o(1)}\frac{V_{\widetilde{g}}(B_R^e)}{(\widetilde{l_1^{R}})^4},\label{jieyongyixia}
\end{align}
where $o(1)\to 0$ as $R\to\infty$.
Cancelling the denominators in (\ref{jieyongyixia}), we obtain
\begin{align}
V_g(B_R^e)=e^{o(1)}V_{\widetilde{g}}(B_R^e).\label{xiashiyaoyong}
\end{align}
Hence, in particular, for the $\widetilde{R}$ given in (\ref{ZhengHaoxiangdeng}) and using (\ref{xiashiyaoyong}), we have
\begin{align}
    V_g(B_{\widetilde{R}}^e\backslash N)&=e^{o(1)}V_{\widetilde{g}}(B_{\widetilde{R}}^e\backslash N)=e^{o(1)}(V_{\widetilde{g}}({{\widetilde{B}}_{\widetilde{l_z}}}\backslash N)).\label{Houmiandeshizigaigai}
\end{align}

Therefore, using (\ref{BallDeBaoHan}), (\ref{ZhengHaoxiangdeng}), (\ref{xiashiyaoyong}), $o(1) < \varepsilon$, (\ref{line:324}), the relation $\widetilde{l_z} = \widetilde{l_1^{\widetilde{R}}}$, and Lemma \ref{line:399lemmaGeodesic}, we obtain
\begin{align}
    V_g\left(B_{\frac{l_z}{e^{2\varepsilon}}}\backslash N\right)
    &\leq V_g\left(\widetilde{B}_{\widetilde{{l_z}}}\backslash N\right)=V_{g}\left(B^e_{\widetilde{R}}\backslash N\right)= e^{o(1)}V_{\widetilde{g}}(\widetilde{B}_{\widetilde{{l_z}}}\backslash N)\leq
    e^{\varepsilon}V_{\widetilde{g}}(\widetilde{B}_{\widetilde{{l_z}}})
    \nonumber\\
    &\leq e^\varepsilon e^\varepsilon (C_0+\varepsilon)(\widetilde{l_z})^4\leq e^{2\varepsilon}(C_0+\varepsilon)\cdot e^{4\varepsilon}(l_z)^4.\label{line80880}
\end{align} 
Replacing $l_z$ by $e^{2\varepsilon}l_z$, equation (\ref{line80880}) becomes 
\begin{align}
    V_g(B_{l_z}\backslash N)\leq e^{14\varepsilon}\cdot (C_0+\varepsilon) (l_z)^4.\label{xiaoyudengyu111}
\end{align}
Thus, 
\begin{align}
    V_g(B_{l_z})&=V_g(B_{l_z}\backslash N)+V_g(B_{l_z}\cap N)\leq e^{14\varepsilon}\cdot (C_0+\varepsilon) (l_z)^4+V_g(N)\nonumber\\
    &\leq e^{15\varepsilon}\cdot (C_0+\varepsilon) (l_z)^4\label{xiaoyudengyu}
\end{align}
\par
On the other hand, we prove the converse inclusion. For any $\varepsilon>0$, there exists an $R_0$ such that (\ref{Cedixianchabuduo}) is valid if $|x|>R_0$. Following an argument similar to the proof of (\ref{Limianyougeoushiqiu}), for $\widetilde{l_z}$ large enough, we have 
\begin{align}
    B^e_{R_*}\subseteq {\widetilde{B}}_{\widetilde{l_z}},\label{Limianyougeoushiqiu2}
\end{align}
where $R_*:=K(R_0)+1$.

We assert that the geodesic cylinder $\widetilde{B}_{\widetilde{l_z}}$ and $B_{l_z}$ satisfy the inclusion
\begin{align}
\widetilde{B}_{\widetilde{l_z}} \subseteq B_{e^{2\varepsilon}l_z}. \label{Baohanana}
\end{align}
For every $y\in\widetilde{B}_{\widetilde{l_z}}\backslash B^e_{R_*}$, we have $\widetilde{l_y}\leq\widetilde{l_z}$. Thus,
 $$l_y\leq e^{\varepsilon}\widetilde{l_y}\leq e^{\varepsilon} \widetilde{l_z}\leq e^{2\varepsilon} l_z.$$
 Thus, we obtain
 \begin{align}
     \widetilde{B}_{\widetilde{l_z}}\backslash B^e_{R_*}\subseteq B_{e^{2\varepsilon}l_z}.\label{x1}
 \end{align}
 For every $y_1\in B_{R_*}^e$, let $\displaystyle C:=\max_{y\in\overline{B^e_{R_*}}}l_y$. Since $l_z\to\infty$, for sufficiently large $|z|$ we have
 \begin{align} 
     l_{y_1}\leq C\leq e^{2\varepsilon} l_z.
 \end{align}
 Thus, we have
 \begin{align}
     B_{R_*}^e\subseteq B_{e^{2\varepsilon}l_z}.\label{x2}
 \end{align}
Thus, from (\ref{x1}) and (\ref{x2}), the assertion holds.
Now, for sufficiently large $|x|$, using (\ref{Baohanana}), Lemma \ref{line:399lemmaGeodesic}, (\ref{y1}), and (\ref{line:324}), we obtain
\begin{align}
    \frac{V_g(B_{e^{2\varepsilon}l_z}\backslash N)}{(e^{2\varepsilon}l_z)^4}
    &\geq\frac{V_g(\widetilde{B}_{\widetilde{l_z}}\backslash N)}{(e^{2\varepsilon}l_z)^4}
    \geq\frac{V_g(\widetilde{B}_{\widetilde{l_z}}\backslash N)}{(e^{3\varepsilon}\widetilde{l_z})^4}\nonumber\\
    &=\frac{V_g(B^e_{\widetilde{R}})-V_g(B^e_{\widetilde R}\cap N)}{(e^{3\varepsilon}\widetilde{l_1^{\widetilde{R}}})^4}\to e^{-12\varepsilon}C_0\nonumber\\
    &\geq C_0e^{-13\varepsilon}.\label{Newliine:567}
\end{align}
Replacing $l_z$ by $l_z/e^{2\varepsilon}$, we obtain 
\begin{align}
    V_g(B_{l_z})\geq V_g(B_{l_z}\backslash N)\geq \frac{C_0}{e^{13\varepsilon}}(l_z)^4.\label{xiajiexiajie}
\end{align}\par
Finally, combining (\ref{xiaoyudengyu}) with (\ref{xiajiexiajie}), we obtain
\begin{align}
    \lim_{|x|\to\infty}\frac{V_g(B_{l_z})}{(l_z)^4}=C_0=\omega_3\frac{(1+m)^3}{4},  \label{temme}
\end{align}
As shown above, $l_z$ takes every sufficiently large value along a fixed direction, with the corresponding parameter tending to infinity. Thus, setting $r=l_z$ in (\ref{temme}) proves the theorem when both curvatures are non-negative on the entire end.

We now remove this extra restriction. Write $E_1=\{x:|x|\geq1\}$ and choose $a>1$ so that $R[g],Q[g]\geq0$ on $E_a=\{x:|x|\geq a\}$. After the coordinate change $x=a\xi$, the conformal factor on $E_a$ is $u_a(\xi)=u(a\xi)+\log a$, so the case just proved applies to this end. In particular, $Q[g]$ is integrable on $E_a$ by Lemma \ref{NormalMetric}; it is integrable on $E_1\backslash E_a$ by smoothness and compactness, and hence is absolutely integrable on $E_1$.

For $z\in E_a$, define
\[
    l_z^{(a)}:=d_g\left(z,a\frac{z}{|z|}\right),\qquad
    B_r^{(a)}:=\{z\in E_a:l_z^{(a)}\leq r\},\qquad
    C_a:=\max_{\theta\in S^3}d_g(\theta,a\theta)<\infty.
\]
All distances here are measured in the ambient complete manifold. The triangle inequality gives $|l_z-l_z^{(a)}|\leq C_a$, and consequently, for $r>C_a$,
\[
    B_{r-C_a}^{(a)}\subseteq B_r\cap E_a\subseteq B_{r+C_a}^{(a)}.
\]
The omitted compact annulus has finite volume. Therefore the limit of $r^{-4}V_g(B_r)$ is the same as that of $r^{-4}V_g(B_r^{(a)})$.

It remains to check that the curvature expression for $1+m$ is also unchanged. Set $\Sigma_s=\{x:|x|=s\}$, and continue to compute $T_+$ with the normal pointing toward infinity. The compact annulus $A=\{1\leq|x|\leq a\}$ has Euler characteristic zero and is conformally flat. Its boundary spheres are umbilic, so the additional boundary invariant $L$ in the Chern--Gauss--Bonnet formula vanishes. The outward normals of $A$ are $-N_+$ at $\Sigma_1$ and $N_+$ at $\Sigma_a$. Thus the boundary formula \cite[(0.7)]{ChangQingYangInvention2000} yields
\[
    0=\int_A Q[g]\,dV[g]
      +\int_{\Sigma_a}T_+[g]\,d\sigma[g]
      -\int_{\Sigma_1}T_+[g]\,d\sigma[g].
\]
It follows that
\begin{align*}
    \int_{\Sigma_a}T_+[g]\,d\sigma[g]-\int_{E_a}Q[g]\,dV[g]
    =\int_{\Sigma_1}T_+[g]\,d\sigma[g]-\int_{E_1}Q[g]\,dV[g].
\end{align*}
The rescaling is only a change of coordinates and preserves these geometric integrals. Hence the constant obtained for $E_a$ is exactly the $C_0$ stated for the original end $E_1$. Together with the distance and volume comparison, this proves Theorem \ref{MainTheorem} under its stated near-infinity assumptions. 
\end{proof}

\begin{lemma}[\cite{MaQing2025LinearPotential}.Theorem 1.5]\label{FiniteEnd11}
    Suppose that $(\Omega,g=e^{2u}\bar g)$ is complete, that $\Omega\subseteq S^4$ is a subdomain, that $R[g]\geq0$ in $\Omega$, and that $Q[g]\geq0$ near $S^4\backslash \Omega$. Then $\Omega$ is conformally equivalent to  $S^4\backslash\{p_i\}_{i=1}^k$, where $\{p_i\}_{i=1}^k$ is a finite set of points on $S^4$.
\end{lemma}
\begin{proof}
    This is a direct corollary of Theorem 1.5 in Ma-Qing \cite{MaQing2025LinearPotential}. From the condition, we know that $(\Omega,g=e^{2u}\bar g)$ is complete. Then, applying this theorem with $D:=S^4$ and $S:=S^4\backslash\Omega$, we note that $Q[g]$ is non-negative near $S$. Consequently, $Q^-[g]$ is supported in a non-singular compact set and is therefore absolutely integrable on $D\backslash S$. Now, applying \cite[Theorem 1.5]{MaQing2025LinearPotential}, we obtain that $\Omega $ is conformally equivalent to $ S^4\backslash\{p_i\}_{i=1}^k$.
\end{proof}
As for the complete, simply connected, locally conformally flat manifold $(M^4,g)$, there exists an immersion $\Phi: M^4 \to S^4$ which is called the developing map such that the locally conformally flat structure is induced by $\Phi$. If $R[g]$ is non-negative, the results of Schoen--Yau \cite{SchoenYau1988conformally} and Chodosh--Li \cite[Corollary 4]{LiChaoChodosh2024DevelopingMapInjective} imply that the developing map $\Phi:M^4 \to S^4$ is injective. It means that we can regard such a kind of $M^4$ as a subdomain of $S^4$ in the conformal sense. 

Moreover, we can get the Chern-Gauss-Bonnet formula for such a manifold.
\begin{theorem}\label{thm:IthIsoperimetricRatio}
    Suppose that $(M^4,g)$ is a complete, simply connected, locally conformally flat manifold, the scalar curvature $R[g]$ is non-negative, and the $Q$-curvature $Q[g]$ is non-negative near infinity. Then, we have the Chern-Gauss-Bonnet formula
    $$
    \chi(M)-\frac{1}{32\pi^2}\int_{M}8Q[g]dV[g]=\sum_{i=1}^m\nu_i,
    $$
    where 
    \begin{align}
            \nu_i=\lim_{r\to\infty} \left(\frac{2}{\pi^2}\right)^{1/3}\frac{V_g(B^i_r)^{1/3}}{r^{4/3}} \label{IthIsoperimetricRatio}
    \end{align}
    is the asymptotic ratio of the i-th end on $M$. Here $V_g(\cdot)$ is the four-dimensional volume with respect to the metric $g$, and $B^i_{r}$ is the radial geodesic cylinder on the $i$-th end defined in (\ref{def:B_r}).
\end{theorem}
\begin{proof}
From the results of Schoen–Yau \cite{SchoenYau1988conformally} and Chodosh–Li \cite[Corollary 4]{LiChaoChodosh2024DevelopingMapInjective} discussed ahead of the statement of Theorem \ref{thm:IthIsoperimetricRatio}, we may regard $M$ conformally as a subdomain $\Omega \subseteq S^4$. First, for the complete manifold $(\Omega, g)$, obviously $R[g]\geq 0$. Since $Q[g]\geq0$ near infinity on $(M,g)$, there exists a compact set $N\subset M$ such that $Q[g]\geq 0$ on $M\backslash N$. Still, $\Phi(N)$ is a compact subset of $\Omega$. So, $Q[g]\geq 0$ on $\Omega \backslash \Phi(N)$, i.e. $Q[g]\geq 0$ near $S^4\backslash \Omega$. Then, by Lemma \ref{FiniteEnd11} of Ma–Qing, $M$ is conformally equivalent to $S^4\backslash\{p_i\}_{i=1}^k$.

Indeed, a small neighbourhood of each ${p_i}$ corresponds to an end of the manifold. We take a sufficiently small neighbourhood of ${p_i}$ where $Q[g] \ge 0$, and apply a stereographic projection of suitable scale. This neighbourhood becomes $(\mathbb{R}^4\backslash B^e_M(0), g = e^{2w_i}|dy|^2)$ with $Q[g] \ge 0$ and $R[g]\ge 0$. Rescaling the coordinates by $x=y/M$, we obtain precisely the conformally simple end in the form $(\mathbb{R}^4\backslash B, g = e^{2u_i}|dx|^2)$ introduced in the previous section with $Q[g] \ge 0$ and $R[g]\ge 0$, which is an end of the complete manifold $(M,g)$.  Thus, we can use all the results we mentioned previously in Section \ref{ChapterMaintheorem}. 

Consequently, we know that
    $$
        M=K\bigcup \left\{ \bigcup^m_{i=1}E_i\right\}
    $$
    where $(K,g)$ is a compact Riemannian manifold with boundary
    $$
        \partial K= \bigcup^m_{i=1}\partial E_i,
    $$
    And each $E_i$ satisfies
    $$
        (E_i,g)=(\mathbb{R}^4 \backslash B,~g=e^{2u_i}|dx|^2).
    $$
On each component $\partial E_i=\partial K\cap E_i$,
the outward unit normal of the compact core $K$ agrees
with $N_+[g]$. Hence the boundary curvature appearing
in the Chern--Gauss--Bonnet formula for $K$ is $T_+[g]$.
Applying the Chern–Gauss–Bonnet formula to the compact manifold $K$ with boundary, we obtain
\begin{align}
\chi(K)=\frac{1}{32\pi^2}\int_{K}\bigl(|W[g]|^2+8Q[g]\bigr)\,dV[g]
+\sum_{i=1}^{m}\frac{1}{4\pi^2}\int_{\partial E_i}(T_+[g]+L[g])\,d\sigma[g].\label{NewLiNe:608}
\end{align}
where
\begin{align*}
L[g]&=L_1+L_2\\
L_1 &= \frac{1}{3} R_{ijij} L_{aa} - R_{aNaN} L_{bb} + R_{aNbN} L_{ab} - R_{acbc} L_{ab}, \\
L_2 &= \frac{2}{9} L_{aa} L_{bb} L_{cc} - L_{aa} L_{bc} L_{bc} + L_{ab} L_{bc} L_{ca}.
\end{align*}
where $i,j=1,2,3,4$, $a,b,c=1,2,3$, and $N$ is the outward normal direction.
We now show that the boundary contribution of $L[g]$ vanishes.
On each end $E_i$, write
\[
g=e^{2u_i(x)}|dx|^2,\qquad x=r\theta,
\qquad r\ge1,\quad \theta\in S^3.
\]
Setting $t=\log r$, we obtain
\[
|dx|^2=e^{2t}(dt^2+g_{S^3}),
\qquad
g=e^{2U_i(t,\theta)}g_0,
\]
where
\[
g_0:=dt^2+g_{S^3},
\qquad
U_i(t,\theta):=u_i(e^t\theta)+t.
\]
The boundary $\partial E_i$ corresponds to the slice
$\Sigma_0:=\{0\}\times S^3$.

With respect to the product metric $g_0$, the slice $\Sigma_0$
is totally geodesic. Hence its second fundamental form satisfies
\[
L_{ab}[g_0]=0.
\]
Every term in the expressions for $L_1$ and $L_2$ contains
at least one factor of the second fundamental form.
It follows that
\[
L[g_0]\big|_{\Sigma_0}=0.
\]

By the conformal invariance of the boundary density
$L[g]\,d\sigma[g]$
\cite{ChangQing1997BoundaryOperator,ChangQingYangInvention2000},
we have
\[
L[g]\,d\sigma[g]
=
L[g_0]\,d\sigma[g_0]
=0
\qquad\text{on }\partial E_i.
\]
Here the normals for $g_0$ and $g$ are chosen with the same
orientation, namely, the outward orientation of the compact
core $K$.
Consequently,
\[
\int_{\partial E_i}L[g]\,d\sigma[g]=0
\qquad\text{for every }i=1,\ldots,m.
\]
For the Chern--Gauss--Bonnet formula used here, see
\cite[formula (0.7)]{ChangQingYangInvention2000}.
Next, we use Theorem \ref{MainTheorem},
    \begin{align}
    &~~~~\frac{1}{32\pi^2}\int_{E_i}8Q[g]dV[g]+\nu_i=
        \frac{1}{32\pi^2}\int_{\mathbb{R}^4\backslash B}8Q[g]e^{4u_i}dx+\nu_i\nonumber\\
        &=\frac{1}{4\pi^2}\int_{\partial B}T_+[g]e^{3u_i}d\sigma =\frac{1}{4\pi^2}\int_{\partial E_i}T_+[g]d\sigma[g], \label{NewLIne:614}
    \end{align}
    where 
    $$\nu_i=\lim_{r\to\infty}\left(\frac{2}{\pi^2}\frac{V_g(B^i_{r})}{r^4}\right)^{1/3},$$
where $B^i_{r}$ is the radial geodesic cylinder on the $i$-th end defined in (\ref{def:B_r}). Since the Weyl tensor vanishes on a locally conformally flat manifold, we have $|W[g]|^2 = 0$ on $(M,g)$. Finally, substituting the expression for the boundary integral from (\ref{NewLIne:614}) into (\ref{NewLiNe:608}), we obtain
    \begin{align}
        \chi(M)=\frac{1}{32\pi^2}\int_{M}8Q[g]dV[g]+\sum_{i=1}^m\nu_i.\nonumber
    \end{align}
Here $\chi(K)=\chi(M)$ because $K$ is homotopy equivalent to $M$.
\end{proof}\par

This theorem also admits an intrinsic formulation for the quantities $\nu_i$:
\begin{theorem}[Intrinsic characterization of the end invariants]\label{Maaainnini}
    Suppose that $(M^4,g)$ is a complete, simply connected, locally conformally flat manifold, the scalar curvature $R[g]$ is non-negative, and the $Q$-curvature $Q[g]$ is non-negative near infinity. Then, we have the Chern-Gauss-Bonnet formula
        \begin{align}
            \chi(M)-\frac{1}{32\pi^2}\int_{M}8Q[g]dV[g]=\sum_{i=1}^m\nu_i,\label{chernGaussBonnet}
        \end{align}
    where 
    \begin{align}
            \nu_i=\lim_{r\to\infty} \left(\frac{2}{\pi^2}\right)^{1/3}\frac{V_g(B_{r}(p)\cap M_i)^{1/3}}{r^{4/3}} \label{IthIsoperimetricRatioGeodsc11}
    \end{align}
    gives an intrinsic characterization of the invariant of the $i$-th end $M_i$. Here $p\in M$ is arbitrary, $B_r(p)$ is the geodesic ball of radius $r$ centred at $p$, and $V_g$ is the Riemannian volume. Thus $\nu_i$ is the cube root of the asymptotic volume ratio $V_g(B_r(p)\cap M_i)/((\pi^2/2)r^4)$. The limit is independent of $p$ and of changing $M_i$ within a compact set; in particular, it depends only on the intrinsic geometry of the end and does not require a choice of conformal coordinates.
\end{theorem}
\begin{proof}
Formula (\ref{chernGaussBonnet}) is a direct  result of Theorem \ref{thm:IthIsoperimetricRatio}. 

What remains is to prove (\ref{IthIsoperimetricRatioGeodsc11}).
Fix $p\in M$. On the $i$-th end $M_i = (\mathbb{R}^4\backslash B, g = e^{2u_i}|dx|^2)$, to avoid confusion we write $A^i_r := \{ z \in M_i: d_g(z, z/|z|) \leq r \}$ (this $ A^i_r$ is exactly the set $B_r$ defined in (\ref{def:B_r})). Let $m_0 = \max\{ d_g(x,p) : x \in \partial M_i \}$. Here $B_r(p)=\{z\in M:d_g(z,p)<r\}$ is the open geodesic ball. We assert that for $r>m_0+1$,
\begin{align}
A_{r-m_0-1}^i \subseteq B_r(p) \cap M_i \subseteq A^i_{r+m_0}. \label{haoguanxi}
\end{align}
    For the first inclusion “$\subseteq$” in (\ref{haoguanxi}), take $z \in A^i_{r-m_0-1}$. By the triangle inequality,
    \begin{align}
                d_g(z,p)\leq& d_g(z,z/|z|)+d_g(z/|z|,p) \leq r-m_0-1+d_g(z/|z|,p)\nonumber\\\leq& r-m_0-1+m_0=r-1<r.\nonumber
    \end{align}
    Naturally, there holds $A^i_{r-m_0-1}\subseteq M_i$. So $A_{r-m_0-1}^i\subseteq B_r(p)\cap M_i$.\par
    Regarding the second inclusion “$\subseteq$” in (\ref{haoguanxi}), we  take $y \in B_r(p)\cap M_i$ arbitrarily. Then, $d_g(y,y/|y|)\leq d_g(y,p)+d_g(p,y/|y|)\leq r+m_0$.
    Thus, the second inclusion in (\ref{haoguanxi}) is valid. \par
    Finally, combine (\ref{haoguanxi}) with the limit (\ref{IthIsoperimetricRatio}). Since $(r-m_0-1)/r\to1$ and $(r+m_0)/r\to1$, these fixed radius shifts do not change the normalized volume limits. The squeeze theorem therefore gives (\ref{IthIsoperimetricRatioGeodsc11}).
    To make the intrinsic nature explicit, let $q\in M$ and set $c=d_g(p,q)$. For $r>c$, the triangle inequality gives
    \[
        B_{r-c}(p)\cap M_i\subseteq B_r(q)\cap M_i
        \subseteq B_{r+c}(p)\cap M_i.
    \]
    Dividing volumes by $r^4$ and taking limits proves independence of the base point. Replacing $M_i$ by another representative of the same end changes its intersection with a geodesic ball by a subset of a fixed compact set, whose volume is finite. This contributes $O(r^{-4})$ to the normalized volume and hence does not change the limit. All quantities in the resulting formula are defined by $g$, without reference to conformal coordinates. 
\end{proof}
\begin{corollary}\label{COroll}
            Suppose that $(M^4,g)$ is a complete, connected,   locally conformally flat manifold, the scalar curvature $R[g]$ is non-negative, and the $Q$-curvature $Q[g]$ is non-negative near infinity. If the fundamental group $\pi_1(M)$ is finite, then the result of Theorem  \ref{Maaainnini} is still valid for $\widetilde{M}$, where $p:\widetilde{M}\to M$ is the universal covering map of $M$.
\end{corollary}
\begin{proof}
     It remains to verify that $\widetilde{M}$ satisfies all the assumptions of Theorem \ref{Maaainnini}. A Riemannian manifold is a topological manifold and is therefore locally path-connected and semilocally simply connected. From \cite[Corollary 82.2]{TopologyMunkrus}, $M$ has a universal cover $p \colon \widetilde{M} \to M$. By \cite[Theorem 54.4]{TopologyMunkrus}, we know that
    \begin{align}
        \phi:\pi_1(M,x_0)\to p^{-1}(x_0)
    \end{align}
    is bijective, since the universal cover is simply connected. Hence, $\widetilde{M}$ is a universal covering space with finite fibers $p^{-1}(x_0)$. Besides, we know that $(\widetilde{M},p^*g)$ is complete since $(M,g)$ is complete.  \par
    Choose a compact set $K\subseteq M$ such that $Q[g]\geq0$ on $M\backslash K$. Since $p$ is a finite covering, $p^{-1}(K)$ is compact. Equip $\widetilde M$ with $\widetilde g:=p^*g$. Curvature is preserved under the local isometry $p$, so $R[\widetilde g]=R[g]\circ p\geq0$ everywhere and $Q[\widetilde g]=Q[g]\circ p\geq0$ on $\widetilde M\backslash p^{-1}(K)$. The pullback metric is also locally conformally flat and complete. Thus $(\widetilde M,\widetilde g)$ satisfies all the assumptions of Theorem \ref{Maaainnini}, which establishes the corollary. 
\end{proof}
\begin{corollary}
Under the assumptions of Corollary \ref{COroll}, let
\[
p:\widetilde M\longrightarrow M
\]
be the universal covering map, and equip $\widetilde M$
with the pullback metric $\widetilde g:=p^*g$.
Let $\mathcal C$ be the group of deck transformations of $p$.
Then $\mathcal C\cong\pi_1(M)$ is finite and acts on
$(\widetilde M,\widetilde g)$ by isometries.
Moreover, with its quotient metric $g_{\mathrm{quot}}$,
the quotient $\widetilde M/\mathcal C$ is globally isometric
to $(M,g)$ via the natural map
\[
k:\widetilde M/\mathcal C\longrightarrow M,
\qquad k([x])=p(x).
\]
\end{corollary}

\begin{proof}
Since $\widetilde M$ is simply connected, the universal
covering $p$ is regular, and its deck transformation group
satisfies
\[
\mathcal C\cong\pi_1(M).
\]
In particular, $\mathcal C$ is finite by assumption.

For every $\gamma\in\mathcal C$, the identity
$p\circ\gamma=p$ implies
\[
\gamma^*\widetilde g
=\gamma^*(p^*g)
=(p\circ\gamma)^*g
=p^*g
=\widetilde g.
\]
Thus every deck transformation is an isometry of
$(\widetilde M,\widetilde g)$.

The action of $\mathcal C$ is free and properly discontinuous.
Consequently, $\widetilde M/\mathcal C$ is a smooth manifold
and carries a unique quotient metric $g_{\mathrm{quot}}$
such that the quotient projection
\[
q:\widetilde M\longrightarrow\widetilde M/\mathcal C
\]
is a local isometry; equivalently,
\[
q^*g_{\mathrm{quot}}=\widetilde g.
\]

The map $k([x]):=p(x)$ is well defined because $p$ is
constant on every $\mathcal C$-orbit.
Since $p$ is a regular covering, its fibers are precisely
the $\mathcal C$-orbits. Hence $k$ is bijective.
Furthermore, the identity
\[
p=k\circ q
\]
and local covering charts show that $k$ is a local
diffeomorphism. Therefore, $k$ is a diffeomorphism.

Finally,
\[
q^*(k^*g)
=(k\circ q)^*g
=p^*g
=\widetilde g
=q^*g_{\mathrm{quot}}.
\]
Since $q$ is a surjective local diffeomorphism, this implies
\[
k^*g=g_{\mathrm{quot}}.
\]
Thus $k$ is a global Riemannian isometry.
\end{proof}

\section{Example}\label{ChapterExample}
From \cite[Theorem 1.2]{liWeiXu2025geometry}, it is known that on a complete, conformally flat manifold $(\mathbb{R}^4, g = e^{2u}|dx|^2)$, if the scalar curvature is non-negative and the $Q$-curvature is non-negative, then the Ricci curvature is also non-negative (semi-positive definite). The condition on the $Q$-curvature in that result is a global one, and the background space is the full Euclidean space $\mathbb{R}^4$. In contrast, our assumption in this paper is that the $Q$-curvature is non-negative only near infinity, which is a non-global condition, and the background space of the objects under study is not $\mathbb{R}^4$. To highlight the difference between our work and that of \cite{liWeiXu2025geometry}, we present the following example, which shows that on a locally conformally flat manifold, if the scalar curvature and the $Q$-curvature are non-negative locally, the Ricci curvature may not be locally non-negative (semi-positive definite).

We consider the conformal metric  $g=e^{2u}|dx|^2$ on a punctured ball $B^e_1(0)\backslash\{0\}$, where
$$u(x)=\frac{x_1}{|x|^2}.$$
This function is biharmonic for $x=(x_1,x_2,x_3,x_4)\in B^e_1(0)\backslash\{0\}$. Hence $Q[g]=0$. Moreover, both $R[g]>0$ and $Ric_g(\partial_1,\partial_1)<0$ hold near the point $(3/8,1/100,0,0)$. \par
Let us make a brief computation. Denote
\begin{align}
    u_{i,j}:=&\partial^2_{ij}u-\varGamma_{ij}^k\partial_ku-\partial_iu\partial_ju+\frac{1}{2}g_{ij}g^{kl}\partial_ku\partial_lu\nonumber\\
    =&\partial^2_{ij}u-\partial_iu\partial_ju+\frac{1}{2}\delta_{ij}\sum_{k=1}^{4}(\partial_ku)^2,\nonumber
\end{align}
where $\Gamma_{i,j}^k$ is the Christoffel symbol of $|dx|^2$. 
Using the transformation formula for the components of the Ricci tensor, $R_{ij}=R^e_{ij}-2u_{i,j}-g_{ij}g^{kl}u_{k,l}$, we obtain, after a direct calculation, that 
\begin{align}
    Ric_g(\partial_1,\partial_1)=8\frac{x_1}{|x|^4}\left(2-\frac{x_1}{|x|^2}\right)\left(1-\frac{x_1^2}{|x|^2} \right)\nonumber
\end{align}
is negative near $x=(3/8,1/100,0,0)$. Thus, the Ricci tensor is not locally semi-positive definite, while the scalar curvature 
\begin{align}
    R_g:=&e^{-2u}(R_{g_e}-6\delta^{ij}u_{i,j})\nonumber\\
    =&-6e^{-2u}\sum_{k=1}^{4}u_{k,k}\nonumber\\
     =&6e^{-2u}\left(4\frac{x_1}{|x|^4}-\frac{1}{|x|^4}\right)\nonumber
\end{align}
is positive near $x=(3/8,1/100,0,0)$. To verify it, note that $4x_1-1>0$. As for the $Q$-curvature, from (\ref{Def:QCurvature}), we have 
$$
Q[g]=\frac{1}{2}e^{-4u}\Delta^2u=0,
$$
since this $u$ is biharmonic.
\section{Appendix}\label{ChapterAppendix}
For the reader’s convenience and to keep the paper self-contained, we provide proofs of the lemmas, which are heavily inspired by the works \cite{ChangQingYangDuke2000,ChangQingYangInvention2000,grafakos2008classical,lee1977generalizationsLHôpital}.

\begin{lemma}\label{Lemma:intfgeqf}
        Let $f:[a,\infty)\to(0,\infty)$ be continuous, where $a\geq1$. Suppose that $f(s)=s^{p+o(1)}$ as $s\to\infty$, with $p\in\mathbb{R}$. Then, for all sufficiently large $s$,
        \[
            \int_a^s f(r)\,dr>f(s).
        \]
\end{lemma}
\begin{proof}
To arrive at this conclusion, consider the asymptotic behavior of $  f(s) = s^{p + \varepsilon(s)}  $ where $  \varepsilon(s) \to 0  $ as $  s \to \infty  $, and analyze the integral $  I(s) = \int_a^s f(r) \, dr  $ in different regimes based on $  p  $.
For any $  \delta > 0  $, there exists $  S>a  $ such that for $  r > S  $, $  |\varepsilon(r)| < \delta  $, so $  r^{p - \delta} < f(r) < r^{p + \delta}  $. Thus, for large $  s  $,
\begin{align}
    \int_S^s r^{p - \delta} \, dr < \int_S^s f(r) \, dr < \int_S^s r^{p + \delta} \, dr,\nonumber
\end{align}
and $  I(s)  $ equals a constant plus the middle term.

Case $  p > -1  $: The integrals are asymptotically $  s^{p - \delta + 1}/(p - \delta + 1)  $ (lower) and $  s^{p + \delta + 1}/(p + \delta + 1)  $ (upper). Since $  \delta  $ is arbitrary, $  \log I(s) / \ln s \to p + 1  $, so $  I(s) = s^{p + 1 + o(1)}  $. Then $  I(s)/f(s) = s^{1 + o(1)} \to \infty  $; hence, $  I(s) > f(s)  $ for large $  s  $.\par
Case $  p < -1  $: Both bounding integrals converge as $  s \to \infty  $ (since $  p \pm \delta + 1 < 0  $ for small $  \delta  $), so $  I(s) \to L  $ for some finite $  L > 0  $. Meanwhile, $  f(s) = s^{p + o(1)} \to 0  $ (exponent $  \to p < 0  $), so $  I(s) > f(s)  $ for large $  s  $.\par
Case $  p = -1  $: The lower bound converges, but the upper diverges like $  s^\delta / \delta \to \infty  $. Thus, $  I(s)  $ may converge to a finite limit $  L > 0  $ or diverge, depending on $  \varepsilon(s)  $. If it converges, $  f(s) = s^{-1 + o(1)} \to 0  $, so $  I(s) > f(s)  $ eventually. If it diverges (e.g., like $  \ln s  $, $  (\ln s)^2/2  $, etc., in examples where $  \varepsilon(s)  $ is small or positive), $  I(s) \to \infty  $ while $  f(s) \to 0  $, so again $  I(s) > f(s)  $ for large $  s  $.\par

In all cases, the positivity of $  f  $ ensures $  I(s)  $ is at least a positive constant for large $  s  $, and the asymptotics confirm the strict inequality holds eventually.
\end{proof}
\begin{lemma}[\cite{limingxiang2024VolumeEntropy}.Lemma 2.3]
    Let \( f \) satisfy \( f \in L^\infty_{\text{loc}}(\mathbb{R}^4\backslash B) \) and \( f \in L^1(\mathbb{R}^4\backslash B) \). Then for \( |x| \gg 1 \), \(0<\delta<1\), there holds
\begin{align}
    \mathcal{L}(f)(x) = (-\alpha_0 + o(1)) \log |x|+ o(1)\log\frac{1}{\delta} + \frac{1}{8\pi^2} \int_{B_\delta^e(x)} \log \frac{1}{|x-y|} f(y) \, dy\label{ALogEquRepresentation}
\end{align}
    where \( o(1) \to 0 \) as \( |x| \to \infty \), $\alpha_0=\frac{1}{8\pi^2}\int_{\mathbb{R}^4\backslash B}f(y)dy$ and $\mathcal{L}(f)$ is defined in (\ref{defofLF}). 
\end{lemma}
\begin{proof}
   This is a minor adaptation of \cite[Lemma 2.3]{limingxiang2024VolumeEntropy}. Although the domain of integration is $\mathbb{R}^4\backslash B$ instead of $\mathbb{R}^4$, the same argument used in \cite[Lemma 2.3]{limingxiang2024VolumeEntropy} for $n=4$ can be applied to show that 
 \begin{align}
         \mathcal{L}(f)(x) = (-\alpha_0 + o(1)) \log |x| + \frac{1}{8\pi^2} \int_{B_{1}^e(x)} \log \frac{1}{|x-y|} f(y)dy. \label{lineq:835}
 \end{align}
    Then for $|x|>>1$, split  $B_{1}^e(x)$ into $\left(B_{1}^e(x)\backslash B_\delta^e(x)\right)\bigcup B_\delta^e(x)$.
    For every $y\in B_{1}^e(x)\backslash B_\delta^e(x)$, we have $\delta \leq |x-y|<1$, and hence $0<\log\frac{1}{|x-y|}\leq \log\frac{1}{\delta}$. Therefore, the integral over $B_{1}^e(x)\backslash B_\delta^e(x)$ satisfies 
    \begin{align}
        \left|\frac{1}{8\pi^2}\int_{B_{1}^e(x)\backslash B_\delta^e(x)} \log \frac{1}{|x-y|} f(y)dy\right|\nonumber
        \leq&\frac{1}{8\pi^2}\int_{B_{1}^e(x)\backslash B_\delta^e(x)} \log \frac{1}{|x-y|}|f(y)|dy\nonumber\\
        \leq&\frac{1}{8\pi^2}\log\frac{1}{\delta}\int_{B_{1}^e(x)\backslash B_\delta^e(x)}|f(y)|dy\nonumber\\
         \leq&\frac{1}{8\pi^2}\log\frac{1}{\delta}\int_{B_{1}^e(x)}|f(y)|dy\nonumber\\=&\frac{1}{8\pi^2}\log\frac{1}{\delta}\cdot o(1), ~as~ |x|\to\infty,\label{Line11:843}
    \end{align}
since \( f \in L^1(\mathbb{R}^4\backslash B) \) so that \begin{align}
    \int_{B_{1}^e(x)}|f(y)|dy\to 0,
\end{align}
as $|x|\to \infty$.

Thus, for $\delta\in(0,1)$, combining (\ref{Line11:843}) with (\ref{lineq:835}) gives
 \begin{align}
         \mathcal{L}(f)(x) =& (-\alpha_0 + o(1)) \log |x| + \frac{1}{8\pi^2}\left( \int_{B_{1}^e(x)\backslash B_{\delta}^e(x)}+\int_{B_{\delta}^e(x)}\right) \log \frac{1}{|x-y|} f(y)dy\nonumber\\
         =&(-\alpha_0 + o(1)) \log |x| + o(1)\log\frac{1}{\delta} + \frac{1}{8\pi^2} \int_{B_{\delta}^e(x)}\log \frac{1}{|x-y|} f(y)dy,\label{lineq:847}
 \end{align}
where $o(1)\to 0$, as $|x|\to+\infty$. 
\end{proof}
The following lemma provides a bound for the spherical average of the logarithmic potential. Define $$\mathcal {A}(r):=\fint_{\partial B^e_r}\mathcal{L}(f)(y)d\sigma_y,~r>0.$$ We refer to $\mathcal A$ as the spherical average in the following lemma.
\begin{lemma}\label{Rwidetildegeq0}
    Suppose that $f\in L^1 (\mathbb{R}^4\backslash B)$ with $f$ non-negative. Then $$-C_1\log R-C_2\leq \mathcal{A}(R)\leq C_3,$$
    for $R$ large enough, where $\mathcal{A}(R)$ is defined in (\ref{defofA}) and $C_1$, $C_2$ and $C_3$ are constants.
\end{lemma}
\begin{proof}
 Decompose $\mathbb{R}^4\backslash B$ as $B_2^e\cup(B_2^e)^c$, using the end-truncated open-ball notation fixed above. Here $(B_2^e)^c$ denotes the complement relative to $\mathbb{R}^4\backslash B$. Then
    \begin{align}
        \mathcal{A}(R)&=\frac{1}{8\pi^2}\fint_{\partial B_{R}^{e}}{\int_{B_{2}^{e} \bigcup{(B_{2}^{e} )^c}}{\log \frac{|y|}{|x-y|}f\left( y \right) dyd\sigma_x }}\nonumber
        \\&:=\mathcal{A}_1(R)+\mathcal{A}_2(R).\label{Newline:579}
    \end{align}
    We now derive lower bounds for $\mathcal{A}_1(R)$ and $\mathcal{A}_2(R)$.\par
    The first term is
    \begin{align}
        \mathcal{A}_1(R)&=\frac{1}{8\pi^2}\fint_{\partial B_{R}^{e}}{\int_{{B_{2}^{e}}}{\log \frac{|y|}{|x-y|}f\left( y \right) dyd\sigma_x }}.
    \end{align}
    For any $y \in B_2^e$, we choose $R$ large enough so that $|x|/|x-y|$ is close to $1$ and consequently $\bigl|\log\frac{|x|}{|x-y|}\bigr| \le 2$ for $x \in \partial B_{R}^{e}$. Then, 
    \begin{align}
        \left|\log |x|+\log\frac{|y|}{|x-y|}\right|=\left|\log\frac{|x|}{|x-y|}+\log|y| \right|\leq |\log|y||+2,\nonumber
    \end{align}
    So
    \begin{align}
        \log \frac{|y|}{|x-y|}\geq-\log |x|-|\log|y||-2.\label{Newline:591}
    \end{align}
Set $\alpha:=\frac{1}{8\pi^2}\int_{B_2^e}f(y)\,dy$. Since $|\log|y||$ is bounded on $B_2^e$ and $f$ is integrable on the end, we obtain, for a constant $C$ independent of $R$,
    \begin{align}
        \mathcal{A}_1(R)&=\frac{1}{8\pi^2}\fint_{\partial B_{R}^{e}}{\int_{{B_2^e\backslash B}}{\log \frac{|y|}{|x-y|}f\left( y \right) dyd\sigma_x }}\nonumber\\
        &\geq\frac{1}{8\pi^2}\fint_{\partial B_{R}^{e}}{\int_{B_2^e\backslash B}{-(\log |x|+|\log|y||+2)f\left( y \right) dyd\sigma_x }}\nonumber\\
        &\geq\fint_{\partial B_{R}^{e}}{{(-\alpha\log |x|-C)d\sigma_x }}\nonumber\\
        &\geq-\alpha\log|x|-C=-\alpha \log R-C, \label{A_1result}
    \end{align}
    for sufficiently large $R$. \par
    Then we consider $\mathcal{A}_2(R)$. Using Fubini's theorem and the computational techniques in \cite[P593,D.3]{grafakos2008classical}, we obtain
    \begin{align}
        \mathcal{A}_2(R)&=\frac{1}{8\pi^2}\fint_{\partial B_{R}^{e}}{\int_{{(B_{2}^{e})^c}}{\log \frac{|y|}{|x-y|}f\left( y \right) dyd\sigma_x }}\nonumber\\
        &=\frac{1}{8\pi^2}\int_{{(B_{2}^{e})^c}}\fint_{\partial B_{R}^{e}}{\log \frac{|y|}{|x-y|}d\sigma_x f\left( y \right) dy}\nonumber\\
        &=\frac{1}{16\pi^2}\int_{{(B_{2}^{e})^c}}\fint_{\partial B_{R}^{e}}{\log \frac{|y|^2}{|x-y|^2}d\sigma_x f\left( y \right) dy}\nonumber\\
        &=\frac{1}{16\pi^2}\int_{{(B_{2}^{e})^c}}\fint_{\partial B_{1}^{e}\left( 0 \right)}{\log \frac{|y|^2}{|R\theta-y|^2}d\sigma_\theta f\left( y \right) dy}\nonumber\\
        &=\frac{1}{16\pi^2}\int_{{(B_{2}^{e})^c}}\fint_{\partial B_{1}^{e}\left( 0 \right)}{\log \frac{|y|^2}{R^2+|y|^2-2R\theta\cdot y}d\sigma_\theta f\left( y \right) dy}.\label{Newline:606}
    \end{align}
Applying the result of \cite[P593, D.3]{grafakos2008classical} to (\ref{Newline:606}), we obtain
    \begin{align}
        \mathcal{A}_2(R)&=\frac{1}{16\pi^2}\int_{{(B_{2}^{e})^c}}\fint_{\partial B_{1}^{e}\left( 0 \right)}{\log  \frac{|y|^2}{R^2+|y|^2-2R\theta\cdot y}d\sigma_\theta f\left( y \right) dy}\nonumber\\
        &=\frac{1}{16\pi^2}\frac{\omega_{2}}{\omega_{3}}\int_{{(B_{2}^{e})^c}}\int_{-1}^{1}{\log \frac{|y|^2}{R^2+|y|^2-2Rs|y|}\sqrt{1-s^2}ds f\left( y \right) dy}.\label{Newline:611}
    \end{align}
Since $R^2 + |y|^2 - 2Rs|y| \le 2R^2 + 2|y|^2$, it follows from (\ref{Newline:611}) that
        \begin{align}
        \mathcal{A}_2(R)&\geq-\frac{1}{16\pi^2}\frac{\omega_{2}}{\omega_{3}}\int_{({B_{2}^{e})^ c}}\int_{-1}^{1}{\log \frac{2R^2+2|y|^2}{|y|^2}\cdot 1~ds f\left( y \right) dy}\nonumber\\
        &=-\frac{1}{16\pi^2}\frac{\omega_{2}}{\omega_{3}}\int_{({B_{2}^{e})^c}}\int_{-1}^{1}{\log (\frac{2R^2}{|y|^2}+2)ds f\left( y \right) dy}\nonumber\\
        &\geq-C\log R,\label{Result_A_2}
    \end{align}
    where the last inequality is obtained from $|y|\geq2$. Consequently, combining (\ref{Newline:579}) with the lower bounds (\ref{A_1result}) and (\ref{Result_A_2}) for $\mathcal{A}_1(R)$ and $\mathcal{A}_2(R)$, we obtain
    \begin{align}
    \mathcal{A}(R) \geq -C_1 \log R - C_2. \label{LowerBoundOfA}
    \end{align}
    Now we establish a constant upper bound for $\mathcal{A}$. Using the same technique as in the result of \cite[P593, D.3]{grafakos2008classical}, we obtain
    \begin{align}
         \mathcal{A}(R)&=\frac{1}{16\pi^2}\int_{\mathbb{R}^4\backslash B}\fint_{\partial B_{1}^{e}\left( 0 \right)}{\log \frac{|y|^2}{R^2+|y|^2-2R\theta\cdot y}d\sigma_\theta f\left( y \right) dy}\nonumber\\
         &=-\frac{1}{16\pi^2}\frac{\omega_{2}}{\omega_{3}}\int_{\mathbb{R}^4\backslash B}\int_{-1}^{1}{\log \frac{R^2+|y|^2-2Rs|y|}{|y|^2}(1-s^2)^{1/2}ds f\left( y \right) dy}.\nonumber
    \end{align}
    Since $R^2 - 2Rs|y| \ge -s^2|y|^2$, it follows that
    \begin{align}
        \mathcal{A}(R)&=-\frac{1}{16\pi^2}\frac{\omega_{2}}{\omega_{3}}\int_{\mathbb{R}^4\backslash B}\int_{-1}^{1}{\log \frac{R^2+|y|^2-2Rs|y|}{|y|^2}(1-s^2)^{1/2}ds f\left( y \right) dy}\nonumber\\
        &\leq -\frac{1}{16\pi^2}\frac{\omega_{2}}{\omega_{3}}\int_{\mathbb{R}^4\backslash B}\int_{-1}^{1}{\log \frac{(1-s^2)|y|^2}{|y|^2}(1-s^2)^{1/2}ds f\left( y \right) dy}\nonumber\\
        &\leq\tilde C\int_{\mathbb{R}^4\backslash B}f(y)dy\leq C_3.\label{UpperBoundOfA}
    \end{align}
    The last two inequalities follow from the finiteness of 
    $$-\int_{-1}^{1}(1-s^2)^{1/2}\log(1-s^2)ds=\frac{\pi}{2}(\log 4-1)$$
    and the assumption $f\in L^1 (\mathbb{R}^4\backslash B)$.\par
    In summary, from the bounds (\ref{LowerBoundOfA}) and (\ref{UpperBoundOfA}), we obtain  
    $$-C_1\log R-C_2\leq \mathcal{A}(R)\leq C_3,$$
    for all sufficiently large $R$.
\end{proof}
\begin{lemma}\label{MainLemma}
    Suppose $(\mathbb{R}^4\backslash B,g=e^{2u}|dx|^2)$ is a conformally flat simple end of a complete manifold $(M,g)$, and $R[g]\geq0$ and $Q[g]\geq0$ on the end. By Lemma \ref{NormalMetric}, $Q[g]\in L^1(\mathbb{R}^4\backslash B,g)$ and the normal representation defining $D$ is available. Then 
    $$\lim_{r\to \infty}rD'(r)=-1+\frac{1}{4\pi^2}\int_{\partial B}T_+[g]e^{3u(y)}d\sigma_y-\frac{1}{4\pi^2}\int_{\mathbb{R}^4\backslash B}Q_{g}(y)e^{4u(y)}dy,$$
    where $D(r)$ is defined in (\ref{line:272DefOfD}) and $T_+$ is the boundary $T$-curvature associated with $P_3$, computed with respect to $N_+[g]=e^{-u}\partial_r$.
\end{lemma}
\begin{proof}
    Let $\bar g:=e^{2D(|x|)}|dx|^2$ in the following proof. It is a direct result of the calculation that 
    \begin{equation}
        \Delta^2\mathcal{L}(f)(z)=2Q_g(z)e^{4u(z)}.\label{line:552}
    \end{equation}
    Using the Paneitz operator equation (\ref{Def:PaneitzOperator}), we have 
    \begin{equation}
        \Delta^2D(|x|)=2Q_{\bar g}(x)e^{4D(|x|)}.\label{line:555}
    \end{equation}
    Set $d\sigma_r:=r^3d\sigma$, where $d\sigma$ is the standard volume form on the unit sphere. Thus, we have
    \begin{align}
        &2\int_{\mathbb{R}^4\backslash B}Q_{\bar g}(y)e^{4D(|y|)}dy=\int_{\mathbb{R}^4\backslash B}\Delta^2D(|y|)dy\nonumber\\
        &=\int_{1}^{\infty}\int_{\partial B_r^e}\Delta^2D(r)d\sigma_r dr=\int_{1}^{\infty}\Delta^2D(r)\cdot\left|\partial B_r^e\right| dr\nonumber\\&=\int_{1}^{\infty}\int _{\partial B_r^e}\Delta ^2\mathcal{L}(f)(z)d\sigma_z dr=\int_{\mathbb{R}^4\backslash B}\Delta^2\mathcal{L}(f)(y)dy\nonumber\\
        &=2\int_{\mathbb{R}^4\backslash B}Q_{ g}(y)e^{4u(y)}dy,\label{line:562}
    \end{align}
    since the spherical average  commutes with $\Delta$, $\gamma\log r$ and $b(x/|x|^2)$ are bi-harmonic, and $D(r)$ is a radially symmetric function.\par
    Next, we can apply the method used by Chang-Qing-Yang in \cite[2. Symmetric case]{ChangQingYangDuke2000} to handle the radially symmetric term $D(r)$ and show that
    $\lim_{t\to\infty}tD'(t)$ exists. For the reader’s convenience, we present the proof here. Introducing the cylindrical coordinate $t=\log r$ with $r=|x|$ and setting $v(t):=D(e^t)+t$, we rewrite the Paneitz operator equation (\ref{line:555}). Since
    \[
     \frac{\partial}{\partial r} = e^{-t} \frac{\partial}{\partial t}, \quad \frac{\partial^2}{\partial r^2} = e^{-2t} \left( \frac{\partial^2}{\partial t^2} - \frac{\partial}{\partial t} \right),
    \]
    the Laplacian in standard  cylindrical coordinates is
    \[
        \Delta = \frac{\partial^2}{\partial r^2} + \frac{3}{r} \frac{\partial}{\partial r} = e^{-2t} \left( \frac{\partial^2}{\partial t^2} + 2 \frac{\partial}{\partial t} \right).
    \]
    Equation (\ref{line:555}) now becomes
    \begin{align}
        \left(\frac{\partial^2}{\partial t^2}-2\frac{\partial}{\partial t}\right)\left(\frac{\partial^2}{\partial t^2}+2\frac{\partial}{\partial t}\right)(D+t)=2Q[\bar g]e^{4D+4t}.\label{zanshi}
    \end{align}
    Here we have replaced $D$ with $D+t$ on the left side of (\ref{zanshi}). Consequently, 
    \begin{equation}
        v^{(4)}(t)-4v^{(2)}(t)=F(t),~0\leq t <\infty, \label{line:567}
    \end{equation}
    where $F(t):=2Q[\bar g]e^{4v}$. By (\ref{line:562}), $F\in L^1([0,\infty))$. Equation (\ref{line:567}) is a non-homogeneous ODE. First we find a particular solution to (\ref{line:567}); then the general solution can be written as
    \begin{align}
        v(t)=c_0+c_1t+c_2e^{-2t}+c_3e^{2t}+f(t).\label{line:571}
    \end{align}
    Using the method of variation of constants, we analyze the properties of the particular solution $f(t)$. Setting $f''(t) := C(t) e^{-2t}$ and substituting into (\ref{line:567}), we find that $C(t)$ satisfies 
    \begin{align}
         C''(t)-4C'(t)=F(t)e^{2t}.\label{zanshide}
    \end{align}
Applying the integrating-factor method and multiplying by $e^{-4t}$ in (\ref{zanshide}), we obtain, after a short computation,
    \begin{align*}
        C(t) &= -\int_{0}^{t} e^{4x} \int_{x}^{\infty} F(y)e^{-2y}dydx\nonumber\\
        &= -\frac{1}{4} e^{4t} \int_{t}^{\infty} F(x)e^{-2x} \, dx - \frac{1}{4} \int_{0}^{t} F(x)e^{2x} \, dx+\frac{1}{4}K_1,
    \end{align*}
    where $K_1$ is defined in (\ref{line:593}).
    Thus, $f''(t):=C(t)e^{-2t}$ can be represented as
    \begin{align}
        f''(t) = -\frac{1}{4} e^{2t} \int_{t}^{\infty} F(x)e^{-2x} \, dx - \frac{1}{4} e^{-2t} \int_{0}^{t} F(x)e^{2x} \, dx+\frac{1}{4}K_1e^{-2t}.\label{line:584}
    \end{align}
    Applying the Newton–Leibniz formula to each term that has a finite limit at the relevant endpoint, either zero or infinity, we can express $f'(t)$ as
    \begin{align}
    f'(t) =& \frac{1}{8}  e^{-2t} \int_{0}^{t} F(x)e^{2x} \, dx - \frac{1}{8}e^{2t} \int_{t}^{\infty} F(x)e^{-2x} \, dx\nonumber\\
    &+\frac{1}{8}(K_1 - K_2) + \frac{1}{8}\left\{\int_{t}^{\infty} F(x) \, dx - \int_{0}^{t} F(x) \, dx \right\}-\frac{1}{8}K_1 e^{-2t},
    \label{line:589}
     \end{align}
    where 
    \begin{align}
        K_1 &= \lim_{t \to 0} e^{2t} \int_{t}^{\infty} F(x)e^{-2x} \, dx= \int_{0}^{\infty}F(x)e^{-2x}dx,\label{line:593}\\
        K_2 &= \lim_{t \to \infty} e^{-2t} \int_{0}^{t} F(x) e^{2x} \, dx.\label{line:594}
    \end{align}
    It is easy to see that 
    \begin{align}
    |K_1|\leq \int_{0}^{\infty} |F(x)|dx< \infty.    \end{align}
    Next, we show that 
    \begin{align}
        K_2=0.\label{line:599}
    \end{align}
    In detail, we prove that $K_2=0$. Since $F$ is absolutely integrable, we have 
    \begin{align}
        e^{-2t} \int_{0}^{t} F(x) e^{2x} \, dx &= e^{-2t} \left\{ \int_{0}^{\log t} F(x) e^{2x} \, dx + \int_{\log t}^{t} F(x) e^{2x} \, dx \right\}\nonumber\\
        &\leq e^{-2t+2\log t} \int_{0}^{\log t} |F(x)| \, dx + e^{-2t}\int_{\log t}^{t} |F(x)|e^{2t} \, dx\nonumber\\
        &\leq e^{-2t+2\log t} \int_{0}^{\infty} |F(x)| \, dx + \int_{\log t}^{\infty} |F(x)| \, dx.\label{line:599111}
    \end{align}
Thus, letting $t \to \infty$, we obtain $K_2 = 0$. \par
    To clarify the coefficients in (\ref{line:571}), we need the asymptotic behavior of $f'$ and $f''$ near infinity. Following the argument that gives $K_2=0$ and taking the limits $t\to0$ and $t\to+\infty$ in (\ref{line:589}), respectively, we obtain
    \begin{align}
        \lim_{t\to0}f'(t)&=\frac{1}{8}\left(-K_1+\int_{0}^{+\infty}F(s)ds\right),\label{line:606}\\
        \lim_{t\to+\infty}f'(t)&=\frac{1}{8}\left(K_1-\int_{0}^{+\infty}F(s)ds\right).\label{line:607}
    \end{align}
From the expression (\ref{line:584}) for $f''(t)$, letting $t \to 0$ and comparing the definition (\ref{line:593}) of $K_1$ with (\ref{line:599111}), we obtain 
    \begin{align}
        \lim_{t\to0}f''(t)=0.\label{line:613}
    \end{align}
    Differentiating (\ref{line:584}) and using an argument similar to the one that established $K_2 = 0$, we obtain 
    \begin{align}
        \lim_{t\to\infty} f'''(t) =&\lim_{t\to\infty} -\frac{1}{2} e^{2t} \int_{t}^{\infty} F(x) e^{-2x} \, dx 
        + \frac{1}{2} e^{-2t} \int_{0}^{t} F(x) e^{2x} \, dx \nonumber\\
        =& \lim_{t\to\infty} -\frac{1}{2} e^{2t} \int_{t}^{\infty} F(x) e^{-2x} \, dx+\frac{1}{2}K_2   
        =0, \label{DerivativeThreeTimesIsZero}
    \end{align}
    where the limit vanishes since
    \begin{align}
        \left|e^{2t} \int_{t}^{\infty} F(x) e^{-2x} \, dx\right|\leq  \int_{t}^{\infty} |F(x)| \, dx \to 0, 
    \end{align}
    as $t\to\infty$, since $F(t)\in L^1[0,\infty)$.
\par
    Here, we must determine $c_3$. Applying Lemma \ref{Rwidetildegeq0}, we know that
    $$-(C_1-\gamma)\log s-C_2-C_4\leq D(s)=\mathcal{A}(s)+\gamma \log s+\bar h(s)\leq C_3+\gamma \log s+C_4,$$
    for $s$ large. Recalling the definition of $D(s)$ in (\ref{line:272DefOfD}), we see that
    \begin{align}
        -(C_1-\gamma-1)t-C_2-C_4\leq v(t)=D(e^t)+t\leq C_3+C_4+(\gamma+1)t,\label{GoodBound}
    \end{align}
     for sufficiently large $t$. The bound (\ref{GoodBound}) for $v(t)$ implies that $c_3=0$. In summary, 
    \begin{equation}
        v(t)=c_0+c_1t+c_2e^{-2t}+f(t).\label{line:639}
    \end{equation}
    \par
     From (\ref{DerivativeThreeTimesIsZero}) and (\ref{line:639}), we have 
    \begin{align}
        \lim_{t\to+\infty}v'''(t)=\lim_{t\to+\infty}f'''(t)=0.\label{ANewLine:773}
    \end{align}    
As for the special case of the standard cylinder $S^3 \times [0,t]$, since $v(t)$ is radially symmetric (constant on each slice $S^3$) and the background metric of $\bar g$  is the cylinder metric  $g_{0}:=dt^2 + d\theta^2$ (with $R=6$, $H=0$, $L_{\alpha\beta}=0$,\(\tilde{\nabla} H=0\), \(R_{\alpha N\beta N}=0\)), the terms \(-\tilde{\Delta} \frac{\partial}{\partial n} v\), \(-\frac{2}{3} H \tilde{\Delta} v\), \(L_{\alpha\beta} \tilde{\nabla}^\alpha \tilde{\nabla}^\beta v\), and \(\frac{1}{3} \tilde{\nabla} H \cdot \tilde{\nabla} v\) all vanish when applied to \(v\).

On each slice $\Sigma_s=\{s\}\times S^3$, let $T_+[\bar g]$
denote the boundary curvature computed with respect to
$N_+=e^{-v(s)}\partial_s$.
The outward unit normals of $[0,b]\times S^3$ are $-N_+$
at $s=0$ and $N_+$ at $s=b$.
Consequently, the two boundary contributions have opposite signs.

Let $P_{3,+}[g_0]$ denote the boundary operator computed with respect to the $g_0$-unit normal $\partial_s$. By the formula (\ref{DefOfP_3}) for $P_3$, we have  \(P_3[g_{0}] v = -\frac{1}{2} \frac{\partial}{\partial n} \Delta v + \left(\frac{1}{3} R\right) \frac{\partial}{\partial n} v = -\frac{1}{2} \frac{\partial}{\partial n} \Delta v + 2 \frac{\partial}{\partial n} v\). For the cylinder metric, \(\Delta v = v''\) and \(\frac{\partial}{\partial n} = \partial_t\), so \(\frac{\partial}{\partial n} \Delta v = \partial_t (v'') = v'''\) and \(\frac{\partial}{\partial n} v = v'(t)\).

Thus, \(P_{3,+}[g_0] v = -\frac{1}{2} v''' + 2 v'\). Hence it follows from (\ref{BoundaryOpt}) that
    \begin{align}
        T_+[\bar g]e^{3v}=P_{3,+}[g_0]v=-\frac{1}{2}v'''+2v'.\label{Newline:781A}
    \end{align}    

    Using the Chern--Gauss--Bonnet formula, we obtain
    \begin{align}
        \frac{1}{4\pi^2}\int_0^t\int_{S^3} Q[\bar g]e^{4v(s)}\,d\sigma\,ds=&-\frac{1}{4\pi^2}\int_{S^3}T_+[\bar g]e^{3v(t)}\,d\sigma\nonumber\\~~~~&+\frac{1}{4\pi^2}\int_{S^3}T_+[\bar g]e^{3v(0)}\,d\sigma.\label{Newline:777A}
    \end{align}
Substituting (\ref{Newline:781A}) into (\ref{Newline:777A}), we obtain 
    \begin{align}
        \frac{1}{8\pi^2}\int_0^t\int_{S^3}2Q[\bar g]e^{4v(s)}\,d\sigma\,ds=\frac{1}{4}v'''(t)-v'(t)+\frac{1}{4\pi^2}\int_{S^3}T_+[\bar g]e^{3v(0)}\,d\sigma.\label{NewLine:785A}
    \end{align}
Using the property (\ref{ANewLine:773}) of $v'''$ and letting $t \to \infty$, we obtain 
    \begin{align}
        \lim_{t\to+\infty}v'(t)=\frac{1}{4\pi^2}\int_{S^3}T_+[\bar g]e^{3v(0)}\,d\sigma-\frac{1}{4\pi^2}\int_0^{+\infty}\int_{S^3}Q[\bar g]e^{4v(s)}\,d\sigma\,ds.\label{NewLine:790}
    \end{align}\par
We now express the curvature integrals associated with $\bar g$ in \eqref{NewLine:790} in terms of the original metric $g$. We aim to replace the integral with respect to the metric $|dx|^2=dr^2+r^2d\theta^2=e^{2t}(dt^2+d\theta^2)=e^{2t}g_0$ in (\ref{NewLine:790}) by an integral with respect to the original metric $g = e^{2u(x)}|dx|^2=e^{2u(r\theta)}(dr^2+r^2d\theta^2)=e^{2u(e^t \theta)+2t}(dt^2+d\theta^2):=e^{2U(t,\theta)}g_0$. Similarly $\bar g:=e^{2D(|x|)}|dx|^2=e^{2D(r)}(dr^2+r^2d\theta^2)=e^{2D(e^t)+2t}g_0=e^{2v(t)}g_0$. With the preparation above, note that \begin{align}v(t)=D(e^t)+t=\fint_{S^3}\left(u(e^t \theta)+t\right)d\sigma=\fint_{S^3}U(t,\theta)d\sigma.\label{20260915}\end{align} Applying (\ref{BoundaryOpt}) and (\ref{20260915}) to the metrics $\bar g$ and $g$ with the common background metric $g_0$, we obtain 
    \begin{align}
T_+[\bar g]e^{3v(t)}=P_{3,+}[g_0]v(t)=\fint_{S^3}P_{3,+}[g_0]U(t,\theta)d\sigma=\fint_{S^3}T_+[g](e^t\theta)e^{3U(t,\theta)}d\sigma,\label{20260914111}
    \end{align}
    where the integration over $S^3$ eliminates the tangential
Laplacian terms in $P_{3,+}[g_0]$ and leaves only the radial term.
Setting $t=0$ in (\ref{20260914111}), we have 
    \begin{align}
        \frac{1}{4\pi^2}\int_{\partial B}T_+[\bar g]e^{3D(|y|)}d\sigma_y=\frac{1}{4\pi^2}\int_{\partial B}T_+[g]e^{3u(y)}d\sigma_y. \label{NewLine:798}
    \end{align}
    Finally, substituting (\ref{line:562}) and (\ref{NewLine:798}) into (\ref{NewLine:790}) yields directly
    \begin{align}
          \lim_{t\to \infty}v'(t)=\frac{1}{4\pi^2}\int_{\partial B}T_+[g]e^{3u(y)}d\sigma_y-\frac{1}{4\pi^2}\int_{\mathbb{R}^4\backslash B}Q_{g}(y)e^{4u(y)}dy.\label{20260824}  
    \end{align}
    Differentiating $v(t)=D(e^t)+t$ with respect to $t$, we obtain
    \begin{align}
        v'(t)=e^tD'(e^t)+1=rD'(r)+1.\label{202608241}
    \end{align}
    Thus, from (\ref{20260824}) and (\ref{202608241}), it follows that 
    \begin{align}
        \lim_{r\to \infty}rD'(r)&=\lim_{t\to \infty}v'(t) -1\nonumber\\
        &=\frac{1}{4\pi^2}\int_{\partial B}T_+[g]e^{3u(y)}d\sigma_y-\frac{1}{4\pi^2}\int_{\mathbb{R}^4\backslash B}Q_{g}(y)e^{4u(y)}dy-1,
    \end{align}
     which completes the proof.  
\end{proof}

\newpage
\bibliographystyle{unsrt}
\bibliography{reference}

\end{document}